\documentclass[10pt,aps,prb, preprint,reqno]{amsart}
\usepackage{amsmath}
\usepackage{amssymb}
\usepackage{yfonts}
\usepackage{hyperref}
\usepackage{mathrsfs}
\usepackage {latexsym}
\usepackage{graphics}
\usepackage{txfonts}
\usepackage{enumitem}
\allowdisplaybreaks

\DeclareMathSymbol{:}{\mathord}{operators}{"3A}

\newtheorem{theorem}{Theorem}[section]
\newtheorem{remark}{Remark}[section]
\newtheorem{lemma}[theorem]{Lemma}

\newtheorem{proposition}[theorem]{Proposition}

\newtheorem{define}{Definition}[section]

\allowdisplaybreaks 

\begin{document}
\title[Non-uniqueness in law up to the Lions' exponent]{Remarks on the non-uniqueness in law of the Navier-Stokes equations up to the J.-L. Lions' exponent} 
 
\author{Kazuo Yamazaki}  
\address{Texas Tech University, Department of Mathematics and Statistics, Lubbock, TX, 79409-1042, U.S.A.; Phone: 806-834-6112; Fax: 806-742-1112; E-mail: (kyamazak@ttu.edu)}
\date{}
\maketitle

\begin{abstract} 
Lions (1959, Bull. Soc. Math. France, \textbf{87}, 245--273) introduced the Navier-Stokes equations with a viscous diffusion in the form of a fractional Laplacian; subsequently, he (1969, Dunod, Gauthiers-Villars, Paris) claimed the uniqueness of its solution when its exponent is not less than five quarters in case the spatial dimension is three. Following the work of Hofmanov$\acute{\mathrm{a}}$, Zhu and Zhu (2019, arXiv:1912.11841 [math.PR]), we prove the non-uniqueness in law for the three-dimensional stochastic Navier-Stokes equations with the viscous diffusion in the form of a fractional Laplacian with its exponent less than five quarters. 
\vspace{5mm}

\end{abstract}
\footnote{2010MSC : 35A02; 35Q30; 35R60}

\section{Introduction}\label{Introduction}

After the pioneering work on the Navier-Stokes (NS) equations by Leray \cite{L34}, Lions in \cite[Remark 8.1 on pg. 263]{L59} introduced a generalized NS (GNS) equations \eqref{GNS equations} via an addition of a fractional Laplacian ``$(-1)^{m} \epsilon \Delta^{m}$'' (see \cite[Equation (8.7) on pg. 263]{L59}) and remarkably in \cite[Remark 6.11 on pg. 96]{L69} already claimed the uniqueness of the solution precisely when ``$m \geq \frac{n+2}{4}$'' in case the spatial variable $x$ is an $n$-dimensional vector (see \cite[Equation (6.164) on pg. 97]{L69}). To this day, except the logarithmic extension by Tao in \cite{T09} (see also \cite{BMR14} for further logarithmic improvement), it remains unknown whether or not for any initial condition sufficiently smooth, there exists a unique smooth solution in case the exponent is smaller than $\frac{n+2}{4}$ for $n \geq 3$ (we refer to \cite{W04, W05} for such results under the constraint on the smallness of initial condition). 

A lot of effort has been devoted toward verifying the uniqueness of the weak solution, which is known to exist globally in time (e.g., \cite[Chapter 3, Section 3]{T77}), as a potential path toward a successful proof of the global well-posedness of the three-dimensional (3-d) NS equations. Here and hereafter, by a weak solution we refer to a vector field $u \in C_{t}L_{x}^{2}$ that is merely weakly divergence-free and satisfies the NS equations in the sense of distributions; a well-known Leray-Hopf solution, originally constructed by Leray \cite{L34} and Hopf \cite{H51}, is additionally required to satisfy the energy inequality. As a culmination of a recent series of breakthroughs, we now know that this strategy fails, at least in the class of weak solutions aforementioned. Let us next introduce the GNS equations of our main concern and  describe the ground-breaking results that led to such a surprising outcome. 

Throughout this paper, except when we refer to the works of others, we consider $x = (x_{1}, x_{2}, x_{3}) \in \mathbb{T}^{3} = [-\pi, \pi]^{3}$, and thus $(-\Delta)^{m}$ may be considered as a Fourier operator with a symbol $\lvert \xi \rvert^{2m}$ so that $\widehat{(-\Delta)^{m} f}(\xi) = \lvert \xi \rvert^{2m} \hat{f}(\xi)$. Let us denote the velocity and pressure fields respectively by $u = (u^{1}, u^{2}, u^{3})$ and $\pi$ which map from $\mathbb{R}_{+} \times \mathbb{T}^{3}$, and additionally the viscosity coefficient by $\nu \geq 0$, so that the GNS equations read 
\begin{equation}\label{GNS equations}
\partial_{t} u + \nu (-\Delta)^{m} u + \text{div} (u\otimes u) + \nabla \pi = 0, \hspace{3mm} \nabla\cdot u =0, \hspace{3mm} t > 0,   
\end{equation}
which reduces to the classical NS equations when $\nu > 0$ and $m = 1$ while the Euler equations when $\nu = 0$. This system \eqref{GNS equations} has scaling-invariance such that if $(u,\pi)$ is a solution, then so is $(u_{\lambda}, \pi_{\lambda})(t,x) \triangleq (\lambda^{2m -1} u, \lambda^{4m -2} \pi)(\lambda^{2m} t, \lambda x)$. Consequently, \eqref{GNS equations} is supercritical with respect to (w.r.t.) the $L^{2}_{x}$-norm, which represents energy, if $m < \frac{5}{4}$, implying that the global existence of a unique solution by Lions applied only to the critical and subcritical cases. 

In 1949 Onsager \cite{O49} conjectured that in any spatial dimension strictly larger than one, in case $\alpha > \frac{1}{3}$ every weak solution $u$ to the Euler equations with spatial regularity in H$\ddot{\mathrm{o}}$lder class with exponent $\alpha$ conserves energy, while in case $\alpha \leq \frac{1}{3}$ there exists a weak solution $u$ with the  spatial regularity in H$\ddot{\mathrm{o}}$lder class with exponent $\alpha$ that does not conserve the energy. The positive direction of the Onsager's conjecture when $\alpha > \frac{1}{3}$ was settled in 1994 partially by Eyink \cite{E94}, although the result therein required a norm that is slightly stronger than H$\ddot{\mathrm{o}}$lder continuity, and then completely by Constantin, E and Titi \cite{CET94} using Littlewood-Paley theory techniques. Concerning the negative direction of Onsager's conjecture, in 1993 Scheffer \cite{S93} constructed a weak solution for the two-dimensional (2-d) Euler equations with compact support in space and time and thus the energy was certainly not conserved; however, Scheffer's solution had the regularity of $u \in L^{2}(\mathbb{R} \times \mathbb{R}^{2})$ and thus it was rougher than H$\ddot{\mathrm{o}}$lder class of any non-negative exponent. Subsequently in 1997, Shnirelman \cite{S97} constructed another weak solution $u \in L^{2}(\mathbb{R} \times \mathbb{T}^{2})$ to the Euler equations that vanishes whenever the temporal variable satisfies $\lvert t \rvert > C$ for some $C > 0$. Importantly, the main idea of the Shnirelman's proof was to construct a sequence of solutions to not the Euler equations but the Euler equations with an external force that in the limit oscillates infinitely fast in space such that it becomes unnoticed in the sense of distributions and consequently the limiting solution indeed solves the Euler equations weakly. It is of interest to emphasize that the works by Scheffer \cite{S93} and Shnirelman \cite{S97} consisted of constructing one weak solution with a specific feature while the works we shall discuss next demonstrate the existence of infinitely many weak solutions in the spirit of Gromov's homotopy-principle (h-principle) (e.g., \cite[pg. 3]{G86}). 

In 1954 Nash \cite{N54} proved the isometric embedding theorem by introducing a notion of ``short'' isometric embedding, which is an embedding that fails to be isometric by the error measured by a symmetric positive-definite tensor, and constructing a sequence of short isometric embeddings such that the error vanishes in the limit and consequently the limiting embedding indeed becomes isometric. This isometric embedding theorem of Nash, along with that of Kuiper \cite{K55}, was considered by Gromov to be one of the primary examples of h-principle, for which he initiated the method of convex integration \cite[Part 2.4]{G86}. M$\ddot{\mathrm{u}}$ller and $\check{\mathrm{S}}$ver$\acute{\mathrm{a}}$k \cite{MS03} extended the convex integration of Gromov to Lipschitz mappings and that led to  the proof of the the existence of the weak solution $u \in L^{\infty} (\mathbb{R} \times \mathbb{R}^{n})$ and $\pi \in L^{\infty} (\mathbb{R} \times \mathbb{R}^{n})$, where $n \geq 2$, to the Euler equations with compact support in space and time by De Lellis and Sz$\acute{\mathrm{e}}$kelyhidi Jr. \cite{DS09} using both convex integration and Baire category arguments. For more discussions on connections with the h-principle of Gromov \cite{G86}, we refer to \cite{DS12}. Subsequently, utilizing Beltrami flows and considering Reynolds stress, De Lellis and Sz$\acute{\mathrm{e}}$kelyhidi Jr. \cite{DS13} proved that given any smooth function $e \hspace{0.2mm}: \hspace{0.5mm} [0,1] \mapsto \mathbb{R}_{+}$, there exists a weak solution $u \in C(\mathbb{R}_{+} \times \mathbb{T}^{3})$ and $\pi\in C(\mathbb{R}_{+} \times \mathbb{T}^{3})$ to the Euler equations such that $e(t) = \lVert u(t) \rVert_{L_{x}^{2}}^{2}$ for all $t \in [0,1]$; in particular, this implies not only the lack of energy conservation but even dissipation (see also \cite{DS10}). For a discussion comparing the short embedding of Nash \cite{N54} and the Reynolds stress, we refer to \cite{T19}. Building on this idea of iterative constructions of Euler-Reynolds system, Buckmaster, De Lellis, Isett, and Sz$\acute{\mathrm{e}}$kelyhidi Jr. \cite{BDIS15} proved the negative direction of Onsager's conjecture with exponent up to $\frac{1}{5} - \epsilon$ for any $\epsilon > 0$. Subsequently, with a new strategy utilizing Mikado flows and a certain gluing scheme, Isett \cite{I18} proved the negative direction of Onsager's conjecture for any H$\ddot{\mathrm{o}}$lder exponent less than $\frac{1}{3}$ in any spatial dimension strictly larger than two (see \cite[Theorem 1 and page 877]{I18}). Making use of intermittency, Buckmaster and Vicol \cite{BV19a} proved the non-uniqueness of weak solutions to the 3-d NS equations; specifically, they proved that there exists $\gamma > 0$ such that for any non-negative smooth function $e(t)$ on $[0,T]$ and any $\nu \in (0, 1]$, the existence of a weak solution $u \in C([0,T]; H^{\gamma}(\mathbb{T}^{3}))$ such that $e(t) = \lVert u(t) \rVert_{L_{x}^{2}}^{2}$ for all $t \in [0,T]$ holds. A consequence of this result from \cite{BV19a} is the failure of the strategy of proving the global existence of a unique solution via weak solutions. It should be emphasized that the non-uniqueness proved in \cite{BV19a} was for a weak solution $u \in C_{t}L_{x}^{2}$, and not the Leray-Hopf weak solution, to which an extension is a challenging open problem. Concerning the GNS equations \eqref{GNS equations}, Luo and Titi \cite[Theorem 1]{LT20} extended the result of \cite{BV19a} to \eqref{GNS equations} with $m \in [1, \frac{5}{4})$ by adapting the method of using intermittent Beltrami flows in \cite{BV19a}. Moreover, via the introduction of intermittent jets, Buckmaster, Colombo, and Vicol \cite[Theorem 1.5]{BCV18} proved non-uniqueness of weak solutions to the GNS equations with $m \in [1, \frac{5}{4})$ which have bounded energy, integrable vorticity, and are smooth outside a fractal set of singular times with Hausdorff dimension strictly less than one. 

Our main focus in this manuscript will be the GNS equations forced by random noise, so-called generalized stochastic NS (GSNS) equations: 
\begin{equation}\label{[Equation (1.1), HZZ19]}
du + \nu (-\Delta)^{m} u dt + \text{div}(u\otimes u) dt + \nabla \pi dt = G(u) dB, \hspace{3mm} \nabla\cdot u =0, \hspace{3mm} t > 0,  
\end{equation}
where two different cases of $G(u)dB$ will be considered following \cite{HZZ19}: additive noise in the form of a cylindrical Wiener process and linear multiplicative noise driven by $\mathbb{R}$-valued Wiener process (see \eqref{[Equation (1.4), HZZ19]} and \eqref{[Equation (1.7), HZZ19]}). As we did for \eqref{GNS equations}, we refer to \eqref{[Equation (1.1), HZZ19]} when $m = 1$ as the stochastic NS (SNS) equations. Let us also recall that the uniqueness in law holds if for any solutions $(u, B)$ and $(\tilde{u}, \tilde{B})$ with same initial distributions, which may be defined on different filtered probability space, the law of $u$ coincides with that of $\tilde{u}$, while path-wise uniqueness holds if for any solutions $(u, B)$ and $(\tilde{u}, B)$ with common initial conditions defined on same probability space, $u(t) = \tilde{u}(t)$ for all $t$ with probability one. We refer to e.g., \cite[Example 2.2]{C03} for an instance in which uniqueness in law holds but path-wise uniqueness fails. 

The study of the SNS equations may be traced back to the work of Bensoussan and Temam \cite{BT73}. In particular,  Flandoli and Gatarek \cite{FG95} established the global existence of a weak solution to the 3-d SNS equations via Galerkin approximation (see \cite{F08} for a survey of results on the 3-d SNS equations). Path-wise uniqueness of the 2-d SNS equations is well-known (e.g., \cite{CM10}) and they can be extended to the GSNS equations for any $n \geq 3$ as long as $m \geq \frac{n+2}{4}$ analogously to the deterministic results (see also \cite{BF17} in the case of Boussinesq system). While the path-wise uniqueness of the 3-d classical SNS equations seemed to be as difficult as the deterministic case, Cherny \cite[Theorem 3.2]{C03} proved that the uniqueness in law and the existence of a strong solution, which requires that the solution is adapted to the completed natural filtration of $B$, together imply path-wise uniqueness. Therefore, many works (e.g., \cite{DD03, FR08, GRZ09}) were devoted toward proving the uniqueness in law for the 3-d SNS equations, especially via martingale problem (see Definitions \ref{[Definition 3.1, HZZ19]}, \ref{[Definition 3.2, HZZ19]}, \ref{[Definition 5.1, HZZ19]} and \ref{[Definition 5.2, HZZ19]}) formulated by Stroock and Varadhan \cite{SV97} because the uniqueness of the solution to the martingale problem is equivalent to the uniqueness in law (e.g., \cite[Corollary 4.9 in Chapter 5]{KS91}). In this perspective, it is truly stunning that Hofmanov$\acute{\mathrm{a}}$, Zhu, and Zhu \cite{HZZ19} proved the non-uniqueness in law for the 3-d SNS equations. Let us point out that first, by a contrapositive of the celebrated Yamada-Watanabe theorem, the work of \cite{HZZ19} implies a lack of path-wise uniqueness for the 3-d SNS equations. Second, the important tool in \cite{HZZ19} was the convex integration approach following \cite{BV19a}, more precisely \cite{BV19b}; therefore, the failure of an approach to prove path-wise uniqueness via uniqueness in law for the SNS equations may be related to the failure of an approach to prove the global well-posedness of the NS equations via weak solutions. Third, to the author's understanding, the work of Hofmanov$\acute{\mathrm{a}}$, Zhu, and Zhu \cite{HZZ19} makes a significant contribution to the Open Problem 9 in \cite[Section 8]{DS12}, which suggested probabilistic approach to convex integration. Fourth, the solution to the SNS equations constructed in \cite{HZZ19} does not have the regularity of $L_{t}^{2}H_{x}^{1}$ and thus are not at the level of Leray-Hopf solutions, identically to the current status in the deterministic case. 

\section{Statement of main results}
As we explained through discussions of the works \cite{BMR14, L69, T09, W04, W05}, the general consensus is that higher the exponent of the fractional Laplacian in the diffusive term of \eqref{GNS equations} and \eqref{[Equation (1.1), HZZ19]}, the easier to prove the global existence of a unique solution. Vice versa, the general consensus is that the higher the exponent, the more difficult to prove the non-uniqueness of a weak solution. Let us hereafter consider $\nu = 1$ in \eqref{[Equation (1.1), HZZ19]} for simplicity. 

\subsection{The case of additive noise}
In the case of an additive noise, we consider 
\begin{equation}\label{[Equation (1.4), HZZ19]}
du + (-\Delta)^{m} u dt + \text{div} (u\otimes u) dt + \nabla \pi dt = dB, \hspace{3mm} \nabla\cdot u =0, \hspace{3mm} t > 0, 
\end{equation}
where $B$ is a $GG^{\ast}$-Wiener process on a probability space $(\Omega, \mathcal{F}, \textbf{P})$ and $G$ is a certain Hilbert-Schmidt operator to be described in more detail subsequently in \eqref{[Equation (1.7c), HZZ19]}-\eqref{[Equation (1.7d), HZZ19]}, and the asterisk denotes the adjoint operator. We denote by $(\mathcal{F}_{t})_{t\geq 0}$ the filtration generated by $B$.   
\begin{theorem}\label{[Theorem 1.1, HZZ19]} 
Suppose that $m \in (\frac{13}{20}, \frac{5}{4})$ and Tr$((-\Delta)^{\frac{5}{2} - m + 2 \sigma} GG^{\ast} ) < \infty$ for some $\sigma > 0$. Then given $T> 0, K > 1$ and $\iota \in (0,1)$, there exists $\gamma \in (0,1)$ and a $\textbf{P}$-almost surely (a.s.) strictly positive stopping time $\mathfrak{t}$ such that $\textbf{P} ( \{ \mathfrak{t} \geq T \}) > \iota$ and the following is additionally satisfied. There exists an $(\mathcal{F}_{t})_{t\geq 0}$-adapted process $u$ which is a weak solution to \eqref{[Equation (1.4), HZZ19]} starting from $u^{\text{in}}$ that is deterministic and satisfies 
\begin{equation}\label{[Equation (1.5), HZZ19]}
\text{esssup}_{\omega \in \Omega} \text{sup}_{s \in [0, \mathfrak{t}]} \lVert u(s,\omega) \rVert_{H_{x}^{\gamma}} < \infty, 
\end{equation}
and on the set $\{ \mathfrak{t} \geq T \}$, 
\begin{equation}\label{[Equation (1.6), HZZ19]}
\lVert u(T) \rVert_{L_{x}^{2}} > K \lVert u^{\text{in}} \rVert_{L_{x}^{2}} + K ( T \text{Tr} (GG^{\ast} ))^{\frac{1}{2}}.
\end{equation} 
\end{theorem} 
The precise value of $\gamma \in (0,1)$ is required to be quite small; it comes from $\gamma \in (0, \frac{\beta}{4+ \beta})$ upon deriving \eqref{[Equation (4.12a), HZZ19]} in the proof of Theorem \ref{[Theorem 1.1, HZZ19]} where $\beta > 0$ is taken to be very small (see Sub-subsection \ref{Subsubsection 4.3.1}). By applying Theorem \ref{[Theorem 1.1, HZZ19]}, the following main result in the case of an additive noise can be deduced. 
\begin{theorem}\label{[Theorem 1.2, HZZ19]} 
Suppose that $m \in (\frac{13}{20}, \frac{5}{4})$ and Tr$((-\Delta)^{\frac{5}{2} - m + 2 \sigma} GG^{\ast} ) < \infty$ for some $\sigma > 0$. Then non-uniqueness in law holds for \eqref{[Equation (1.4), HZZ19]} on $[0,\infty)$. Moreover, for all $T > 0$ fixed, non-uniqueness in law holds for \eqref{[Equation (1.4), HZZ19]} on $[0, T]$. 
\end{theorem} 
\begin{remark}
As we will describe in Proposition \ref{[Proposition 3.6, HZZ19]}, if we demand that the process $z$ defined in \eqref{[Equation (3.10), HZZ19]} has regularity of $L_{\omega}^{1} C_{T}^{\frac{1}{2} - \delta} H_{x}^{\frac{3+ \sigma}{2}}$ for any $\delta \in (0, \frac{1}{2})$ similarly to \cite[Proposition 3.6]{HZZ19}, then it seems that we need a stronger hypothesis of $Tr ((-\Delta)^{\max \{ \frac{3}{2} + 2 \sigma, \frac{5}{2} - m + 2 \sigma} GG^{\ast} ) < \infty$. The author is grateful to the referee who suggested that under the current hypothesis of $Tr ((-\Delta)^{\frac{5}{2} - m + 2 \sigma} GG^{\ast}) < \infty$, we can still deduce that $z \in L_{\omega}^{1} C_{T}^{\frac{2}{5} - \delta} H_{x}^{\frac{3+ \sigma}{2}}$ for any $\delta \in (0, \frac{2}{5})$ and that actually suffices for the proof of the convex integration to go  through upon an appropriate adjustment of parameters, as it turns out that indeed that is the case (cf. \eqref{new 15} and \eqref{new 16}).
\end{remark}

\subsection{The case of linear multiplicative noise}
In this case we consider specifically 
\begin{equation}\label{[Equation (1.7), HZZ19]}
du + (-\Delta)^{m} u dt + \text{div} (u\otimes u) dt + \nabla \pi dt = u dB, \hspace{3mm} \nabla\cdot u = 0, \hspace{3mm} t > 0, 
\end{equation} 
where $B$ is a $\mathbb{R}$-valued Wiener process on $(\Omega, \mathcal{F}, \textbf{P})$. 
\begin{theorem}\label{[Theorem 1.3, HZZ19]}
Suppose that $m \in (\frac{13}{20}, \frac{5}{4})$ and $B$ is a $\mathbb{R}$-valued Wiener process on $(\Omega, \mathcal{F}, \textbf{P})$. Given $T> 0, K > 1$ and $\iota \in (0, 1)$, there exists $\gamma \in (0,1)$ and a $\textbf{P}$-a.s. strictly positive stopping time $\mathfrak{t}$ such that $\textbf{P} (\{ \mathfrak{t} \geq T \}) > \iota$ and the following is additionally satisfied. There exists an $(\mathcal{F}_{t})_{t\geq 0}$-adapted process $u$ which is a weak solution to \eqref{[Equation (1.7), HZZ19]} starting from $u^{\text{in}}$ that is deterministic and satisfies 
\begin{equation}\label{[Equation (1.7a), HZZ19]}
\text{esssup}_{\omega \in \Omega} \text{sup}_{s\in [0, \mathfrak{t}]} \lVert u(s,\omega) \rVert_{H_{x}^{\gamma}} < \infty 
\end{equation} 
and on the set $\{\mathfrak{t} \geq T \}$, 
\begin{equation}\label{[Equation (1.7a'), HZZ19]}
\lVert u(T) \rVert_{L_{x}^{2}} > K e^{\frac{T}{2}} \lVert u^{\text{in}} \rVert_{L_{x}^{2}}. 
\end{equation} 
\end{theorem}  
Similarly to Theorem \ref{[Theorem 1.2, HZZ19]}, applying Theorem \ref{[Theorem 1.3, HZZ19]} deduces the following result: 
\begin{theorem}\label{[Theorem 1.4, HZZ19]}
Suppose that $m \in (\frac{13}{20}, \frac{5}{4})$ and $B$ is a $\mathbb{R}$-valued Wiener process on $(\Omega, \mathcal{F}, \textbf{P})$. Then non-uniqueness in law holds for \eqref{[Equation (1.7), HZZ19]} on $[0,\infty)$. Moreover, for any $T > 0$ fixed, non-uniqueness in law holds for \eqref{[Equation (1.7), HZZ19]} on $[0,T]$.  
\end{theorem}

\begin{remark}
Our proof can be immediately reduced to deterministic case, and we attain non-uniqueness even in case $m < 1$, specifically $m \in (\frac{13}{20}, 1)$ (cf. \cite[Theorem 1.1]{LT20} and \cite[Theorem 1.5]{BCV18} which worked on the case $m \in [1, \frac{5}{4})$). This lower bound of $\frac{13}{20} < m$ is technical and comes from \eqref{[Equation (4.23), HZZ19]}-\eqref{[Equation (B.2a), HZZ19]}. \footnote{The current manuscript was completed in June 2020. In April 2021, the author in \cite{Y21} was able to prove non-uniqueness in law for \eqref{[Equation (1.7), HZZ19]} for $m \in (0, \frac{1}{2})$ by a completely different type of convex integration scheme.} Differently from the proof within \cite{LT20} that relied on intermittent Beltrami flows, we follow the method of intermittent jets within \cite{BCV18, BV19b}. On the other hand, the application of convex integration in \cite{BCV18} is very different from ours because their main result is the non-uniqueness of weak solutions that is almost everywhere smooth and emanates from an initial condition in $C([0,T]; \dot{H}^{3}(\mathbb{T}^{3}))$; it will be an interesting future problem to obtain a stochastic analog of \cite[Theorem 1.1, Remark 1.2]{BCV18}. Strictly speaking of an application of convex integration, we will follow \cite[Chapter 7]{BV19b} as did the authors in \cite{HZZ19}; the difference then is that \cite[Chapter 7]{BV19b} focuses on the NS equations, not the GNS equations. Thus, we adapt the method in \cite{BCV18} on the GNS equations and intermittent jets to the stochastic setting of \cite{HZZ19} and the convex integration scheme in \cite[Chapter 7]{BV19b} which requires finding various appropriate parameters (e.g., see \eqref{alpha}, \eqref{[Equation (4.23), HZZ19]}, \eqref{p}), many of which depend on the value of $m \in (\frac{13}{20}, \frac{5}{4})$ in complex manners. Some estimates must be carefully conducted as well, in comparison to the case $m = 1$ (e.g., Remark \ref{star}); indeed, the value ``$5- 4m$'' appears multiples of times throughout our proofs clearly displaying the threshold of $\frac{5}{4}$. We also mention that it may be possible to further extend Theorems \ref{[Theorem 1.1, HZZ19]}-\ref{[Theorem 1.4, HZZ19]} to the Lagrangian-averaged Navier-Stokes-$\alpha$ model that has caught much attention (e.g., \cite{CMR05}) with filtration in terms of fractional Laplacians (e.g., \cite{BMR14, Y12}). 
\end{remark} 

\section{Preliminaries}\label{Preliminaries}
In this subsection we set up minimum amount of notations and assumptions, while leave the rest to the Appendix for completeness. As we will see in Definitions \ref{[Definition 3.1, HZZ19]}, \ref{[Definition 3.2, HZZ19]}, \ref{[Definition 5.1, HZZ19]} and \ref{[Definition 5.2, HZZ19]}, we will consider martingale problems, although we choose to follow \cite{HZZ19} and not use this terminology, and thus many of the following notations are similar to those in \cite{FR08, GRZ09} (also \cite{Y19}). For brevity we write $A \lesssim_{a,b} B$ and $A \approx_{a,b}B$ to imply the existence of a constant $C = C(a,b) \geq 0$ such that $A \leq CB$ and $A = CB$, respectively. When helpful, we may also write $A \overset{(\cdot)}{\lesssim} B$ e.g. to indicate that this inequality is due to an equation $(\cdot)$. We denote $\mathbb{N} \triangleq \{1, 2, ...\}$ while $\mathbb{N}_{0} \triangleq \mathbb{N} \cup \{0\}$. Following \cite[pg. 105]{BV19a} we write for $p\in[1,\infty]$, 
\begin{equation}\label{estimate 48}
\lVert f \rVert_{L^{p}} \triangleq \lVert f \rVert_{L_{t}^{\infty} L_{x}^{p}}, \hspace{1mm} \lVert f \rVert_{C^{N}} \triangleq \lVert f \rVert_{L_{t}^{\infty} C_{x}^{N}} \triangleq \sum_{0\leq \lvert \alpha \rvert \leq N} \lVert D^{\alpha} f \rVert_{L^{\infty}}, \hspace{1mm} \lVert f \rVert_{C_{t,x}^{N}} \triangleq \sum_{0\leq n + \lvert \alpha \rvert \leq N} \lVert \partial_{t}^{n} D^{\alpha} f \rVert_{L^{\infty}}. 
\end{equation} 
We also denote $L_{\sigma}^{2} \triangleq \{f \in L_{x}^{2}(\mathbb{T}^{3}): \hspace{0.5mm} \nabla\cdot f =0 \}$. We reserve $\mathbb{P}$ for the Leray projection operator and denote by $\mathring{\otimes}$ the trace-free part of a tensor product. We also denote $\mathbb{P}_{<r}$ as a Fourier operator with a Fourier symbol $1_{\{ \lvert \xi \rvert < r \}}(\xi)$. For any Polish space $H$ we write $\mathcal{B}(H)$ to denote the $\sigma$-algebra of Borel sets in $H$. Given any probability measure $P$, we denote a mathematical expectation w.r.t. $P$ by $\mathbb{E}^{P}$. We denote by $\langle \cdot, \cdot \rangle$ the $L^{2}(\mathbb{T}^{3})$-inner product while $\langle\langle A, B \rangle \rangle$ a cross variation of $A$ and $B$, and $\langle \langle A \rangle \rangle \triangleq \langle \langle A, A \rangle \rangle$. We let $\Omega_{0} \triangleq C([0,\infty); H^{-3} (\mathbb{T}^{3})) \cap L_{\text{loc}}^{\infty} ([0,\infty); L_{\sigma}^{2})$. We also denote by $\mathcal{P}(\Omega_{0})$ the set of all probability measures on $(\Omega_{0}, \mathcal{B})$ where $\mathcal{B}$ is the Borel $\sigma$-field of $\Omega_{0}$ from the topology of locally uniform convergence on $\Omega_{0}$. We define $\xi: \hspace{0.5mm}\Omega_{0} \mapsto H^{-3}(\mathbb{T}^{3})$ the canonical process by 
\begin{equation}\label{[Equation (1.7b), HZZ19]}
\xi_{t}(\omega) \triangleq \omega(t).
\end{equation} 
Similarly, for $t \geq 0$ we define $\Omega_{t} \triangleq C([t,\infty); H^{-3}(\mathbb{T}^{3})) \cap L_{\text{loc}}^{\infty} ((t, \infty); L_{\sigma}^{2})$ equipped with Borel $\sigma$-algebra $\mathcal{B}^{t} \triangleq \sigma (\{\xi(s): \hspace{0.5mm} s \geq t \})$. Furthermore, we define $\mathcal{B}_{t}^{0} \triangleq \sigma (\{ \xi(s): \hspace{0.5mm} s \leq t \})$  and $\mathcal{B}_{t} \triangleq \cap_{s > t} \mathcal{B}_{s}^{0}$ for $t \geq 0$. For any Hilbert space $U$, we denote by $L_{2}(U, L_{\sigma}^{2})$ the space of all Hilbert-Schmidt operators from $U$ to $L_{\sigma}^{2}$ with the norm $\lVert \cdot \rVert_{L_{2}(U, L_{\sigma}^{2})}$. We require $G: \hspace{0.5mm} L_{\sigma}^{2} \mapsto L_{2}(U, L_{\sigma}^{2})$ to be $\mathcal{B}(L_{\sigma}^{2}) / \mathcal{B}(L_{2}(U, L_{\sigma}^{2}))$-measurable, that it satisfies 
\begin{equation}\label{[Equation (1.7c), HZZ19]}
\lVert G(\phi) \rVert_{L_{2}(U, L_{\sigma}^{2})} \leq C(1+ \lVert \phi \rVert_{L_{x}^{2}})
\end{equation} 
for any $\phi \in C^{\infty} (\mathbb{T}^{3}) \cap L_{\sigma}^{2}$ and 
\begin{equation}\label{[Equation (1.7d), HZZ19]}
\lim_{n\to\infty} \lVert G(\theta_{n})^{\ast} \phi - G(\theta)^{\ast} \phi \rVert_{U} = 0 
\end{equation} 
if $\lim_{n\to\infty} \lVert \theta_{n} - \theta \rVert_{L_{x}^{2}} = 0$. 

The following additional notations are for the case of a linear multiplicative noise. We assume the existence of another Hilbert space $U_{1}$ such that the embedding $U \subset U_{1}$ is Hilbert-Schmidt. Then we define $\bar{\Omega} \triangleq C([0,\infty); H^{-3}(\mathbb{T}^{3}) \times U_{1}) \cap L_{\text{loc}}^{\infty} ([0,\infty); L_{\sigma}^{2} \times U_{1})$ and $\mathcal{P} (\bar{\Omega})$ to be the set of all probability measures on $(\bar{\Omega}, \bar{\mathcal{B}})$ where $\bar{\mathcal{B}}$ is the Borel $\sigma$-algebra on $\bar{\Omega}$. Furthermore, we define the canonical process on $\bar{\Omega}$ by $(\xi, \theta): \hspace{0.5mm}  \bar{\Omega} \mapsto H^{-3}(\mathbb{T}^{3}) \times U_{1}$ to satisfy $(\xi_{t}(\omega), \theta_{t}(\omega)) \triangleq \omega(t)$. We extend the previous definitions of $\mathcal{B}^{t}, \mathcal{B}_{t}^{0}$ and $\mathcal{B}_{t}$ so that $\bar{\mathcal{B}}^{t} \triangleq \sigma (\{ (\xi(s), \theta(s)): \hspace{0.5mm}  s \geq t \})$, $\bar{\mathcal{B}}_{t}^{0} \triangleq \sigma ( \{ (\xi(s), \theta(s)): \hspace{0.5mm}  s \leq t \})$ and $\bar{\mathcal{B}}_{t} \triangleq \cap_{s> t} \bar{\mathcal{B}}_{s}^{0}$ for $t \geq 0$. 

Next we provide proof of our main results. We aim to elaborate only where there are significant differences from the proofs within \cite{HZZ19}, while including sufficient details to keep this manuscript relatively self-contained. 

\section{Proof in the case of additive noise}

\subsection{Proof of Theorem \ref{[Theorem 1.2, HZZ19]} assuming Theorem \ref{[Theorem 1.1, HZZ19]} }

Let us fix $\gamma \in (0,1)$ for the following general definition. 
\begin{define}\label{[Definition 3.1, HZZ19]}
Let $s \geq 0$ and $\xi^{\text{in}} \in L_{\sigma}^{2}$. Then $P \in \mathcal{P}(\Omega_{0})$ is a martingale solution to \eqref{[Equation (1.1), HZZ19]} with initial condition $\xi^{\text{in}}$ at initial time $s$ if 
\begin{itemize}
\item [] (M1) $P( \{ \xi(t) = \xi^{\text{in}} \hspace{1mm} \forall \hspace{1mm} t \in [0,s] \}) = 1$ and for all $n \in \mathbb{N}$ 
\begin{equation}\label{[Equation (1.7f), HZZ19]}
P( \{ \xi \in \Omega_{0}: \int_{0}^{n} \lVert G(\xi(r)) \rVert_{L_{2}(U, L_{\sigma}^{2})}^{2} dr < \infty \}) = 1, 
\end{equation} 
\item [] (M2) for every $e_{i} \in C^{\infty} (\mathbb{T}^{3} ) \cap L_{\sigma}^{2}$ and $t \geq s$, the process 
\begin{equation}\label{[Equation (1.7g), HZZ19]}
M_{t,s}^{i} \triangleq \langle \xi(t) - \xi(s), e_{i} \rangle + \int_{s}^{t} \langle \text{div} (\xi(r) \otimes \xi(r)) + (-\Delta)^{m} \xi(r), e_{i} \rangle dr
\end{equation}
is a continuous, square-integrable martingale w.r.t. $(\mathcal{B}_{t})_{t\geq s}$ under $P$ with $\langle \langle M_{t,s}^{i} \rangle \rangle = \int_{s}^{t} \lVert G(\xi(r))^{\ast} e_{i} \rVert_{U}^{2} dr$,  
\item [] (M3) for any $q \in \mathbb{N}$, there exists a function $t \mapsto C_{t,q} \in \mathbb{R}_{+}$ for all $t \geq s$ such that  
\begin{equation}\label{[Equation (1.7h), HZZ19]}
\mathbb{E}^{p} [ \sup_{r \in [0,t]} \lVert \xi(r) \rVert_{L_{x}^{2}}^{2q} + \int_{s}^{t} \lVert \xi(r) \rVert_{H_{x}^{\gamma}}^{2} dr ] \leq C_{t,q} (1+ \lVert \xi^{\text{in}} \rVert_{L_{x}^{2}}^{2q}). 
\end{equation}
\end{itemize}
The set of all such martingale solutions with the same constant $C_{t,q}$ in \eqref{[Equation (1.7h), HZZ19]} for every $q \in \mathbb{N}$ and $t \geq s$ will be denoted by $\mathcal{C} (s, \xi^{\text{in}}, \{C_{t,q}\}_{q\in\mathbb{N}, t \geq s})$. 
\end{define}
In case the noise is additive so that $G$ is independent of $\xi$, if $\{e_{i}\}_{i=1}^{\infty}$ is a complete orthonormal system that consists of eigenvectors of $GG^{\ast}$, then $M_{t,s} \triangleq \sum_{i=1}^{\infty} M_{t,s}^{i} e_{i}$  becomes a $GG^{\ast}$-Wiener process starting from initial time $s$ w.r.t. the filtration $(\mathcal{B}_{t})_{t\geq s}$ under $P$. In order to define a martingale solution up to a stopping time $\tau: \hspace{0.5mm}  \Omega_{0} \mapsto [0,\infty]$, we define the space of trajectories stopped at time $\tau$ by 
\begin{equation}\label{[Equation (1.7i), HZZ19]}
\Omega_{0,\tau} \triangleq \{\omega( \cdot  \wedge \tau(\omega)) : \hspace{0.5mm}  \omega \in \Omega_{0} \}
\end{equation} 
and denote by $(\mathcal{B}_{\tau})$ the $\sigma$-field associated to $\tau$. 
\begin{define}\label{[Definition 3.2, HZZ19]}
Let $s \geq 0$, $\xi^{\text{in}} \in L_{\sigma}^{2}$ and $\tau \geq s$ be a stopping time of $(\mathcal{B}_{t})_{t\geq s}$. Then $P \in \mathcal{P} (\Omega_{0,\tau})$ is a martingale solution to \eqref{[Equation (1.1), HZZ19]} on $[s, \tau]$ with initial condition $\xi^{\text{in}}$ at initial time $s$ if 
\begin{enumerate}
\item [] (M1) $P(\{ \xi(t) = \xi^{\text{in}} \hspace{1mm} \forall \hspace{1mm} t \in [0, s]\}) = 1$ and for all $n \in \mathbb{N}$ 
\begin{equation}\label{[Equation (1.7j), HZZ19]}
P(\{ \xi \in \Omega_{0}: \int_{0}^{n \wedge \tau} \lVert G(\xi(r)) \rVert_{L_{2}(U, L_{\sigma}^{2})}^{2} dr < \infty \}) = 1, 
\end{equation}
\item [] (M2) for every $e_{i} \in C^{\infty} (\mathbb{T}^{3}) \cap L_{\sigma}^{2}$ and $t\geq s$, the process 
\begin{equation}\label{[Equation (1.7k), HZZ19]}
M_{t \wedge \tau, s}^{i} \triangleq \langle \xi(t\wedge \tau) - \xi^{\text{in}}, e_{i} \rangle + \int_{s}^{t \wedge \tau} \langle \text{div} (\xi(r) \otimes \xi(r)) + (-\Delta)^{m} \xi(r), e_{i} \rangle dr 
\end{equation}
is a continuous, square-integrable martingale w.r.t. $(\mathcal{B}_{t})_{t\geq s}$ under $P$ with $\langle \langle M_{t\wedge \tau, s}^{i} \rangle \rangle = \int_{s}^{t\wedge \tau} \lVert G(\xi(r))^{\ast} e_{i} \rVert_{U}^{2} dr$, 
\item [] (M3) for any $q \in \mathbb{N}$, there exists a function $t \mapsto C_{t,q} \in \mathbb{R}_{+}$ for all $t \geq s$ such that 
\begin{equation}\label{[Equation (1.7l), HZZ19]}
\mathbb{E}^{P} [ \sup_{r \in [0, t \wedge \tau]} \lVert \xi(r) \rVert_{L_{x}^{2}}^{2q} + \int_{s}^{t \wedge \tau} \lVert \xi(r) \rVert_{H_{x}^{\gamma}}^{2} dr] \leq C_{t,q} (1+ \lVert \xi^{\text{in}} \rVert_{L_{x}^{2}}^{2q}).
\end{equation}
\end{enumerate}
\end{define}

One can trace the proof of \cite[Theorem 3.1]{HZZ19} (also \cite[Theorem 4.1]{FR08}, \cite[Theorem 3.1]{Y19}) to deduce the following result concerning the existence and certain stability of a martingale solution to \eqref{[Equation (1.1), HZZ19]}. 
\begin{proposition}\label{[Theorem 3.1, HZZ19]}
For any $(s, \xi^{\text{in}}) \in [0,\infty) \times L_{\sigma}^{2}$, there exists $P \in \mathcal{P}(\Omega_{0})$ which is a martingale solution to \eqref{[Equation (1.1), HZZ19]} with initial condition $\xi^{\text{in}}$ at initial time $s$  that satisfies Definition \ref{[Definition 3.1, HZZ19]}. Moreover, if there exists a family $\{(s_{n}, \xi_{n})\}_{n\in\mathbb{N}} \subset [0,\infty) \times L_{\sigma}^{2}$ such that $\lim_{n\to\infty} \lVert (s_{n}, \xi_{n}) - (s, \xi^{\text{in}}) \rVert_{\mathbb{R} \times L_{x}^{2}} = 0 $ and $P_{n} \in \mathcal{C} (s_{n}, \xi_{n}, \{C_{t,q}\}_{q\in\mathbb{N}, t \geq s_{n}})$ is the martingale solution corresponding to $(s_{n}, \xi_{n})$, then there exists a subsequence $\{P_{n_{k}}\}_{k\in\mathbb{N}}$ that converges weakly to some $P \in \mathcal{C}(s, \xi^{\text{in}}, \{C_{t,q}\}_{q \in \mathbb{N}, t \geq s})$. 
\end{proposition} 

\begin{proof}[Proof of Proposition \ref{[Theorem 3.1, HZZ19]}]
We sketch its proof in the Appendix only highlighting the difference from the proof of \cite[Theorem 3.1]{HZZ19} due to the diffusive term with a fractional Laplacian in our case for completeness. 
\end{proof}

The proofs of the following results from \cite{HZZ19} do not rely on the specific form of the diffusive term and thus apply directly to our case. 
\begin{lemma}\label{[Proposition 3.2, HZZ19]} 
\rm{(\cite[Proposition 3.2]{HZZ19})} Let $\tau$ be a bounded stopping time of $(\mathcal{B}_{t})_{t\geq 0}$. Then for every $\omega \in \Omega_{0}$, there exists $Q_{\omega} \in \mathcal{P}(\Omega_{0})$ such that 
\begin{subequations}
\begin{align}
&Q_{\omega} (\{ \omega' \in \Omega_{0}: \hspace{0.5mm}  \xi(t, \omega') = \omega(t) \hspace{1mm} \forall \hspace{1mm} t \in [0, \tau(\omega)] \}) = 1, \label{[Equation (3.1), HZZ19]}\\
&Q_{\omega} (A) = R_{\tau(\omega), \xi(\tau(\omega), \omega)} (A) \hspace{1mm} \forall \hspace{1mm} A \in \mathcal{B}^{\tau(\omega)}, \label{[Equation (3.2), HZZ19]}
\end{align}
\end{subequations}
where $R_{\tau(\omega), \xi(\tau(\omega), \omega)} \in \mathcal{P}(\Omega_{0})$ is a martingale solution to \eqref{[Equation (1.1), HZZ19]} with initial condition $\xi(\tau(\omega), \omega)$ at initial time $\tau(\omega)$. Furthermore, for every $B \in \mathcal{B},$ the mapping $\omega \mapsto Q_{\omega}(B)$ is $\mathcal{B}_{\tau}$-measurable.
\end{lemma} 
Let us only mention that in the proof of Lemma \ref{[Proposition 3.2, HZZ19]}, $Q_{\omega}$ is derived from \cite[Lemma 6.1.1]{SV97} as the unique probability measure 
\begin{equation}\label{[Equation (3.4), HZZ19]}
Q_{\omega} = \delta_{\omega} \otimes_{\tau(\omega)} R_{\tau(\omega), \xi(\tau(\omega), \omega)} \in \mathcal{P} (\Omega_{0}),
\end{equation}
where $\delta_{\omega}$ is the Dirac mass, such that \eqref{[Equation (3.1), HZZ19]}-\eqref{[Equation (3.2), HZZ19]} hold. 

\begin{lemma}\label{[Proposition 3.4, HZZ19]}
\rm{(\cite[Proposition 3.4]{HZZ19})} Let $\xi^{\text{in}} \in L_{\sigma}^{2}$ and $P$ be a martingale solution to \eqref{[Equation (1.1), HZZ19]} on $[0,\tau]$ with initial condition $\xi^{\text{in}}$ at initial time $0$ that satisfies Definition \ref{[Definition 3.2, HZZ19]}. Assume the hypothesis of Lemma \ref{[Proposition 3.2, HZZ19]} and additionally that there exists a Borel set $\mathcal{D} \subset \Omega_{0, \tau}$ such that $P(\mathcal{D}) = 0$ and for every $\omega \in \Omega_{0} \setminus \mathcal{D}$, it satisfies 
\begin{equation}\label{[Equation (3.5), HZZ19]}
Q_{\omega} (\{\omega' \in \Omega_{0}: \hspace{0.5mm}  \tau(\omega') = \tau(\omega) \}) = 1. 
\end{equation} 
Then the probability measure $P \otimes_{\tau} R \in \mathcal{P} (\Omega_{0})$ defined by 
\begin{equation}\label{[Equation (3.6), HZZ19]}
P \otimes_{\tau} R(\cdot) \triangleq \int_{\Omega_{0}} Q_{\omega} (\cdot) P(d\omega) 
\end{equation} 
satisfies $P \otimes_{\tau} R \rvert_{\Omega_{0,\tau}} = P \rvert_{\Omega_{0, \tau}}$ and it is a martingale solution to \eqref{[Equation (1.1), HZZ19]} on $[0,\infty)$ with initial condition $\xi^{\text{in}}$ at initial time 0. 
\end{lemma} 

Now we let $\mathcal{B}_{\tau}$ represent the $\sigma$-field associated to the stopping time $\tau$. We split \eqref{[Equation (1.4), HZZ19]} to 
\begin{equation}\label{[Equation (3.10), HZZ19]}
dz + (-\Delta)^{m} z dt + \nabla \pi^{1} dt = dB, \hspace{3mm} \nabla\cdot z = 0 \text{ for } t > 0, \hspace{1mm} z(0,x) = 0,
\end{equation}
\begin{equation}\label{[Equation (3.11), HZZ19]}
\partial_{t} v + (-\Delta)^{m} v + \text{div} ((v+z) \otimes (v+z)) + \nabla \pi^{2} = 0, \hspace{1mm} \nabla\cdot v =0 \text{ for } t > 0, \hspace{1mm} v(0,x) = u^{\text{in}}(x),
\end{equation}
so that $u = v + z$ solves \eqref{[Equation (1.4), HZZ19]} with $\pi = \pi^{1} + \pi^{2}$. We fix a $GG^{\ast}$-Wiener process $B$ on $(\Omega, \mathcal{F}, \textbf{P})$ with $(\mathcal{F}_{t})_{t\geq 0}$ as the canonical filtration of $B$ augmented by all the $\textbf{P}$-negligible sets. We see that 
\begin{equation}\label{new 11}
z(t) = \int_{0}^{t} e^{- (-\Delta)^{m}(t-r)} \mathbb{P} dB_{r} 
\end{equation}
where $e^{(-\Delta)^{m} t}$ denotes a semigroup generated by $(-\Delta)^{m}$.  Concerning regularity of $z$, we have the following result (cf. \cite[Proposition 3.6]{HZZ19}). 
\begin{proposition}\label{[Proposition 3.6, HZZ19]}
Suppose that Tr$((-\Delta)^{\frac{5}{2} - m + 2 \sigma} GG^{\ast} ) < \infty$ for some $\sigma > 0$. Then for all $\delta \in (0, \frac{2}{5})$ and $T > 0$, 
\begin{equation}\label{[Equation (3.11a), HZZ19]}
\mathbb{E}^{\textbf{P}}[ \lVert z \rVert_{C_{T} H_{x}^{\frac{5+\sigma}{2}}} + \lVert z \rVert_{C_{T}^{\frac{2}{5} - \delta } H_{x}^{\frac{3+ \sigma}{2}}}] < \infty.  
\end{equation}
\end{proposition} 
\begin{proof}[Proof of Proposition \ref{[Proposition 3.6, HZZ19]}]
For completeness, let us give details explaining how the hypothesis is used following the proof of \cite[Proposition 34]{D13} that uses factorization method (cf. similar derivation to \cite[Equation (25)]{Y19}). For generality, let us write our hypothesis as 
\begin{equation}\label{new 1}
Tr ((-\Delta)^{m (1+ g)} GG^{\ast}) < \infty. 
\end{equation} 
For $\alpha \in (0,1)$ to be subsequently specified, we define 
\begin{equation}\label{new 2}
Y(s) \triangleq \frac{\sin(\alpha \pi)}{\pi} \int_{0}^{s} e^{-(-\Delta)^{m} (s-r)} (s-r)^{-\alpha} \mathbb{P} dB(r). 
\end{equation} 
It follows from \eqref{new 11} and the identity 
\begin{equation*}
\int_{r}^{t} (t-s)^{\alpha -1} (s-r)^{-\alpha} ds  \frac{\sin(\alpha \pi)}{\pi} = 1 
\end{equation*} 
from \cite[p. 131]{DZ14} that 
\begin{equation}\label{new 3}
z(t) = \int_{0}^{t} (t-s)^{\alpha -1} e^{- (-\Delta)^{m} (t-s)} Y(s) ds. 
\end{equation} 
For $\epsilon$ and $g$ to be specified shortly and $l \in \mathbb{N}$, we compute by Gaussian hypercontractivity theorem (e.g., \cite[Theorem 3.50]{J97}) and It$\hat{\mathrm{o}}$ isometry (e.g., \cite[Equation (4)]{D13}) that if 
\begin{equation}\label{new 4} 
\alpha < \min\{ \frac{g}{2} - \epsilon, \frac{1}{2} \}, 
\end{equation} 
then 
\begin{align}
\mathbb{E}^{\textbf{P}} [ \lVert (-\Delta)^{m(1+ \epsilon)} Y(s) \rVert_{L_{x}^{2}}^{2l} ] \overset{\eqref{new 2}}{\lesssim_{l}}& ( \mathbb{E}^{\textbf{P}} [ \int_{0}^{s} \lVert (-\Delta)^{m (1+ \epsilon)} e^{- (-\Delta)^{m} (s-r)} (s-r)^{-\alpha} G \rVert_{\mathcal{L}_{2}}^{2} dr])^{l} \nonumber \\
\overset{\eqref{new 1}}{\lesssim_{l}}& ( \mathbb{E}^{\textbf{P}} [ \int_{0}^{s} (s-r)^{-2\alpha} ((s-r)^{-(1+ 2 \epsilon - g)} + 1) dr])^{l}  \overset{\eqref{new 4}}{\lesssim_{l}} 1. \label{new 5}  
\end{align}
Integrating over $[0,T]$ and using Fubini's theorem now give us for any $l \in \mathbb{N}$
\begin{equation}\label{new 6} 
\mathbb{E}^{\textbf{P}} [ \int_{0}^{T} \lVert (-\Delta)^{m(1+ \epsilon)} Y(s) \rVert_{L_{x}^{2}}^{2l}ds ] \overset{\eqref{new 5}}{\lesssim_{l}} 1. 
\end{equation} 
By taking $l > \frac{1}{2\alpha}$ we can now deduce 
\begin{equation}\label{new 7}
\mathbb{E}^{\textbf{P}} [ \lVert (-\Delta)^{m(1+ \epsilon)} z \rVert_{C_{T} L_{x}^{2}}^{2l} ] 
\overset{\eqref{new 3}\eqref{new 6}}{\lesssim_{l}} 1.
\end{equation} 
Therefore, in order to deduce our first claim that $\mathbb{E}^{\textbf{P}}[ \lVert z \rVert_{C_{T}H_{x}^{\frac{5+ \sigma}{2}}}] < \infty$, according to \eqref{new 4} and \eqref{new 7} we need $\alpha, \epsilon,$ and $g$ such that 
\begin{equation}\label{new 8} 
m(1+ \epsilon) \geq \frac{5+\sigma}{4} \hspace{2mm} \text{ and } \hspace{2mm} 0 < \alpha < \min\{ \frac{g}{2} - \epsilon, \frac{1}{2}\}. 
\end{equation} 
Thus, let us fix $\epsilon = \frac{5+ \sigma}{4m} - 1$ so that the first inequality of \eqref{new 8} is satisfied and require that $g > 2 \epsilon = \frac{5+\sigma}{2m} -2$. This allows us to choose $g = \frac{5+ 4\sigma}{2m} - 2$ so that 
\begin{equation}\label{new 13} 
\infty \overset{\eqref{new 1}}{>} Tr ( ( -\Delta)^{m(1+ g)} GG^{\ast}) = Tr ((-\Delta)^{m(\frac{5+ 4 \sigma}{2m} -1)} GG^{\ast}) = Tr((-\Delta)^{\frac{5}{2} - m + 2 \sigma} GG^{\ast}) 
\end{equation} 
is satisfied by hypothesis.  Next, let us deduce the second claim $\mathbb{E}^{\textbf{P}} [ \lVert z \rVert_{C_{T}^{\frac{2}{5} - \delta} H_{x}^{\frac{3+ \sigma}{2}}}] < \infty$ for any $\delta \in (0, \frac{2}{5})$. Due to \cite[Proposition A.1.1]{DZ92} and \eqref{new 7} we have for any $\beta \in (0, \alpha - \frac{1}{2l})$ 
\begin{equation}\label{new 9}
\mathbb{E}^{\textbf{P}} [ \sup_{t, t + h \in [0,T]} \lVert (-\Delta)^{m(1+ \epsilon)} (z(t+h) - z(t)) \rVert_{L_{x}^{2}}^{2l} ] \lesssim_{\epsilon, \beta, l, T} \lvert h \rvert^{2\beta l}. 
\end{equation} 
By Kolmogorov's test (e.g., \cite[Theorem 3.3]{DZ14}), \eqref{new 9} implies that $\textbf{P}$-almost every trajectories of $\lVert z (t) \rVert_{H_{x}^{2m(1+ \epsilon)}} \in C_{T}^{\mu}$ for $\mu < \beta - \frac{1}{2l}$. Therefore, in order to deduce $\mathbb{E}^{\textbf{P}} [ \lVert z \rVert_{C_{T}^{\frac{2}{5} - \delta} H_{x}^{\frac{3+ \sigma}{2}}}] < \infty$ for any $\delta \in (0, \frac{2}{5})$, we need $\mu < \beta - \frac{1}{2l}$ arbitrarily close to $\frac{2}{5}$ and thus $\beta < \alpha - \frac{1}{2l}$ arbitrarily close to $\frac{2}{5}$, which in turn implies $\frac{g}{2} - \epsilon \geq \frac{2}{5}$ considering \eqref{new 4} . Hence, we must choose $\epsilon$, and $g$ such that 
\begin{equation}\label{new 10}
m(1+ \epsilon) \geq \frac{3+\sigma}{4}, \hspace{2mm} 0 < \alpha < \min\{\frac{g}{2} - \epsilon, \frac{1}{2}\} = \frac{2}{5}, \hspace{2mm} \text{ and } \hspace{2mm} \frac{g}{2} - \epsilon \geq \frac{2}{5},  
\end{equation} 
in contrast to \eqref{new 8}. For optimality we choose $\epsilon = \frac{3+ \sigma}{4m} - 1$ and realize that a choice of $g = \frac{5+ 4 \sigma}{2m} - 2$ satisfies $\frac{g}{2} - \epsilon = \frac{2+ 3 \sigma}{4m} \geq \frac{2}{5}$ because $m < \frac{5}{4}$. At last, this choice of $g$ leads us again to \eqref{new 13} that is satisfied by hypothesis. 

We point out that repeating this argument to attain $\mathbb{E}^{\textbf{P}}[ \lVert z \rVert_{C_{T}^{\frac{1}{2} - \delta} H_{x}^{\frac{3+ \sigma}{2}}}] < \infty$ instead of $\mathbb{E}^{\textbf{P}}[ \lVert z \rVert_{C_{T}^{\frac{2}{5} - 2 \delta} H_{x}^{\frac{3+ \sigma}{2}}}] < \infty$, which boils down to 
\begin{equation}\label{new 14}
m(1+ \epsilon) \geq \frac{3+\sigma}{4}, \hspace{2mm} 0 < \alpha < \min\{\frac{g}{2} - \epsilon, \frac{1}{2}\} = \frac{1}{2}, \hspace{2mm}  \text{ and } \hspace{2mm} \frac{g}{2} - \epsilon \geq \frac{1}{2}, 
\end{equation} 
shows the necessity of a hypothesis of $\infty > Tr (( - \Delta)^{\frac{3}{2} + 2 \sigma} GG^{\ast})$ instead.
\end{proof}

Next for every $\omega \in \Omega_{0}$ we define 
\begin{equation}\label{[Equation (3.12), HZZ19]}
M_{t,0}^{\omega} \triangleq \omega(t) - \omega(0) + \int_{0}^{t} \mathbb{P} \text{div} (\omega(r) \otimes \omega(r)) + (-\Delta)^{m} \omega(r) dr,
\end{equation} 
\begin{equation}\label{[Equation (3.13), HZZ19]}
Z^{\omega} (t) \triangleq M_{t,0}^{\omega} - \int_{0}^{t} \mathbb{P} (-\Delta)^{m} e^{-(t-r) (-\Delta)^{m}} M_{r,0}^{\omega} dr. 
\end{equation} 
If $P$ is a martingale solution to \eqref{[Equation (1.4), HZZ19]}, then the process $M$ is a $GG^{\ast}$-Wiener process under $P$ and it follows from \eqref{[Equation (3.12), HZZ19]}-\eqref{[Equation (3.13), HZZ19]} that 
\begin{equation}\label{[Equation (3.13a), HZZ19]}
Z(t) =  \int_{0}^{t} \mathbb{P} e^{-(t-r) (-\Delta)^{m}} dM_{r,0}. 
\end{equation} 
It follows from Proposition \ref{[Proposition 3.6, HZZ19]} that $Z \in C_{T} H_{x}^{\frac{5+\sigma}{2}} \cap C_{T}^{\frac{2}{5} - \delta} H_{x}^{\frac{3+ \sigma}{2}}$ $P$-almost surely for any $\delta \in (0, \frac{2}{5})$. For $n \in \mathbb{N}$ and 
\begin{equation}\label{new 15}
\delta \in (0, \frac{1}{30}), 
\end{equation} 
we define 
\begin{align}
\tau_{L}^{n} (\omega) \triangleq& \inf \{t \geq 0: \hspace{0.5mm}  \lVert Z^{\omega}(t) \rVert_{H_{x}^{\frac{5+\sigma}{2}}} > \frac{ (L - \frac{1}{n})^{\frac{1}{4}}}{C_{S}} \} \nonumber\\
& \wedge \inf \{ t \geq 0: \hspace{0.5mm} \lVert Z^{\omega} \rVert_{C_{t}^{\frac{2}{5} - 2 \delta} H_{x}^{\frac{3 + \sigma}{2}}} > \frac{ (L - \frac{1}{n})^{\frac{1}{2}}}{C_{S}} \} \wedge L  \label{[Equation (3.13b), HZZ19]}
\end{align} 
where $C_{S} > 0$ is the Sobolev constant such that $\lVert f \rVert_{L_{x}^{\infty}} \leq C_{s} \lVert f \rVert_{H_{x}^{\frac{3 + \sigma}{2}}}$. This choice of $\delta < \frac{1}{30}$ comes from \eqref{[Equation (4.48o), HZZ19]}-\eqref{[Equation (4.48p), HZZ19]}, specifically to attain $\frac{2}{5} - 2 \delta > \frac{1}{3}$ therein. Clearly $(\tau_{L}^{n})_{n\in\mathbb{N}}$ is non-decreasing in $n$ and we define 
\begin{equation}\label{[Equation (3.14), HZZ19]}
\tau_{L} \triangleq \lim_{n\to\infty} \tau_{L}^{n}. 
\end{equation}
It follows from \cite[Lemma 3.5]{HZZ19} that $\tau_{L}^{n}$ is a stopping time of $(\mathcal{B}_{t})_{t \geq 0}$ and therefore so is $\tau_{L}$. Next we assume Theorem \ref{[Theorem 1.1, HZZ19]} on $(\Omega, \mathcal{F}, (\mathcal{F}_{t})_{t\geq 0}, \textbf{P})$ and denote by $P$ the law of the solution constructed by Theorem \ref{[Theorem 1.1, HZZ19]} to deduce the following result (cf. \cite[Proposition 3.7]{HZZ19}). 
\begin{proposition}\label{[Proposition 3.7, HZZ19]}
Let $\tau_{L}$ be defined by \eqref{[Equation (3.14), HZZ19]}. Then $P$, the law of $u$, is a martingale solution of \eqref{[Equation (1.4), HZZ19]} on $[0, \tau_{L}]$ that satisfies Definition \ref{[Definition 3.2, HZZ19]}. 
\end{proposition} 
\begin{proof}[Proof of Proposition \ref{[Proposition 3.7, HZZ19]}]
We define for $C_{S} > 0$ from \eqref{[Equation (3.13b), HZZ19]}, $L > 1$ and $\delta \in (0, \frac{1}{30})$, 
\begin{align}\label{[Equation (4.2), HZZ19]}
T_{L} \triangleq \inf \{t \geq 0: \hspace{0.5mm}  \lVert z(t) \rVert_{H_{x}^{\frac{5+ \sigma}{2}}} \geq \frac{L^{\frac{1}{4}}}{C_{S}} \}  \wedge \inf \{t \geq 0: \hspace{0.5mm} \lVert z \rVert_{C_{t}^{\frac{2}{5} - 2 \delta} H_{x}^{\frac{3+ \sigma}{2}}} \geq \frac{L^{\frac{1}{2}}}{C_{S}} \} \wedge L.
\end{align} 
We observe that due to Proposition \ref{[Proposition 3.6, HZZ19]}, $T_{L}$ is $\textbf{P}$-a.s. strictly positive and $T_{L} \nearrow + \infty$ as $L \nearrow + \infty$ $\textbf{P}$-a.s.; moreover, the stopping time $\mathfrak{t}$ in the statement of Theorem \ref{[Theorem 1.1, HZZ19]} is actually $T_{L}$ for $L > 0$ sufficiently large, as we will see from the proof of Theorem \ref{[Theorem 1.1, HZZ19]}. The proof of Proposition \ref{[Proposition 3.7, HZZ19]} is similar to that of \cite[Proposition 3.7]{HZZ19}. The proof therein does not rely on the specific form of the diffusive term; however, our definitions of $Z^{\omega}, z$ are slightly different. For completeness we sketch the proof in the Appendix. 
\end{proof}

Next we extend $P$ on $[0,\tau_{L}]$ to $[0,\infty)$ similarly to \cite[Proposition 3.8]{HZZ19}. 
\begin{proposition}\label{[Proposition 3.8, HZZ19]}
The probability measure $P\otimes_{\tau_{L}} R$ in \eqref{[Equation (3.6), HZZ19]} with $\tau_{L}$ from \eqref{[Equation (3.14), HZZ19]} is a martingale solution to \eqref{[Equation (1.4), HZZ19]} on $[0,\infty)$ that satisfies Definition \ref{[Definition 3.1, HZZ19]}. 
\end{proposition}

\begin{proof}[Proof of Proposition \ref{[Proposition 3.8, HZZ19]}]
Due to \eqref{[Equation (3.13b), HZZ19]}, $\tau_{L}$ is a bounded stopping time of $(\mathcal{B}_{t})_{t\geq 0}$ and thus the hypothesis of Lemma \ref{[Proposition 3.2, HZZ19]} is satisfied. Thus, by Proposition \ref{[Proposition 3.7, HZZ19]}, Lemma \ref{[Proposition 3.4, HZZ19]} gives us the desired result once we verify the existence of a Borel set $\mathcal{D} \subset \Omega_{0,\tau}$ such that $P(\mathcal{D}) = 0$ and \eqref{[Equation (3.5), HZZ19]} holds for every $\omega \in \Omega_{0,\tau} \setminus \mathcal{D}$. This verification requires only some modifications of the proof of \cite[Proposition 3.8]{HZZ19} considering that our definition of $Z^{\omega}, z$ are slightly different; we provide a sketch in the Appendix for completeness. 
\end{proof}

\begin{proof}[Proof of Theorem \ref{[Theorem 1.2, HZZ19]} assuming Theorem \ref{[Theorem 1.1, HZZ19]}]
The proof of \cite[Theorem 1.2]{HZZ19} essentially applies to our case as the specific form of the diffusive term plays no role therein. In short, we can fix $T > 0$ arbitrarily, any $\iota \in (0,1)$ and $K > 1$ such that $\iota K^{2} \geq 1$, rely on Theorem \ref{[Theorem 1.1, HZZ19]} and Proposition \ref{[Proposition 3.8, HZZ19]} to deduce the existence of $L > 1$ and a measure $P \otimes_{\tau_{L}} R$ that is a martingale solution to \eqref{[Equation (1.4), HZZ19]} on $[0,\infty)$ and coincides with $P$, the law of the solution constructed in Theorem \ref{[Theorem 1.1, HZZ19]}, over a random interval $[0, \tau_{L} ]$. Therefore, $P \otimes_{\tau_{L}} R$ starts  with a deterministic initial condition $\xi^{\text{in}}$ from the proof of Theorem \ref{[Theorem 1.1, HZZ19]}. It follows that 
\begin{align}
P \otimes_{\tau_{L}} R( \{ \tau_{L} \geq T \}) \overset{\eqref{[Equation (3.6), HZZ19]}\eqref{[Equation (3.19), HZZ19]}  }{=}& P \otimes_{\tau_{L}} R ( \{ \tau_{L} (\omega) \geq T \}) \overset{\eqref{[Equation (3.15), HZZ19]}}{=} \textbf{P} ( \{ T_{L} \geq T \})  > \iota   \label{[Equation (3.19a), HZZ19]}
\end{align} 
so that 
\begin{align}
\mathbb{E}^{P \otimes_{\tau_{L}} R} [ \lVert \xi(T) \rVert_{L_{x}^{2}}^{2}] 
\overset{\eqref{[Equation (1.6), HZZ19]} \eqref{[Equation (3.19a), HZZ19]} }{>} \iota [ K \lVert \xi^{\text{in}} \rVert_{L_{x}^{2}} + K (T \text{Tr} (GG^{\ast} ))^{\frac{1}{2}}]^{2}   
\geq  \iota K^{2}(\lVert \xi^{\text{in}} \rVert_{L_{x}^{2}}^{2} + T \text{Tr} (GG^{\ast})). 
\end{align} 
On the other hand, the classical method of Galerkin approximation gives us another martingale solution $\Theta$ (e.g., \cite{FR08, Y19}) which starts from the same initial condition $\xi^{\text{in}}$ and satisfies 
\begin{equation*}
\mathbb{E}^{\Theta} [ \lVert \xi(T)\rVert_{L_{x}^{2}}^{2}] \leq \lVert \xi^{\text{in}} \rVert_{L_{x}^{2}}^{2} + T\text{Tr}(GG^{\ast}).
\end{equation*} 
Because $\iota K^{2} \geq 1$, this implies $P \otimes_{\tau_{L}} R\neq \Theta$ and hence \eqref{[Equation (1.4), HZZ19]} fails the uniqueness in law. 
\end{proof}

\subsection{Proof of Theorem \ref{[Theorem 1.1, HZZ19]} assuming Proposition \ref{[Proposition 4.2, HZZ19]}}
Considering \eqref{[Equation (3.11), HZZ19]} we will construct solutions $(v_{q}, \mathring{R}_{q})$ for $q \in \mathbb{N}_{0}$ to 
\begin{equation}\label{[Equation (4.1), HZZ19]}
\partial_{t} v_{q} + (-\Delta)^{m} v_{q} + \text{div} ((v_{q} + z) \otimes ((v_{q} + z)) + \nabla \pi_{q} = \text{div} \mathring{R}_{q}, \hspace{3mm} \nabla\cdot v_{q} = 0, \hspace{3mm} t > 0. 
\end{equation} 
We set for $a > 0, b \in \mathbb{N}$ and $\beta \in (0,1)$, 
\begin{equation}\label{[Equation (4.1a), HZZ19]}
\lambda_{q} \triangleq a^{b^{q}} \text{ and } \delta_{q} \triangleq \lambda_{q}^{-2\beta}. 
\end{equation} 
We see that by \eqref{[Equation (4.2), HZZ19]}, for any $\delta \in (0, \frac{1}{30})$ and $t \in [0,T_{L}]$, 
\begin{equation}\label{[Equation (4.3), HZZ19]}
\lVert z(t) \rVert_{L_{x}^{\infty}} \leq L^{\frac{1}{4}}, \hspace{1mm} \lVert \nabla z(t) \rVert_{L_{x}^{\infty}} \leq L^{\frac{1}{4}}, \hspace{1mm} \lVert z \rVert_{C_{t}^{\frac{2}{5} - 2 \delta} L_{x}^{\infty}} \leq L^{\frac{1}{2}}. 
\end{equation} 
Now we make the assumption of $a^{\beta b} > 3$, which will be formally stated in \eqref{[Equation (4.7), HZZ19]} of Proposition \ref{[Lemma 4.1, HZZ19]},  so that because $b^{r-1} \geq r$ for all $b \geq 2$ and $r \geq 1$, one can deduce $\sum_{1 \leq r \leq q} \delta_{r}^{\frac{1}{2}} < \frac{1}{2}$ for any $q \in \mathbb{N}$. We furthermore define 
\begin{equation}\label{estimate 4}
M_{0}(t) \triangleq L^{4} e^{4Lt}, 
\end{equation} 
set the convention of $\sum_{1 \leq r \leq 0} \triangleq 0$, let $c_{R} > 0$ be a universal constant that is sufficiently small with specific conditions subsequently given (e.g.,  \eqref{[Equation (4.28), HZZ19]}, \eqref{[Equation (4.29), HZZ19]}, \eqref{[Equation (4.37), HZZ19]}), and assume the following bounds over $t \in [0, T_{L}]$ inductively: 
\begin{subequations}\label{[Equation (4.4), HZZ19]}
\begin{align}
& \lVert v_{q} \rVert_{C_{t}L_{x}^{2}} \leq M_{0}(t)^{\frac{1}{2}} (1+ \sum_{1 \leq r \leq q} \delta_{r}^{\frac{1}{2}}) \leq 2 M_{0}(t)^{\frac{1}{2}}, \label{[Equation (4.4a), HZZ19]}\\
& \lVert v_{q} \rVert_{C_{t,x}^{1}} \leq M_{0}(t)^{\frac{1}{2}} \lambda_{q}^{4}, \label{[Equation (4.4b), HZZ19]}\\
& \lVert \mathring{R}_{q} \rVert_{C_{t}L_{x}^{1}} \leq M_{0}(t) c_{R} \delta_{q+1}. \label{[Equation (4.4c), HZZ19]}
\end{align}
\end{subequations}
Our choice of $\mathring{R}_{0}$ in the following result differs from that of \cite[Lemma 4.1]{HZZ19} (cf. \cite[Equation (7.4)]{BV19b}). 
\begin{proposition}\label{[Lemma 4.1, HZZ19]}
Define an operator $\mathcal{R}$ by Lemma \ref{divergence inverse operator}. For $L > 1$ let
\begin{equation}\label{[Equation (4.5a), HZZ19]}
v_{0} (t,x) \triangleq \frac{L^{2} e^{2Lt}}{(2\pi)^{\frac{3}{2}}} 
\begin{pmatrix}
\sin(x_{3}) & 0 & 0
\end{pmatrix}^{T},
\end{equation} 
where $T$ denotes a transpose. Then together with 
\begin{align}
\mathring{R}_{0}(t,x) =& 
\frac{ 2 L^{3} e^{2L t}}{(2\pi)^{\frac{3}{2}}} 
\begin{pmatrix}
0 & 0 & - \cos(x_{3})\\
0 & 0 & 0 \\
-\cos(x_{3}) & 0 & 0 
\end{pmatrix} 
+ \mathcal{R} (-\Delta)^{m} v_{0} + v_{0} \mathring{\otimes} z + z \mathring{\otimes} v_{0} + z \mathring{\otimes} z, \label{[Equation (4.6), HZZ19]}
\end{align} 
it satisfies \eqref{[Equation (4.1), HZZ19]} at level $q= 0$. Moreover, \eqref{[Equation (4.4a), HZZ19]}-\eqref{[Equation (4.4c), HZZ19]} are satisfied at level $q = 0$ provided 
\begin{equation}\label{[Equation (4.7), HZZ19]}
16 \leq L \text{ and } (17) (2\pi)^{\frac{3}{2}} (9) < (17) (2\pi)^{\frac{3}{2}} a^{2 \beta b} \leq c_{R} L \leq c_{R} \left( \frac{ (2\pi)^{\frac{3}{2}} a^{4} - 2}{2} \right), 
\end{equation}
where the inequality of $9 < a^{2\beta b}$ in \eqref{[Equation (4.7), HZZ19]} is assumed for the sake of second inequality in \eqref{[Equation (4.4a), HZZ19]}. Furthermore, $v_{0}(0,x)$ and $\mathring{R}_{0}(0,x)$ are both deterministic. 
\end{proposition} 

\begin{proof}[Proof of Proposition \ref{[Lemma 4.1, HZZ19]}]
Using that $\text{div} \mathcal{R} ((-\Delta)^{m} v_{0}) = (-\Delta)^{m} v_{0}$ by Lemma \ref{divergence inverse operator} and that divergence of a matrix $(a_{ij})_{1 \leq i, j \leq 3}$ is a 3-$d$ vector of which $k$-th component is $\sum_{j=1}^{3} \partial_{j} a_{kj}$, it is immediately verified that $(v_{0}, \mathring{R}_{0})$ indeed satisfies \eqref{[Equation (4.1), HZZ19]} at level $q= 0$ if we choose $\pi_{0} = -\frac{1}{3} (2 v_{0} \cdot z + \lvert z \rvert^{2})$. Next we can immediately deduce 
\begin{subequations}
\begin{align}
& \lVert v_{0}(t) \rVert_{L_{x}^{2}} = \frac{ L^{2} e^{2Lt}}{\sqrt{2}} \overset{\eqref{estimate 4}}{\leq} M_{0}(t)^{\frac{1}{2}}, \label{[Equation (4.8a), HZZ19]}\\
& \lVert v_{0} \rVert_{C_{t,x}^{1}} = \frac{L^{2} e^{2L t} (2L +1)}{(2\pi)^{\frac{3}{2}}} \overset{\eqref{[Equation (4.7), HZZ19]}}{\leq} M_{0}(t)^{\frac{1}{2}} \lambda_{0}^{4}, \label{[Equation (4.8b), HZZ19]}
\end{align} 
\end{subequations}
and thus \eqref{[Equation (4.4a), HZZ19]} -  \eqref{[Equation (4.4b), HZZ19]} are both satisfied at level $q= 0$. Next we see that 
\begin{equation*}
\lVert \mathring{R}_{0}(t) \rVert_{L_{x}^{1}} \leq 16 L^{3} e^{2L t} (2\pi)^{\frac{1}{2}} + (2\pi)^{\frac{3}{2}} \lVert \mathcal{R} (-\Delta)^{m} v_{0} \rVert_{L_{x}^{2}} + 20 M_{0}(t)^{\frac{1}{2}} L^{\frac{1}{4}} (2\pi)^{\frac{3}{2}} + 10L^{\frac{1}{2}} (2\pi)^{3} 
 \end{equation*} 
due to \eqref{[Equation (4.8a), HZZ19]} and \eqref{[Equation (4.3), HZZ19]}. We use the fact that $\Delta v_{0} = -v_{0}$; moreover, it is immediate to estimate directly using \eqref{estimate 5} that $\lVert \mathcal{R} (-\Delta)^{m} f \rVert_{L_{x}^{2}} \leq 24 \lVert \Delta f \rVert_{L_{x}^{2}}$. Hence, 
\begin{align*}
\lVert \mathring{R}_{0}(t) \rVert_{L_{x}^{1}} \leq& 16 L M_{0}(t)^{\frac{1}{2}} (2\pi)^{\frac{1}{2}} + (2\pi)^{\frac{3}{2}} 24 M_{0}(t)^{\frac{1}{2}} \\
&+ 20M_{0}(t)^{\frac{1}{2}} L^{\frac{1}{4}} (2\pi)^{\frac{3}{2}} + 10 L^{\frac{1}{2}} (2\pi)^{3} \overset{\eqref{estimate 4} \eqref{[Equation (4.7), HZZ19]}}{\leq} 17 (2\pi)^{\frac{3}{2}} M_{0}(t) L^{-1}.
\end{align*} 
Therefore, due to \eqref{[Equation (4.7), HZZ19]} we are able to deduce \eqref{[Equation (4.4c), HZZ19]} at level $q = 0$. Finally, it is clear from \eqref{[Equation (4.5a), HZZ19]} that $v_{0}(0,x)$ is deterministic and consequently, it follows from \eqref{[Equation (4.6), HZZ19]} and that $z(0,x) = 0$ from \eqref{[Equation (3.10), HZZ19]} that $\mathring{R}_{0}(0,x)$ is also deterministic. 
\end{proof}

Before we state the next result, let us point out that $L > 16 \vee (c_{R}^{-1}153 (2\pi)^{\frac{3}{2}} )$ is not sufficient but necessary to guarantee \eqref{[Equation (4.7), HZZ19]}. 
\begin{proposition}\label{[Proposition 4.2, HZZ19]}
Let $L > 16 \vee (c_{R}^{-1}153 (2\pi)^{\frac{3}{2}} )$ and suppose that $(v_{q}, \mathring{R}_{q})$ is an $(\mathcal{F}_{t})_{t\geq 0}$-adapted process that solves \eqref{[Equation (4.1), HZZ19]} and satisfies \eqref{[Equation (4.4), HZZ19]}.  Then there exists  a choice of parameters $a, b, \beta$ such that \eqref{[Equation (4.7), HZZ19]} is fulfilled and there exists an $(\mathcal{F}_{t})_{t\geq 0}$-adapted process $(v_{q+1}, \mathring{R}_{q+1})$ that solves \eqref{[Equation (4.1), HZZ19]}, satisfies \eqref{[Equation (4.4), HZZ19]} at level $q+1$ and 
\begin{equation}\label{[Equation (4.12), HZZ19]}
\lVert v_{q+1}(t) - v_{q}(t) \rVert_{L_{x}^{2}} \leq M_{0}(t)^{\frac{1}{2}} \delta_{q+1}^{\frac{1}{2}}.
\end{equation} 
Finally, if $v_{q}(0,x)$ and $\mathring{R}_{q}(0,x)$ are deterministic, then so are $v_{q+1}(0,x)$ and $\mathring{R}_{q+1}(0,x)$. 
\end{proposition}

We assume Proposition \ref{[Proposition 4.2, HZZ19]} and prove Theorem \ref{[Theorem 1.1, HZZ19]} now. 
\begin{proof}[Proof of Theorem \ref{[Theorem 1.1, HZZ19]}  assuming Proposition \ref{[Proposition 4.2, HZZ19]}]

This proof is similar to that of \cite[Theorem 1.1]{HZZ19}; we sketch it for completeness. Given any $T > 0, K > 1$ and $\iota \in (0,1)$, starting from $(v_{0}, \mathring{R}_{0})$ in Proposition \ref{[Lemma 4.1, HZZ19]}, Proposition \ref{[Proposition 4.2, HZZ19]} gives us $(v_{q}, \mathring{R}_{q})$ for $q \geq 1$ that satisfies \eqref{[Equation (4.4), HZZ19]} and \eqref{[Equation (4.12), HZZ19]}. Then, for all $\gamma \in (0, \frac{\beta}{4+ \beta})$ and $t \in [0, T_{L}]$, by Gagliardo-Nirenberg's inequality and the fact that $b^{q+1} \geq b(q+1)$ for all $q \geq 0$ and $b \geq 2$, we can obtain
\begin{align}
\sum_{q \geq 0} \lVert v_{q+1}(t) - v_{q}(t) \rVert_{H_{x}^{\gamma}} \overset{\eqref{[Equation (4.4), HZZ19]}  \eqref{[Equation (4.12), HZZ19]}   }{\lesssim} M_{0}(t)^{\frac{1}{2}} \sum_{q\geq 0} \delta_{q+1}^{\frac{1-\gamma}{2}} \lambda_{q+1}^{4\gamma} \lesssim M_{0}(t)^{\frac{1}{2}}. \label{[Equation (4.12a), HZZ19]}
\end{align} 
Thus, we can deduce the limiting solution $\lim_{q\to\infty} v_{q} \triangleq v \in C([0, T_{L}]; H^{\gamma}(\mathbb{T}^{3}))$ and an existence of a deterministic constant $C_{L}> 0$ such that 
\begin{equation}\label{[Equation (4.13), HZZ19]}
\sup_{t \in [0, T_{L}]} \lVert v(t) \rVert_{H_{x}^{\gamma}} \leq C_{L}. 
\end{equation} 
Since each $v_{q}$ is $(\mathcal{F}_{t})_{t\geq 0}$-adapted, $v$ is also $(\mathcal{F}_{t})_{t \geq 0}$-adapted. Moreover, as $\lim_{q\to \infty} \lVert \mathring{R}_{q} \rVert_{C_{T_{L}}L_{x}^{1}} = 0$ due to \eqref{[Equation (4.4c), HZZ19]}, we see that $v$ is a weak solution to \eqref{[Equation (3.11), HZZ19]} over $[0, T_{L}]$, and consequently from \eqref{[Equation (3.10), HZZ19]}-\eqref{[Equation (3.11), HZZ19]} we see that $u$ solves \eqref{[Equation (1.4), HZZ19]}. Now for $c_{R} > 0$ from the proof of Proposition \ref{[Proposition 4.2, HZZ19]}, we can choose $L = L(T, K, c_{R}, \text{Tr}(GG^{\ast})) > c_{R}^{-1} 153 (2\pi)^{\frac{3}{2}}$ such that 
\begin{equation}\label{[Equation (4.15), HZZ19]} 
(\frac{3}{2} + \frac{1}{L}) < (\frac{1}{\sqrt{2}} - \frac{1}{2}) e^{LT} \text{ and } L^{\frac{1}{4}} (2\pi)^{\frac{3}{2}} + K (T \text{Tr} (GG^{\ast}))^{\frac{1}{2}} \leq (e^{LT} - K) \lVert u^{\text{in}} \rVert_{L_{x}^{2}}+ L e^{LT},
\end{equation} 
where $u^{\text{in}}(x) = v(0,x)$ by \eqref{[Equation (3.10), HZZ19]}. For $T > 0$ and $\iota > 0$ fixed, increasing $L$ larger if necessary allows us to deduce $\textbf{P} (\{T_{L} \geq T \}) > \iota$. From \eqref{[Equation (3.10), HZZ19]} we see that  $z$ is $(\mathcal{F}_{t})_{t \geq 0}$-adapted, and thus $u$ is $(\mathcal{F}_{t})_{t \geq 0}$-adapted. Moreover, \eqref{[Equation (4.13), HZZ19]} and \eqref{[Equation (4.3), HZZ19]} imply \eqref{[Equation (1.5), HZZ19]}. Next, as $b^{q+1} \geq b(q+1)$ for all $q \in \mathbb{N}_{0}$ if $b \geq 2$, we can compute for all $t \in [0, T_{L}]$, 
\begin{align}
\lVert v(t) - v_{0}(t) \rVert_{L_{x}^{2}} \leq \sum_{q \geq 0} \lVert v_{q+1}(t) - v_{q}(t) \rVert_{L_{x}^{2}} 
\overset{\eqref{[Equation (4.1a), HZZ19]}\eqref{[Equation (4.12), HZZ19]} }{\leq}  M_{0}(t)^{\frac{1}{2}} \sum_{q\geq 0} \lambda_{q+1}^{-\beta}\overset{\eqref{[Equation (4.7), HZZ19]}}{<}& M_{0}(t)^{\frac{1}{2}}(\frac{1}{2}). \label{[Equation (4.14a), HZZ19]} 
\end{align}
Due to \eqref{[Equation (4.15), HZZ19]} this gives us
\begin{align}\label{[Equation (4.14), HZZ19]}
( \lVert v(0) \rVert_{L_{x}^{2}} + L)e^{LT} \overset{\eqref{[Equation (4.8a), HZZ19]}\eqref{[Equation (4.14a), HZZ19]}}{\leq} ( \frac{3}{2} M_{0} (0)^{\frac{1}{2}} +L) e^{LT} \overset{\eqref{[Equation (4.8a), HZZ19]} \eqref{[Equation (4.15), HZZ19]} \eqref{[Equation (4.14a), HZZ19]}}{<}& \lVert v(T) \rVert_{L_{x}^{2}}. 
\end{align} 
Therefore, on $\{T_{L} \geq T \}$, 
\begin{align}\label{[Equation (4.15a), HZZ19]}
\lVert u(T) \rVert_{L_{x}^{2}} \overset{\eqref{[Equation (4.14), HZZ19]}}{>} (\lVert v(0) \rVert_{L_{x}^{2}} + L)e^{LT} - \lVert z(T) \rVert_{L_{x}^{\infty}} (2\pi)^{\frac{3}{2}} \overset{ \eqref{[Equation (4.3), HZZ19]} \eqref{[Equation (4.15), HZZ19]}}{\geq} K \lVert u^{\text{in}} \rVert_{L_{x}^{2}} + K (T \text{Tr} (GG^{\ast}))^{\frac{1}{2}}.   
\end{align} 
This implies \eqref{[Equation (1.6), HZZ19]}. Finally, because $v_{0}(0,x)$ is deterministic by Proposition \ref{[Lemma 4.1, HZZ19]}, Proposition \ref{[Proposition 4.2, HZZ19]} implies that $v(0,x)$ is deterministic; by \eqref{[Equation (3.10), HZZ19]} this implies that $u^{\text{in}}$ is deterministic. 
\end{proof}

\subsection{Proof of Proposition \ref{[Proposition 4.2, HZZ19]}}
We now proceed with the proof of Proposition \ref{[Proposition 4.2, HZZ19]}. We refer to the Appendix for preliminaries on intermittent jets. 
\subsubsection{Choice of parameters}\label{Subsubsection 4.3.1}
We fix $L$ sufficiently large so that in particular it satisfies $L > 16 \vee (c_{R}^{-1} 153 (2\pi)^{\frac{3}{2}})$. We will need $\alpha < \frac{5-4m}{(24)(16)}$; e.g., for simplicity we can fix 
\begin{equation}\label{alpha}
\alpha \triangleq \frac{5-4m}{480}. 
\end{equation} 
For the $\alpha> 0$ fixed, we fix $b \in \mathbb{N}$ such that $b >  L^{2}\vee \frac{16}{\alpha}$. We will need 
\begin{equation}\label{estimate 49}
b a^{b(\frac{5-4m}{48})} a^{b [ -\frac{49 \alpha}{24} + \frac{4m-5}{24}]} \ll 1,  
\end{equation} 
which is possible because $4m - 5 < 0$. We point out that that by requiring $a > 0$ to be sufficiently large so that 
\begin{equation}
e^{2} \leq a^{\frac{5-4m}{48}}, 
\end{equation}
and $\beta> 0$ to be sufficiently small so that $\alpha > 16 \beta b$, \eqref{estimate 49} gives us 
\begin{equation}\label{[Equation (4.18), HZZ19]}
M_{0}(L)^{\frac{1}{2}} \lambda_{q+1}^{\frac{ 4m -5  - 52 \alpha}{24}} \delta_{q+2}^{-1} \ll 1. 
\end{equation} 
Let us point out that our condition of $\alpha > 16 \beta b$ is slightly more stringent than the analogous condition in \cite{HZZ19}, which was $\alpha > 8 \beta b$. The only reason for this change is upon deriving \eqref{[Equation (4.48h), HZZ19]} from \eqref{estimate 46}. Next we fix 
\begin{equation}\label{[Equation (4.17), HZZ19]}
l \triangleq \lambda_{q+1}^{- \frac{3\alpha}{2}} \lambda_{q}^{-2};
\end{equation} 
together with $b > L^{2} \vee \frac{16}{\alpha}$, we obtain by taking $a > 0$ sufficiently large  
\begin{equation}\label{[Equation (4.16), HZZ19]}
l \lambda_{q}^{4} \leq \lambda_{q+1}^{-\alpha}, \hspace{1mm} 4L \leq l^{-1} \hspace{1mm} \text{ and } \hspace{1mm} l^{-1} \leq \lambda_{q+1}^{2\alpha}. 
\end{equation} 
Finally, the last inequality of $c_{R} L \leq c_{R} (\frac{ (2\pi)^{\frac{3}{2}} a^{4} -2}{2})$ in \eqref{[Equation (4.7), HZZ19]} is immediately satisfied by taking $a$ sufficiently large. The inequalities of $(17) (2\pi)^{\frac{3}{2}} (9) < (17) (2\pi)^{\frac{3}{2}} a^{2\beta b} \leq c_{R} L$ of \eqref{[Equation (4.7), HZZ19]} can be satisfied by taking $\beta > 0$ sufficiently small. This is possible because we chose $L$ such that $16 \vee (c_{R}^{-1} 153 (2\pi)^{\frac{3}{2}}) < L$. We shall hereafter consider such $\alpha, L$ and $b$ fixed, and take $\beta > 0$ as small as we wish while $a > 0$ as large as we wish satisfying $a^{\frac{ 25 - 20 m}{24}} \in \mathbb{N}$, which will be necessary at the selection of $r_{\bot}, r_{\lVert}$ and $\mu$ in \eqref{[Equation (4.23), HZZ19]} . We chose these parameters carefully; we point out that $\alpha$ depends on $m$ and consequently, $b, l$, $\beta$ and $a$ all depend on $m$. 

\subsubsection{Mollification}\label{Mollification 1}
We let $\{ \phi_{\epsilon} \}_{\epsilon > 0}$ and $ \{ \varphi_{\epsilon}\}_{\epsilon > 0}$ respectively be families of standard mollifiers  on $\mathbb{R}^{3}$ and $\mathbb{R}$ with mass one where the latter is equipped with a compact support on $\mathbb{R}_{+}$. Then we mollify $v_{q}, \mathring{R}_{q}$ and $z$ in space and time to obtain 
\begin{equation}\label{[Equation (4.18c), HZZ19]}
v_{l} \triangleq (v_{q} \ast_{x} \phi_{l}) \ast_{t} \varphi_{l}, \hspace{3mm} \mathring{R}_{l} \triangleq (\mathring{R}_{q} \ast_{x} \phi_{l}) \ast_{t} \varphi_{l}, \hspace{3mm} z_{l} \triangleq (z \ast_{x} \phi_{l}) \ast_{t} \varphi_{l}, 
\end{equation} 
where $\phi_{l} \triangleq \frac{1}{l^{3}} \phi(\frac{\cdot}{l})$ and $\varphi_{l} \triangleq \frac{1}{l} \varphi( \frac{\cdot}{l})$, respectively. It follows from \eqref{[Equation (4.1), HZZ19]} that $(v_{l}, \mathring{R}_{l})$ solves 
\begin{equation}\label{[Equation (4.19), HZZ19]}
\partial_{t}v_{l} + (-\Delta)^{m} v_{l} + \text{div} ((v_{l} + z_{l}) \otimes (v_{l} + z_{l})) + \nabla \pi_{l} = \text{div} (\mathring{R}_{l} + R_{\text{commutator 1}}), \hspace{3mm} \nabla\cdot v_{l} = 0 
\end{equation}
for $t > 0$, where 
\begin{subequations}\label{[Equation (4.19a), HZZ19]}
\begin{align} 
R_{\text{commutator 1}} \triangleq& (v_{l} + z_{l}) \mathring{\otimes} (v_{l} + z_{l}) - (( ( v_{q} + z) \mathring{\otimes} (v_{q} + z)) \ast_{x} \phi_{l}) \ast_{t}\varphi_{l}, \\
\pi_{l} \triangleq& (\pi_{q} \ast_{x} \phi_{l}) \ast_{t} \varphi_{l} - \frac{1}{3} ( \lvert v_{l} + z_{l} \rvert^{2} - ( \lvert v_{q} + z \rvert^{2} \ast_{x} \phi_{l} ) \ast_{t} \varphi_{l}). 
\end{align}
\end{subequations}
For the mollified velocity field, we can verify by taking $a > 0$ sufficiently large and relying on that $\beta < \alpha$, which is from $16 \beta b < \alpha$, that for all $t \in [0, T_{L}]$ and $N \geq 1$, 
\begin{subequations}
\begin{align}
&\lVert v_{q} - v_{l} \rVert_{C_{t}L_{x}^{2}} \overset{\eqref{[Equation (4.4b), HZZ19]}}{\lesssim} l M_{0}(t)^{\frac{1}{2}} \lambda_{q}^{4}  \overset{\eqref{[Equation (4.16), HZZ19]}}{\lesssim} M_{0}(t)^{\frac{1}{2}} \lambda_{q+1}^{-\alpha}\leq  \frac{1}{4} M_{0}(t)^{\frac{1}{2}} \delta_{q+1}^{\frac{1}{2}}, \label{[Equation (4.20), HZZ19]}\\
& \lVert v_{l} \rVert_{C_{t}L_{x}^{2}} \leq \lVert v_{q} \rVert_{C_{t}L_{x}^{2}} \overset{\eqref{[Equation (4.4a), HZZ19]}}{\leq} M_{0}(t)^{\frac{1}{2}} (1+ \sum_{1 \leq r \leq q} \delta_{r}^{\frac{1}{2}}), \label{[Equation (4.21), HZZ19]}\\
&\lVert v_{l}  \rVert_{C_{t,x}^{N}} \lesssim l^{-(N-1)} \lVert v_{q} \rVert_{C_{t,x}^{1}} 
\overset{\eqref{[Equation (4.4b), HZZ19]}}{\lesssim} l^{-N+1} M_{0}(t)^{\frac{1}{2}} \lambda_{q}^{4} 
\overset{\eqref{[Equation (4.17), HZZ19]}}{\leq} l^{-N} M_{0}(t)^{\frac{1}{2}} \lambda_{q+1}^{-\alpha}. \label{[Equation (4.22), HZZ19]} 
\end{align} 
\end{subequations}
Now we choose the following parameters, which depend on the value of $m$:
\begin{equation}\label{[Equation (4.23), HZZ19]} 
r_{\lVert} \triangleq \lambda_{q+1}^{\frac{13 - 20m}{12}}, \hspace{3mm} r_{\bot} \triangleq \lambda_{q+1}^{\frac{ 1 - 20m}{24}}, \hspace{3mm} \mu \triangleq \frac{ \lambda_{q+1}^{2m-1} r_{\lVert}}{r_{\bot}} = \lambda_{q+1}^{2m-1} \lambda_{q+1}^{\frac{ 25 - 20 m}{24}}. 
\end{equation} 
Interestingly, they do not recover the choice from \cite[Equation (4.23)]{HZZ19} or \cite[Equation (7.25)]{BV19b} when $m = 1$; we refer to \cite{BCV18} for some heuristics on these choices. Considering that $m \in (\frac{13}{20}, \frac{5}{4})$, it can be readily verified that these choices of $r_{\bot}$ and $r_{\lVert}$ satisfy 
\begin{equation}\label{[Equation (B.2a), HZZ19]}
r_{\bot} \ll r_{\lVert} \ll 1 \text{ and } r_{\bot}^{-1} \ll \lambda_{q+1}
\end{equation} 
once we take $a> 0$ sufficiently large. In order to ensure the appropriate periodicity of $W_{(\xi)}, V_{(\xi)}, \Phi_{(\xi)}, \phi_{(\xi)}$ and $\psi_{(\xi)}$ defined respectively in \eqref{[Equation (B.3), HZZ19]}, \eqref{[Equation (B.6), HZZ19]}, \eqref{[Equation (B.2c), HZZ19]}, we make sure that $\lambda_{q+1} r_{\bot} = a^{b^{q+1} (\frac{25 - 20m}{24})} \in \mathbb{N}$ where we carefully recall that $m \in (\frac{13}{20}, \frac{5}{4})$. For this purpose, it suffices that $a^{\frac{25 - 20m}{24}} \in \mathbb{N}$ because $b \in \mathbb{N}$ as we described already in Sub-subsection \ref{Subsubsection 4.3.1}. 

\subsubsection{Perturbation}
Now we let $\chi$ be a smooth function such that 
\begin{equation}\label{[Equation (4.23a), HZZ19]}
\chi(z) \triangleq 
\begin{cases}
1 & \text{ if } z \in [0,1], \\
z & \text{ if } z \in [2, \infty),  
\end{cases}
\end{equation} 
and $z \leq 2 \chi(z) \leq 4z$ for $z \in (1,2)$. We define for $t \in [0, T_{L}]$ and $\omega \in \Omega$, 
\begin{equation}\label{[Equation (4.23b), HZZ19]}
\rho(\omega, t, x) \triangleq 4c_{R} \delta_{q+1} M_{0}(t) \chi ((c_{R} \delta_{q+1} M_{0}(t))^{-1} \lvert \mathring{R}_{l} (\omega, t, x) \rvert).
\end{equation} 
It follows immediately from \eqref{[Equation (4.23a), HZZ19]} that 
\begin{equation}\label{[Equation (4.23c), HZZ19]}
\lvert \frac{ \mathring{R}_{l}(\omega, t, x)}{\rho(\omega, t, x) } \rvert =  \frac{ \lvert \mathring{R}_{l}(\omega, t, x) \rvert}{4 c_{R} \delta_{q+1} M_{0}(t) \chi ((c_{R} \delta_{q+1} M_{0}(t))^{-1} \lvert \mathring{R}_{l} (\omega, t, x) \rvert )} \leq \frac{1}{2}.
\end{equation} 
We may compute directly from \eqref{[Equation (4.23a), HZZ19]} and \eqref{[Equation (4.23b), HZZ19]} for any $p \in [1, \infty]$ and $t \in [0, T_{L}]$, 
\begin{align}
\lVert \rho \rVert_{C_{t}L_{x}^{p}} \leq& \sup_{s \in [0,t]} 4 c_{R} \delta_{q+1} M_{0}(s) \lVert 1 + 3 (c_{R} \delta_{q+1} M_{0}(s))^{-1} \lvert \mathring{R}_{l} (\omega, s, x) \rvert \rVert_{L_{x}^{p}} \nonumber\\
\leq& 12 ((8\pi^{3})^{\frac{1}{p}} c_{R} \delta_{q+1} M_{0}(t) + \lVert \mathring{R}_{l} \rVert_{C_{t}L_{x}^{p}}). \label{[Equation (4.24), HZZ19]}
\end{align} 
Next, for any $N \geq 0$ and $t \in [0,T_{L}]$, we have due to  $W^{4,1}(\mathbb{T}^{3}) \hookrightarrow L^{\infty}(\mathbb{T}^{3})$, 
\begin{equation}\label{[Equation (4.24a), HZZ19]}
\lVert \mathring{R}_{l} \rVert_{C_{t,x}^{N}} \overset{\eqref{estimate 48}}{\lesssim} \sum_{0 \leq n + \lvert \alpha \rvert \leq N} \lVert \partial_{t}^{n} D^{\alpha} (-\Delta)^{2} \mathring{R}_{l} \rVert_{L_{t}^{\infty} L_{x}^{1}} \overset{\eqref{[Equation (4.4c), HZZ19]}}{\lesssim} l^{-4-N} M_{0}(t) c_{R} \delta_{q+1}.
\end{equation} 
This leads us to for any $N \geq 0$, 
\begin{subequations}\label{[Equation (4.25), HZZ19]}
\begin{align}
& \lVert \rho \rVert_{C_{t}C_{x}^{N}} \lesssim c_{R} \delta_{q+1} M_{0}(t) l^{-4-N}, \hspace{10mm} \lVert \rho \rVert_{C_{t}^{1}C_{x}} \lesssim c_{R} \delta_{q+1} M_{0}(t) l^{-5}, \\
& \lVert \rho \rVert_{C_{t}^{1}C_{x}^{1}} \lesssim c_{R} \delta_{q+1} M_{0}(t) l^{-10}, \hspace{13mm} \lVert \rho \rVert_{C_{t}^{1}C_{x}^{2}} \lesssim c_{R} \delta_{q+1} M_{0}(t) l^{-15}. 
\end{align}
\end{subequations}
The first inequality in \eqref{[Equation (4.25), HZZ19]} may be computed relying on \eqref{[Equation (4.24), HZZ19]} in case $N = 0$ and \cite[Proposition C.1]{BDIS15} in case $N \geq 1$ (cf. \cite[Section 7.5.2]{BV19b}), while one can directly apply derivatives $\partial_{t}, \nabla$, use \eqref{[Equation (4.24a), HZZ19]}, that $\partial_{t} M_{0}(t) = 4L M_{0}(t)$ and $4L \leq l^{-1}$ from \eqref{[Equation (4.16), HZZ19]} to deduce the other three inequalities. Next we define the amplitude function
\begin{equation}\label{[Equation (4.26), HZZ19]}
a_{(\xi)} (\omega, t,x) \triangleq a_{\xi, q+1} (\omega, t,x) \triangleq \rho(\omega, t,x)^{\frac{1}{2}} \gamma_{\xi} ( \text{Id} - \frac{ \mathring{R}_{l} (\omega, t, x)}{\rho (\omega, t, x)} ) (2\pi)^{-\frac{3}{4}}, 
\end{equation} 
where $\text{Id} - \rho^{-1} \mathring{R}_{l}$ lies in the domain of $\gamma_{\xi}$ defined in Lemma \ref{[Lemma B.1, HZZ19]} due to \eqref{[Equation (4.23c), HZZ19]}. Moreover, 
\begin{equation}\label{[Equation (4.27), HZZ19]}
(2\pi)^{\frac{3}{2}} \sum_{\xi \in \Lambda} a_{(\xi)}^{2} \fint_{\mathbb{T}^{3}} W_{(\xi)} \otimes W_{(\xi)} dy  \overset{\eqref{[Equation (4.26), HZZ19]} \eqref{[Equation (B.5), HZZ19]}}{=} \rho \text{Id} - \mathring{R}_{l}
\end{equation} 
where $W_{(\xi)}$ is the intermittent jet defined in \eqref{[Equation (B.3), HZZ19]}. For all $t \in [0, T_{L}]$, with $\Lambda$ from Lemma \ref{[Lemma B.1, HZZ19]}, along with $C_{\Lambda}$ and $M$ from \eqref{[Equation (A.9e), HZZ19]} we can estimate
\begin{equation}\label{[Equation (4.28), HZZ19]}
\lVert a_{(\xi)} \rVert_{C_{t}L_{x}^{2}} \overset{\eqref{[Equation (4.23c), HZZ19]} \eqref{[Equation (4.24), HZZ19]} \eqref{[Equation (A.9e), HZZ19]}}{\lesssim} 4((8\pi^{3}) c_{R} \delta_{q+1} M_{0} (t) + \lVert \mathring{R}_{l} \rVert_{C_{t} L^{1}})^{\frac{1}{2}} (\frac{M}{C_{\Lambda}})\overset{\eqref{[Equation (4.4c), HZZ19]} \eqref{[Equation (A.9e), HZZ19]}}{\leq} \frac{c_{R}^{\frac{1}{4}} M_{0}(t)^{\frac{1}{2}} \delta_{q+1}^{\frac{1}{2}}}{2 \lvert \Lambda \rvert} 
\end{equation} 
where we required $c_{R}^{\frac{1}{4}} \ll \frac{1}{M}$ in the last inequality. Using the fact that $\rho(t) \geq 2c_{R} \delta_{q+1} M_{0}(t)$ from \eqref{[Equation (4.23b), HZZ19]} and \eqref{[Equation (4.23a), HZZ19]}, we have 
\begin{subequations}\label{[Equation (4.29), HZZ19]} for any $N \geq 0$
\begin{align}
& \lVert a_{(\xi)} \rVert_{C_{t} C_{x}^{N}} \leq c_{R}^{\frac{1}{4}} \delta_{q+1}^{\frac{1}{2}} M_{0}(t)^{\frac{1}{2}} l^{-2 - 5N}, \hspace{5mm} \lVert a_{(\xi)} \rVert_{C_{t}^{1} C_{x}} \leq c_{R}^{\frac{1}{4}} \delta_{q+1}^{\frac{1}{2}} M_{0}(t)^{\frac{1}{2}} l^{-5}, \\
& \lVert a_{(\xi)} \rVert_{C_{t}^{1}C_{x}^{1}} \leq c_{R}^{\frac{1}{4}} \delta_{q+1}^{\frac{1}{2}} M_{0}(t)^{\frac{1}{2}} l^{-10}, \hspace{9mm}  \lVert a_{(\xi)} \rVert_{C_{t}^{1}C_{x}^{2}} \leq c_{R}^{\frac{1}{4}} \delta_{q+1}^{\frac{1}{2}} M_{0}(t)^{\frac{1}{2}} l^{-15}. 
\end{align} 
\end{subequations}
These inequalities can be estimated very similarly to \eqref{[Equation (4.25), HZZ19]} while we also actually use \eqref{[Equation (4.25), HZZ19]}; e.g., the first is done by using \eqref{[Equation (4.26), HZZ19]}, product estimate (e.g., \cite[Lemma A,4]{W04}) and \cite[Proposition (C.1)]{BDIS15} while the rest can be seen by directly applying $\partial_{t}, \nabla$,  and using \eqref{[Equation (4.23c), HZZ19]} and \eqref{[Equation (4.24a), HZZ19]}. We also took $c_{R}$ sufficiently small depending on $C_{\lambda}$ and $M$. E.g., 
\begin{align*}
 \lVert a_{(\xi)} \rVert_{C_{t}^{1}C_{x}} 
\lesssim \sup_{s \in [0,t]} [  \lVert \frac{1}{\rho} \rVert_{L_{x}^{\infty}}^{\frac{1}{2}} \lVert \partial_{s} \rho \rVert_{L_{x}^{\infty}} + ( \lVert \partial_{s} \mathring{R}_{l} \rVert_{L_{x}^{\infty}} \lVert \frac{1}{\rho} \rVert_{L_{x}^{\infty}}^{\frac{1}{2}} + \lVert \partial_{s} \rho \rVert_{L_{x}^{\infty}} \lVert \frac{1}{\rho} \rVert_{L_{x}^{\infty}}^{\frac{1}{2}})] \leq c_{R}^{\frac{1}{4}} \delta_{q+1}^{\frac{1}{2}} M_{0}(t)^{\frac{1}{2}} l^{-5}. 
\end{align*} 
Next we define the principal part, incompressibility corrector, and the temporal corrector of the perturbation $w_{q+1}$ respectively as 
\begin{subequations}
\begin{align}
&w_{q+1}^{(p)} \triangleq \sum_{\xi \in \Lambda} a_{(\xi)} W_{(\xi)}, \label{[Equation (4.30), HZZ19]} \\
&w_{q+1}^{(c)} \triangleq \sum_{\xi \in \Lambda} \text{curl} (\nabla a_{(\xi)} \times V_{(\xi)}) + \nabla a_{(\xi)} \times \text{curl} V_{(\xi)} + a_{(\xi)} W_{(\xi)}^{(c)}, \label{[Equation (4.32), HZZ19]}\\
& w_{q+1}^{(t)} \triangleq - \mu^{-1} \sum_{\xi \in \Lambda} \mathbb{P} \mathbb{P}_{\neq 0} (a_{(\xi)}^{2} \phi_{(\xi)}^{2} \psi_{(\xi)}^{2} \xi), \label{[Equation (4.33), HZZ19]}
\end{align}
\end{subequations} 
where $W_{(\xi)}^{(c)}$ is defined in \eqref{[Equation (B.6), HZZ19]}. They satisfy the following important identities:
\begin{subequations}
\begin{align} 
&w_{q+1}^{(p)} \otimes w_{q+1}^{(p)} + \mathring{R}_{l} \overset{ \eqref{[Equation (4.27), HZZ19]} \eqref{[Equation (4.30), HZZ19]} \eqref{[Equation (B.4), HZZ19]}}{=} \sum_{\xi \in \Lambda}a_{(\xi)}^{2} \mathbb{P}_{\neq 0} (W_{(\xi)} \otimes W_{(\xi)} ) + \rho \text{Id}, \label{[Equation (4.31), HZZ19]} \\
&\sum_{\xi \in \Lambda} \text{curl curl} (a_{(\xi)} V_{(\xi)}) \overset{\eqref{[Equation (4.30), HZZ19]} \eqref{[Equation (4.32), HZZ19]}\eqref{[Equation (B.6), HZZ19]} }{=} w_{q+1}^{(p)} + w_{q+1}^{(c)}, \label{[Equation (4.32a), HZZ19]} \\
& \partial_{t} w_{q+1}^{(t)} + \sum_{\xi \in \Lambda} \mathbb{P}_{\neq 0} (a_{(\xi)}^{2} \text{div} (W_{(\xi)} \otimes W_{(\xi)}) ) \overset{\eqref{[Equation (4.33), HZZ19]} \eqref{estimate 7}}{=} \nabla \pi_{1} - \mu^{-1} \sum_{\xi \in \Lambda} \mathbb{P}_{ \neq 0} [(\partial_{t} a_{(\xi)}^{2}) \phi_{(\xi)}^{2} \psi_{(\xi)}^{2} \xi], \label{[Equation (4.34), HZZ19]}\\
& \hspace{47mm} \text{ where } \pi_{1} \triangleq \Delta^{-1} \text{div } \mu^{-1} \sum_{\xi \in \Lambda} \mathbb{P}_{\neq 0} \partial_{t} (a_{(\xi)}^{2} \phi_{(\xi)}^{2} \psi_{(\xi)}^{2} \xi). \nonumber
\end{align} 
\end{subequations} 
Due to \eqref{[Equation (4.32a), HZZ19]} we see that $\text{div} (w_{q+1}^{(p)} + w_{q+1}^{(c)}) = 0$. At last, we  define 
\begin{equation}\label{[Equation (4.35), HZZ19]}
w_{q+1} \triangleq w_{q+1}^{(p)} + w_{q+1}^{(c)} + w_{q+1}^{(t)} \text{ and } v_{q+1} \triangleq v_{l} + w_{q+1}, 
\end{equation} 
where $w_{q+1}$ is mean zero and divergence-free. Due to \eqref{[Equation (4.28), HZZ19]}, \eqref{[Equation (4.29), HZZ19]} and the fact that $\lVert W_{(\xi)} \rVert_{L^{2}} \overset{\eqref{[Equation (B.4a), HZZ19]}}{=} (2\pi)^{\frac{3}{2}}$, we may apply Lemma \ref{[Lemma 7.4, BV19b]} with ``$\zeta$''  and  ``$\kappa$''  therein being respectively $l^{-8}$ and $\lambda_{q+1} r_{\bot}$ to deduce 
\begin{equation}\label{[Equation (4.37), HZZ19]}
\lVert w_{q+1}^{(p)} \rVert_{C_{t}L_{x}^{2}} \overset{\eqref{[Equation (4.30), HZZ19]}}{\lesssim} \sup_{s \in [0,t]} \sum_{\xi \in \Lambda} \frac{ c_{R}^{\frac{1}{4}} \delta_{q+1}^{\frac{1}{2}} M_{0}(s)^{\frac{1}{2}}}{2 \lvert \Lambda \rvert } \lVert W_{(\xi)}(s) \rVert_{L_{x}^{2}}  \leq \frac{M_{0}(t)^{\frac{1}{2}}}{2} \delta_{q+1}^{\frac{1}{2}} 
\end{equation} 
where we took $c_{R} \ll (2\pi)^{-6}$. We emphasize that here we used the key assumption that $\lambda_{q+1} r_{\bot} = a^{b^{q+1} (\frac{25 - 20 m}{24})} \in \mathbb{N}$ due to our choice of $a$ and $b$; we also had to satisfy one of the hypothesis \eqref{[Equation (7.41), BV19b]}
\begin{equation*}
\frac{ 2 \pi \sqrt{3} l^{-8}}{\lambda_{q+1} r_{\bot}} \leq \frac{1}{3}
\end{equation*}
which required \eqref{[Equation (4.16), HZZ19]}, as well as $\alpha < \frac{25  - 20m}{16 (24)}$ that holds due to our choice from \eqref{alpha}. Next, for all $p \in (1, \infty)$ and $t \in [0, T_{L}]$, we can estimate 
\begin{subequations}\label{estimate 33}
\begin{align}
& \lVert w_{q+1}^{(p)} \rVert_{C_{t}L_{x}^{p}} \overset{\eqref{[Equation (4.30), HZZ19]}}{\lesssim} \sup_{s \in [0,t]} \sum_{\xi \in\Lambda} \lVert a_{(\xi)} (s) \rVert_{L_{x}^{\infty}} \lVert W_{(\xi)}(s) \rVert_{L_{x}^{p}} \overset{\eqref{[Equation (4.29), HZZ19]} \eqref{[Equation (B.7), HZZ19]}}{\lesssim}  M_{0}(t)^{\frac{1}{2}} \delta_{q+1}^{\frac{1}{2}} l^{-2} r_{\bot}^{\frac{2}{p} -1} r_{\lVert}^{\frac{1}{p} - \frac{1}{2}}, \label{[Equation (4.38), HZZ19]}\\
& \lVert w_{q+1}^{(c)} \rVert_{C_{t}L_{x}^{p}} \overset{\eqref{[Equation (4.32), HZZ19]}}{\lesssim} \sup_{s \in [0,t]} \sum_{\xi \in \Lambda} \lVert a_{(\xi)}(s) \rVert_{C_{x}^{2}} \lVert V_{(\xi)}(s) \rVert_{W_{x}^{1,p}} + \lVert a_{(\xi)}(s) \rVert_{C_{x}} \lVert W_{(\xi)}^{(c)}(s) \rVert_{L_{x}^{p}} \nonumber\\
& \hspace{50mm}  \overset{\eqref{[Equation (4.29), HZZ19]} \eqref{[Equation (B.7), HZZ19]}}{\lesssim} M_{0}(t)^{\frac{1}{2}} \delta_{q+1}^{\frac{1}{2}} l^{-12} r_{\bot}^{\frac{2}{p}} r_{\lVert}^{\frac{1}{p} - \frac{3}{2}}, \label{[Equation (4.39), HZZ19]}\\
& \lVert w_{q+1}^{(t)} \rVert_{C_{t}L_{x}^{p}} \overset{\eqref{[Equation (4.33), HZZ19]}}{\lesssim} \sup_{s \in [0,t]} \mu^{-1} \sum_{\xi \in \Lambda} \lVert a_{(\xi)}(s) \rVert_{L_{x}^{\infty}}^{2} \lVert \phi_{(\xi)}\rVert_{L_{x}^{2p}}^{2} \lVert \psi_{(\xi)}(s) \rVert_{L_{x}^{2p}}^{2} \nonumber\\
& \hspace{50mm}  \overset{\eqref{[Equation (4.29), HZZ19]} \eqref{[Equation (B.7), HZZ19]}}{\lesssim} \delta_{q+1} M_{0}(t) l^{-4} r_{\bot}^{\frac{2}{p} -1} r_{\lVert}^{\frac{1}{p} - 2} \lambda_{q+1}^{1-2m}, \label{[Equation (4.40), HZZ19]}
\end{align} 
\end{subequations}
where the first inequality of \eqref{[Equation (4.40), HZZ19]} made use of the fact that $\phi_{(\xi)}$ and $\psi_{(\xi)}$ from \eqref{[Equation (B.2c), HZZ19]} are functions of variables in orthogonal directions. We also emphasize that our estimate of \eqref{[Equation (4.40), HZZ19]} differs from \cite[Equation (4.40)]{HZZ19} in that ``$\lambda_{q+1}^{-1}$'' therein is replaced by ``$\lambda_{q+1}^{1-2m}$ which is crucial. For $t \in [0, T_{L}]$, because \eqref{alpha} gives us $20 \alpha + \frac{20m-25}{24} < 0$ and $4 \alpha + \frac{m}{2} - \frac{5}{8} < 0$, we can take $a> 0$ sufficiently large so that 
\begin{equation}\label{estimate 52}
\lambda_{q+1}^{20 \alpha + \frac{20m-25}{24}} + M_{0}(L)^{\frac{1}{2}} \lambda_{q+1}^{4\alpha + \frac{m}{2} - \frac{5}{8}} < 1
 \end{equation}
 to deduce from \eqref{[Equation (4.39), HZZ19]}, \eqref{[Equation (4.40), HZZ19]} 
\begin{align}
& \lVert w_{q+1}^{(c)} \rVert_{C_{t}L_{x}^{p}} + \lVert w_{q+1}^{(t)} \rVert_{C_{t}L_{x}^{p}} \label{[Equation (4.41), HZZ19]}\\
\overset{\eqref{[Equation (4.16), HZZ19]} \eqref{[Equation (4.23), HZZ19]}}{\lesssim}
& M_{0}(t)^{\frac{1}{2}} \delta_{q+1}^{\frac{1}{2}} l^{-2} r_{\bot}^{\frac{2}{p} - 1} r_{\lVert}^{\frac{1}{p} - \frac{1}{2}}[ \lambda_{q+1}^{20 \alpha + \frac{20m - 25}{24}} + \delta_{q+1}^{\frac{1}{2}} M_{0}(t)^{\frac{1}{2}} \lambda_{q+1}^{4 \alpha + \frac{m}{2} - \frac{5}{8}}] 
\lesssim M_{0}(t)^{\frac{1}{2}} \delta_{q+1}^{\frac{1}{2}} l^{-2} r_{\bot}^{\frac{2}{p} -1} r_{\lVert}^{\frac{1}{p} - \frac{1}{2}}.\nonumber
\end{align} 
Similarly to \eqref{estimate 52}, because $24 \alpha + \frac{20m - 25}{24} < 0$ and $8\alpha + \frac{m}{2} - \frac{5}{8} < 0$ due to \eqref{alpha}, we may take $a> 0$ sufficiently large to deduce from \eqref{[Equation (4.35), HZZ19]}
\begin{align}
&\lVert w_{q+1} \rVert_{C_{t}L_{x}^{2}} \label{[Equation (4.42), HZZ19]}\\
\overset{\eqref{[Equation (4.16), HZZ19]} \eqref{[Equation (4.23), HZZ19]} \eqref{[Equation (4.37), HZZ19]}  \eqref{estimate 33} }{\leq}
& M_{0}(t)^{\frac{1}{2}} \delta_{q+1}^{\frac{1}{2}}[ \frac{1}{2}+  C \lambda_{q+1}^{24 \alpha + \frac{20m -25}{24}} + C M_{0}(t)^{\frac{1}{2}} \delta_{q+1}^{\frac{1}{2}} \lambda_{q+1}^{8\alpha + \frac{m}{2} - \frac{5}{8}}] \leq \frac{3}{4} M_{0}(t)^{\frac{1}{2}} \delta_{q+1}^{\frac{1}{2}}. \nonumber
\end{align} 
At this point, one can immediately verify  \eqref{[Equation (4.4a), HZZ19]} at level $q+ 1$ by \eqref{[Equation (4.35), HZZ19]}, \eqref{[Equation (4.21), HZZ19]} and \eqref{[Equation (4.42), HZZ19]}, as well as \eqref{[Equation (4.12), HZZ19]} by \eqref{[Equation (4.35), HZZ19]}, \eqref{[Equation (4.42), HZZ19]} and \eqref{[Equation (4.20), HZZ19]}. 

Next, for all $t \in [0, T_{L}]$, 
\begin{subequations}\label{estimate 15}
\begin{align}
& \lVert w_{q+1}^{(p)} \rVert_{C_{t,x}^{1}} \overset{ \eqref{[Equation (4.29), HZZ19]} \eqref{[Equation (4.30), HZZ19]} \eqref{[Equation (B.7c), HZZ19]}}{\lesssim} M_{0}(t)^{\frac{1}{2}} l^{-7} r_{\bot}^{-1} r_{\lVert}^{-\frac{1}{2}} \lambda_{q+1} (1+ \frac{r_{\bot}\mu}{r_{\lVert}}) \overset{\eqref{[Equation (4.23), HZZ19]}}{\lesssim} M_{0}(t)^{\frac{1}{2}} l^{-7} r_{\bot}^{-1} r_{\lVert}^{-\frac{1}{2}} \lambda_{q+1}^{2m}, \label{[Equation (4.43), HZZ19]}\\
& \lVert w_{q+1}^{(c)} \rVert_{C_{t,x}^{1}} \overset{ \eqref{[Equation (4.29), HZZ19]} \eqref{[Equation (4.32), HZZ19]} \eqref{[Equation (B.7c), HZZ19]}}{\lesssim} M_{0}(t)^{\frac{1}{2}} l^{-17} r_{\lVert}^{-\frac{3}{2}} \lambda_{q+1} (1+ \frac{r_{\bot} \mu}{r_{\lVert}}) \overset{\eqref{[Equation (4.23), HZZ19]}}{\lesssim} M_{0}(t)^{\frac{1}{2}} l^{-17} r_{\lVert}^{-\frac{3}{2}} \lambda_{q+1}^{2m}.  \label{[Equation (4.44), HZZ19]}
\end{align} 
\end{subequations}
On the other hand, we need to estimate $\lVert w_{q+1}^{(t)}\rVert_{C_{t,x}^{1}}$ with a bit more care. First, because $\mathbb{P} \mathbb{P}_{\neq 0}$ is not bounded in $C_{x}$, we go down to $L^{p}$ for $p < \infty$ in the expense of $\lambda_{q+1}^{\alpha}$, and then use $L^{\infty} (\mathbb{T}^{3}) \hookrightarrow L^{p} (\mathbb{T}^{3})$ to compute 
\begin{align}
\lVert w_{q+1}^{(t)} \rVert_{C_{t,x}^{1}} 
\overset{\eqref{[Equation (4.29), HZZ19]} \eqref{[Equation (4.33), HZZ19]} \eqref{[Equation (B.7), HZZ19]}}{\lesssim}& M_{0}(t) l^{-9} r_{\bot}^{-1} r_{\lVert}^{-2} \lambda_{q+1}^{\alpha} [ l^{2} r_{\bot}^{-1} r_{\lVert} \lambda_{q+1}^{ \frac{ - 28m - 1}{24}} + l^{5} r_{\bot}^{-1} r_{\lVert} \lambda_{q+1}^{ \frac{ - 28m -1}{24}} \lambda_{q+1} (\frac{r_{\bot} \mu}{r_{\lVert}}) \nonumber\\
& \hspace{17mm} + r_{\bot}^{-1} r_{\lVert} \lambda_{q+1}^{\frac{ -28m -1}{24}} + l^{5} r_{\bot}^{-1} r_{\lVert} \lambda_{q+1}^{\frac{ - 28m -1}{24}} \lambda_{q+1} + l^{5}  \lambda_{q+1}^{\frac{ - 28m -1}{24}} \lambda_{q+1}]. \label{estimate 14}
\end{align} 
Due to \eqref{[Equation (4.16), HZZ19]} and \eqref{[Equation (4.23), HZZ19]} we can verify for all $m \in (\frac{13}{20}, \frac{5}{4})$ that 
\begin{subequations}\label{estimate 13}
\begin{align}
&  l^{2} r_{\bot}^{-1} r_{\lVert} \lambda_{q+1}^{\frac{ - 28m -1}{24}} \leq  \lambda_{q+1}^{-2\alpha} \lambda_{q}^{-8} \lambda_{q+1}^{\frac{ - 48m+24}{24}} \lesssim \lambda_{q+1}^{-2m+2} (\frac{r_{\bot} \mu}{r_{\lVert}}), \label{estimate 8} \\
& l^{5} r_{\bot}^{-1} r_{\lVert} \lambda_{q+1}^{\frac{ -28 m-1}{24}} \lambda_{q+1}(\frac{r_{\bot} \mu}{r_{\lVert}})   \leq  \lambda_{q+1}^{-5\alpha} \lambda_{q}^{-20} \lambda_{q+1}^{-2m + 2}(\frac{r_{\bot} \mu}{r_{\lVert}})  \leq \lambda_{q+1}^{-2m+2}(\frac{r_{\bot} \mu}{r_{\lVert}}) , \label{estimate 9} \\
&r_{\bot}^{-1} r_{\lVert} \lambda_{q+1}^{ \frac{-28m-1}{24}} =  \lambda_{q+1}^{-2m+1} \lesssim \lambda_{q+1}^{-2m+2} (\frac{r_{\bot} \mu}{r_{\lVert}}),\label{estimate 10}  \\
& l^{5} r_{\bot}^{-1} r_{\lVert}\lambda_{q+1}^{\frac{ - 28m-1}{24}} \lambda_{q+1}  \leq \lambda_{q+1}^{-5 \alpha} \lambda_{q}^{-20} \lambda_{q+1}^{-2 m + 2} \lesssim \lambda_{q+1}^{-2m+2} (\frac{r_{\bot} \mu}{r_{\lVert}}) , \label{estimate 11} \\
& l^{5}  \lambda_{q+1}^{\frac{ - 28m-1}{24}} \lambda_{q+1} \leq \lambda_{q+1}^{-5\alpha} \lambda_{q}^{-20} \lambda_{q+1}^{\frac{23 - 28m}{24}} \lesssim \lambda_{q+1}^{-2m+2} (\frac{r_{\bot} \mu}{r_{\lVert}}); \label{estimate 12} 
\end{align} 
\end{subequations}
as we will see in Remark \ref{star}, although $\lambda_{q+1}^{-2m+2} \lesssim 1$ in case $m \in [1, \frac{5}{4})$, it was crucial to keep ``$\lambda_{q+1}^{-2m+2}$'' in \eqref{estimate 9}. Applying \eqref{estimate 13} to \eqref{estimate 14} leads us to 
\begin{equation}\label{[Equation (4.45), HZZ19]}
\lVert w_{q+1}^{(t)} \rVert_{C_{t,x}^{1}} \lesssim  M_{0}(t) l^{-9} r_{\bot}^{-1} r_{\lVert}^{-2} \lambda_{q+1}^{-2m + 2 + \alpha} (\frac{r_{\bot} \mu}{r_{\lVert}}) 
\end{equation} 
due to \eqref{[Equation (4.23), HZZ19]}. This allows us to compute for all $t \in [0, T_{L}]$, 
\begin{align}
\lVert v_{q+1} \rVert_{C_{t,x}^{1}} \overset{\eqref{[Equation (4.22), HZZ19]} \eqref{[Equation (4.35), HZZ19]}  \eqref{estimate 15} }{\leq}& M_{0}(t)^{\frac{1}{2}} [l^{-1} \lambda_{q+1}^{-\alpha} + Cl^{-7} r_{\bot}^{-1} r_{\lVert}^{-\frac{1}{2}} \lambda_{q+1}^{2m} \nonumber\\
& \hspace{8mm} +  Cl^{-17} r_{\lVert}^{-\frac{3}{2}} \lambda_{q+1}^{2m} +  CM_{0}(t)^{\frac{1}{2}} l^{-9} r_{\bot}^{-1} r_{\lVert}^{-2} \lambda_{q+1}^{-2m + 2 + \alpha} (\frac{r_{\bot} \mu}{r_{\lVert}})]. \label{[Equation (4.45a), HZZ19]}
\end{align} 
Now our choice of $\alpha = \frac{5-4m}{480}$ from \eqref{alpha} gives us
\begin{equation*}
\alpha < (\frac{1}{14})(\frac{11}{12})(5-4m) \wedge (\frac{1}{34}) (\frac{9}{8}) (5- 4m) \wedge (\frac{1}{19}) (\frac{25}{24}) (5-4m)
\end{equation*}
which imply respectively  
\begin{equation*}
14 \alpha + \frac{44m -7}{12} < 4, \hspace{3mm} 34 \alpha + \frac{36 m - 13}{8} < 4, \hspace{3mm} 19 \alpha + \frac{100m - 29}{24} < 4, 
\end{equation*} 
and  leads us to, for $a > 0$ taken sufficiently large, 
\begin{subequations}
\begin{align}
&l^{-1} \lambda_{q+1}^{-\alpha} \overset{\eqref{[Equation (4.16), HZZ19]}}{\leq} \lambda_{q+1}^{\alpha} \leq \frac{1}{4} \lambda_{q+1}^{4}, \label{estimate 16} \\
&C l^{-7} r_{\bot}^{-1} r_{\lVert}^{-\frac{1}{2}} \lambda_{q+1}^{2m} \overset{\eqref{[Equation (4.16), HZZ19]}}{\leq} C \lambda_{q+1}^{14 \alpha + \frac{44m -7}{12}} \leq \frac{1}{4} \lambda_{q+1}^{4},  \label{estimate 17}\\
&C l^{-17} r_{\lVert}^{-\frac{3}{2}} \lambda_{q+1}^{2m} \overset{\eqref{[Equation (4.16), HZZ19]}}{\leq} C \lambda_{q+1}^{34 \alpha} \lambda_{q+1}^{\frac{36m-13}{8}} \leq  \frac{1}{4} \lambda_{q+1}^{4},   \label{estimate 18}\\
& CM_{0}(t)^{\frac{1}{2}} l^{-9} r_{\bot}^{-1} r_{\lVert}^{-2} \lambda_{q+1}^{-2 m + 2 + \alpha} (\frac{r_{\bot} \mu}{r_{\lVert}})  \overset{\eqref{[Equation (4.16), HZZ19]} \eqref{[Equation (4.23), HZZ19]}}{\leq} CM_{0}(t)^{\frac{1}{2}} \lambda_{q+1}^{19\alpha + \frac{100m - 29}{24}} \leq \frac{1}{4} \lambda_{q+1}^{4}.  \label{estimate 19}
\end{align}
\end{subequations}
\begin{remark}\label{star}
Let us emphasize that in case $m\in [1, \frac{5}{4})$, if we simply bounded ``$\lambda_{q+1}^{-2m+2}$'' in \eqref{estimate 13} by $\lambda_{q+1}^{-2m + 2} \lesssim 1$ which is valid, then the estimate \eqref{estimate 19} would have been impossible. Indeed, if we estimated $\lambda_{q+1}^{-2m + 2} \lesssim 1$ in case $m \in [1, \frac{5}{4})$, then we would have instead of \eqref{estimate 19}
\begin{equation*}
CM_{0}(t)^{\frac{1}{2}} l^{-9} r_{\bot}^{-1} r_{\lVert}^{-2} \lambda_{q+1}^{\alpha} (\frac{r_{\bot} \mu}{r_{\lVert}}) \overset{\eqref{[Equation (4.16), HZZ19]} \eqref{[Equation (4.23), HZZ19]}}{\leq} CM_{0}(t)^{\frac{1}{2}}  \lambda_{q+1}^{18 \alpha} \lambda_{q+1}^{\frac{20m - 1}{24}} \lambda_{q+1}^{\frac{20 m -13}{6}} \lambda_{q+1}^{\alpha} \lambda_{q+1}^{2m-1}
\end{equation*} 
and unfortunately  
\begin{equation*}
18 \alpha + \frac{20m-1}{24} + \frac{20m -13}{6}  + \alpha + 2m - 1 > 4 
\end{equation*}
only for $m > \frac{5}{4} - \frac{3}{37} > 1$.
\end{remark}
Therefore, applying \eqref{estimate 16}-\eqref{estimate 19} to \eqref{[Equation (4.45a), HZZ19]} implies that \eqref{[Equation (4.4b), HZZ19]} at level $q+1$ holds. 

Next we further estimate for $t \in [0, T_{L}]$ and $p \in (1, \infty)$, 
\begin{subequations}
\begin{align}
& \lVert w_{q+1}^{(p)} + w_{q+1}^{(c)} \rVert_{C_{t}W_{x}^{1,p}} 
\overset{\eqref{[Equation (4.32a), HZZ19]}}{\leq} \sum_{\xi \in \Lambda}  \lVert \text{curl curl} (a_{(\xi)} V_{(\xi)} ) \rVert_{C_{t}W_{x}^{1,p}}  \label{[Equation (4.46), HZZ19]}\\
& \hspace{40mm} \overset{\eqref{[Equation (4.29), HZZ19]} \eqref{[Equation (B.7c), HZZ19]}}{\lesssim} M_{0}(t)^{\frac{1}{2}} r_{\bot}^{\frac{2}{p} - 1} r_{\lVert}^{\frac{1}{p} - \frac{1}{2}} l^{-2} \lambda_{q+1}, \nonumber\\
& \lVert w_{q+1}^{(t)} \rVert_{C_{t}W_{x}^{1,p}} \overset{\eqref{[Equation (4.33), HZZ19]}}{\lesssim} \mu^{-1} \sum_{\xi \in \Lambda} \lVert a_{(\xi)} \rVert_{C_{t}C_{x}} \lVert a_{(\xi)} \rVert_{C_{t}C_{x}^{1}} \lVert \phi_{(\xi)} \rVert_{L_{x}^{2p}}^{2} \lVert \psi_{(\xi)} \rVert_{C_{t}L_{x}^{2p}}^{2} \label{[Equation (4.47), HZZ19]}\\
& \hspace{7mm} + \lVert a_{(\xi)} \rVert_{C_{t}C_{x}}^{2} \lVert \phi_{(\xi)} \rVert_{L_{x}^{2p}} \lVert \nabla \phi_{(\xi)} \rVert_{L_{x}^{2p}} \lVert \psi_{(\xi)} \rVert_{C_{t}L_{x}^{2p}}^{2} + \lVert a_{(\xi)} \rVert_{C_{t}C_{x}}^{2} \lVert \phi_{(\xi)} \rVert_{L_{x}^{2p}}^{2} \lVert \psi_{(\xi)} \rVert_{C_{t}L_{x}^{2p}} \lVert \nabla \psi_{(\xi)} \rVert_{C_{t}L_{x}^{2p}} \nonumber\\
&\overset{\eqref{[Equation (4.29), HZZ19]} \eqref{[Equation (B.7), HZZ19]}}{\lesssim} \mu^{-1} M_{0}(t) r_{\bot}^{\frac{2}{p} -2} r_{\lVert}^{\frac{1}{p} -1} [ l^{-9} + l^{-4} \lambda_{q+1} + l^{-4} \lambda_{q+1}^{\frac{ 5(4m-5)}{24}} \lambda_{q+1}] \lesssim \mu^{-1} M_{0}(t) r_{\bot}^{\frac{2}{p} -2} r_{\lVert}^{\frac{1}{p} - 1} l^{-4} \lambda_{q+1} \nonumber 
\end{align}
\end{subequations} 
by H$\ddot{\mathrm{o}}$lder's inequality and that  $\phi_{(\xi)}$ and $\psi_{(\xi)}$ from \eqref{[Equation (B.2c), HZZ19]} are functions of variables in orthogonal directions. In contrast to the case of the NS equations in \cite{HZZ19}, we need higher and lower order estimates. In fact, we can estimate $C_{t}L_{x}^{p}$ and $C_{t}W_{x}^{2,p}$-norms very similarly to \eqref{[Equation (4.46), HZZ19]} and \eqref{[Equation (4.47), HZZ19]} and interpolate between $C_{t}L_{x}^{p}$ and $C_{t}W_{x}^{1,p}$ in case $m \in (\frac{13}{20}, 1)$ and $C_{t}W_{x}^{1,p}$ and $C_{t}W_{x}^{2,p}$ in case $m \in [1, \frac{5}{4})$ to deduce for all $t \in [0, T_{L}]$ and $p\in (1,\infty)$ 
\begin{subequations}
\begin{align}
 \lVert w_{q+1}^{(p)} + w_{q+1}^{(c)} \rVert_{C_{t}W_{x}^{2m-1, p}} \overset{\eqref{[Equation (4.29), HZZ19]} \eqref{[Equation (B.7), HZZ19]} \eqref{[Equation (4.46), HZZ19]} }{\lesssim}& M_{0}(t)^{\frac{1}{2}} r_{\bot}^{\frac{2}{p} -1} r_{\lVert}^{\frac{1}{p} - \frac{1}{2}} l^{-2} \lambda_{q+1}^{2m-1}, \label{estimate 29}\\
 \lVert w_{q+1}^{(t)} \rVert_{C_{t}W_{x}^{2m-1, p}} \overset{\eqref{[Equation (4.29), HZZ19]} \eqref{[Equation (B.7), HZZ19]} \eqref{[Equation (4.47), HZZ19]} }{\lesssim}& \mu^{-1} M_{0}(t) r_{\bot}^{\frac{2}{p} - 2} r_{\lVert}^{\frac{1}{p} -1} l^{-4} \lambda_{q+1}^{2m-1}. \label{estimate 30}
\end{align}
\end{subequations} 

\subsubsection{Reynolds stress}
We see from \eqref{[Equation (4.1), HZZ19]}, \eqref{[Equation (4.19), HZZ19]} and \eqref{[Equation (4.35), HZZ19]} that 
\begin{align}
& \text{div} \mathring{R}_{q+1} - \nabla \pi_{q+1}  \label{estimate 28} \\
=& \underbrace{(-\Delta)^{m} w_{q+1} + \partial_{t} (w_{q+1}^{(p)} + w_{q+1}^{(c)}) + \text{div} ((v_{l}+  z_{l}) \otimes w_{q+1} + w_{q+1} \otimes (v_{l} + z_{l}))}_{\text{div}(R_{\text{linear}}) + \nabla \pi_{\text{linear}}}  \nonumber \\
&+ \underbrace{\text{div} ((w_{q+1}^{(c)} + w_{q+1}^{(t)}) \otimes w_{q+1} + w_{q+1}^{(p)} \otimes (w_{q+1}^{(c)} + w_{q+1}^{(t)}))}_{\text{div} (R_{\text{corrector}}) + \nabla \pi_{\text{corrector}}} + \underbrace{\text{div} (w_{q+1}^{(p)} \otimes w_{q+1}^{(p)} + \mathring{R}_{l})  + \partial_{t} w_{q+1}^{(t)}}_{\text{div}(R_{\text{oscillation}}) + \nabla \pi_{\text{oscillation}}}  \nonumber \\
&+ \underbrace{\text{div} (v_{q+1} \otimes z - v_{q+1} \otimes z_{l} + z \otimes v_{q+1} - z_{l} \otimes v_{q+1} + z \otimes z - z_{l} \otimes z_{l} )}_{\text{div}(R_{\text{commutator 2}}) + \nabla \pi_{\text{commutator 2}}} + \text{div}(R_{\text{commutator 1}}) - \nabla \pi_{l} \nonumber 
\end{align} 
with $R_{\text{commutator 1}}$ and $\pi_{l}$ from \eqref{[Equation (4.19a), HZZ19]} and 
\begin{subequations}\label{estimate 20}
\begin{align}
 R_{\text{linear}} \triangleq& \mathcal{R} ( -\Delta)^{m} w_{q+1} + \mathcal{R}\partial_{t} (w_{q+1}^{(p)} + w_{q+1}^{(c)}) + (v_{l} + z_{l}) \mathring{\otimes} w_{q+1} + w_{q+1} \mathring{\otimes} (v_{l} + z_{l}), \label{estimate 21}\\
 \pi_{\text{linear}} \triangleq& (\frac{2}{3})(v_{l} + z_{l}) \cdot w_{q+1}, \label{estimate 22}\\
 R_{\text{corrector}} \triangleq& (w_{q+1}^{(c)} + w_{q+1}^{(t)}) \mathring{\otimes} w_{q+1} + w_{q+1}^{(p)} \mathring{\otimes} (w_{q+1}^{(c)} + w_{q+1}^{(t)}), \label{estimate 23}\\
 \pi_{\text{corrector}} \triangleq& \frac{1}{3} [ (w_{q+1}^{(c)} + w_{q+1}^{(t)}) \cdot w_{q+1} + w_{q+1}^{(p)} \cdot (w_{q+1}^{(c)} + w_{q+1}^{(t)}) ], \label{estimate 24}\\
R_{\text{commutator 2}} \triangleq& v_{q+1} \mathring{\otimes} (z-z_{l}) + (z-z_{l}) \mathring{\otimes} v_{q+1} + (z-z_{l}) \mathring{\otimes} z + z_{l} \mathring{\otimes} (z-z_{l}), \label{estimate 25}\\
\pi_{\text{commutator 2}} \triangleq& \frac{1}{3} [2 v_{q+1} \cdot (z- z_{l}) + \lvert z \rvert^{2} - \lvert z_{l} \rvert^{2}], \\
R_{\text{oscillation}} \triangleq& R_{\text{oscillation}}^{(x)} + R_{\text{oscillation}}^{(t)} \text{ where } \label{estimate 26}\\
R_{\text{oscillation}}^{(x)} \triangleq& \sum_{\xi \in \Lambda} \mathcal{R}(\nabla a_{(\xi)}^{2} \mathbb{P}_{\neq 0} (W_{(\xi)} \otimes W_{(\xi)})), \hspace{2mm} R_{\text{oscillation}}^{(t)} \triangleq - \mu^{-1} \sum_{\xi \in \Lambda} \mathbb{P}_{\neq 0} (\partial_{t} a_{(\xi)}^{2} (\phi_{(\xi)}^{2} \psi_{(\xi)}^{2} \xi )), \nonumber\\
\pi_{\text{oscillation}} \triangleq& \rho + \pi_{1} \label{estimate 27} 
\end{align}
\end{subequations}
with $\pi_{1}$ from \eqref{[Equation (4.34), HZZ19]} where we used the identity 
\begin{align*}
& \text{div} (w_{q+1}^{(p)} \otimes w_{q+1}^{(p)} + \mathring{R}_{l}) + \partial_{t}w_{q+1}^{(t)} \\
\overset{\eqref{[Equation (4.31), HZZ19]} \eqref{[Equation (4.34), HZZ19]}}{=}& \sum_{\xi \in \Lambda} \mathbb{P}_{\neq 0} ( \nabla a_{(\xi)}^{2} \mathbb{P}_{\neq 0} (W_{(\xi)} \otimes W_{(\xi)})) +\nabla \rho + \nabla \pi_{1} - \mu^{-1} \sum_{\xi \in \Lambda} \mathbb{P}_{\neq 0} (\partial_{t} a_{(\xi)}^{2} (\phi_{(\xi)}^{2} \psi_{(\xi)}^{2} \xi))  
\end{align*} 
and that $\mathcal{R} = \mathcal{R} \mathbb{P}_{\neq 0}$ in \eqref{estimate 26}-\eqref{estimate 27}. From \eqref{estimate 28} we read that $\pi_{q+1} \triangleq \pi_{l} - \pi_{\text{linear}} - \pi_{\text{corrector}} - \pi_{\text{oscillation}} - \pi_{\text{commutator 2}}$ while 
\begin{equation}\label{[Equation (4.48e), HZZ19]}
\mathring{R}_{q+1} \triangleq R_{\text{linear}} + R_{\text{corrector}} + R_{\text{oscillation}} + R_{\text{commutator 2}} + R_{\text{commutator 1}}. 
\end{equation}
Now we shall fix 
\begin{equation}\label{p}
p^{\ast} \triangleq \frac{40m - 14}{170\alpha - 19 + 44m}
\end{equation}
which may be immediately verified to be in the range of $(1,2)$ using \eqref{alpha}. For $t \in [0, T_{L}]$ we may estimate from \eqref{estimate 21} using \eqref{[Equation (4.3), HZZ19]}, \eqref{[Equation (4.32a), HZZ19]} and  \eqref{[Equation (4.4b), HZZ19]}
\begin{align}\label{estimate 32}
\lVert R_{\text{linear}} \rVert_{C_{t}L_{x}^{p^{\ast}}} \lesssim \lVert w_{q+1} \rVert_{C_{t}W_{x}^{2m-1, p^{\ast}}} + \sum_{\xi \in \Lambda} \lVert \partial_{t} \text{curl} (a_{(\xi)} V_{(\xi)}) \rVert_{C_{t}L_{x}^{p^{\ast}}} + M_{0}(t)^{\frac{1}{2}}  \lVert w_{q+1} \rVert_{C_{t}L_{x}^{p^{\ast}}}. 
\end{align} 
First we can estimate from \eqref{[Equation (4.35), HZZ19]} 
\begin{equation}\label{estimate 41}
\lVert w_{q+1} \rVert_{C_{t}W_{x}^{2m-1, p^{\ast}}} \overset{\eqref{estimate 29} \eqref{estimate 30}}{\lesssim} M_{0}(t)^{\frac{1}{2}} \lambda_{q+1}^{-\frac{61\alpha}{6}} + M_{0}(t) \lambda_{q+1}^{\frac{12m - 15 - 148\alpha}{24}}.  
\end{equation} 
Second we can estimate 
\begin{equation}
\sum_{\xi \in \Lambda} \lVert \partial_{t} \text{curl} (a_{(\xi)} V_{(\xi)} ) \rVert_{C_{t}L_{x}^{p^{\ast}}} 
\overset{\eqref{[Equation (4.29), HZZ19]} \eqref{[Equation (B.7), HZZ19]} \eqref{[Equation (4.16), HZZ19]}}{\lesssim} M_{0}(t)^{\frac{1}{2}} l^{-7} r_{\bot}^{\frac{2}{p^{\ast}}} r_{\lVert}^{\frac{1}{p^{\ast}} - \frac{3}{2}} \mu  + M_{0}(t)^{\frac{1}{2}} l^{-10} r_{\bot}^{\frac{2}{p^{\ast}} -1} r_{\lVert}^{\frac{1}{p^{\ast}} - \frac{1}{2}} \lambda_{q+1}^{-1}. \label{estimate 31}
\end{equation} 
Now we can apply \eqref{estimate 29}, \eqref{estimate 30}, \eqref{estimate 41} and  \eqref{estimate 31} to \eqref{estimate 32} and use our choice  \eqref{p} to deduce 
\begin{align}
&\lVert R_{\text{linear}} \rVert_{C_{t}L_{x}^{p^{\ast}}} \nonumber\\
\overset{\eqref{[Equation (4.38), HZZ19]} \eqref{[Equation (4.41), HZZ19]}}{\lesssim}& M_{0}(t)^{\frac{1}{2}} \lambda_{q+1}^{-\frac{61\alpha}{6}} + M_{0}(t) \lambda_{q+1}^{\frac{12m - 15 - 148\alpha}{24}}  + M_{0}(t)^{\frac{1}{2}} l^{-7} r_{\bot}^{\frac{2}{p^{\ast}}} r_{\lVert}^{\frac{1}{p^{\ast}} - \frac{3}{2}} \mu \nonumber\\
& + M_{0}(t)^{\frac{1}{2}} l^{-10} r_{\bot}^{\frac{2}{p^{\ast}} -1} r_{\lVert}^{\frac{1}{p^{\ast}} - \frac{1}{2}} \lambda_{q+1}^{-1} + M_{0}(t)^{\frac{1}{2}}  [M_{0}(t)^{\frac{1}{2}} \delta_{q+1}^{\frac{1}{2}} l^{-2} r_{\bot}^{\frac{2}{p^{\ast}} -1} r_{\lVert}^{\frac{1}{p^{\ast}} - \frac{1}{2}}]\nonumber\\
\overset{\eqref{[Equation (4.16), HZZ19]} \eqref{[Equation (4.23), HZZ19]} }{\lesssim}& M_{0}(t)^{\frac{1}{2}} \lambda_{q+1}^{-\frac{61\alpha}{6}} + M_{0}(t) \lambda_{q+1}^{\frac{12m - 15 - 148\alpha}{24}}   + M_{0}(t)^{\frac{1}{2}} \lambda_{q+1}^{14\alpha} \lambda_{q+1}^{ (\frac{1-20m}{24}) \frac{2}{p^{\ast}}} (\lambda_{q+1}^{\frac{13-20m}{12}})^{\frac{1}{p^{\ast}} - \frac{3}{2}} \lambda_{q+1}^{2m-1} \lambda_{q+1}^{\frac{25-20m}{24}} \nonumber\\
&+ M_{0}(t)^{\frac{1}{2}} \lambda_{q+1}^{20\alpha} (\lambda_{q+1}^{\frac{1-20m}{24}})^{\frac{2}{p^{\ast}} -1} (\lambda_{q+1}^{\frac{13-20m}{12}})^{\frac{1}{p^{\ast}} - \frac{1}{2}} \lambda_{q+1}^{-1} + M_{0}(t)  \lambda_{q+1}^{4\alpha} (\lambda_{q+1}^{\frac{1-20m}{24}})^{\frac{2}{p^{\ast}} -1} (\lambda_{q+1}^{\frac{13-20m}{12}})^{\frac{1}{p^{\ast}} - \frac{1}{2}} \nonumber \\
\lesssim& M_{0}(t)^{\frac{1}{2}} [\lambda_{q+1}^{-\frac{61\alpha}{6}} + \lambda_{q+1}^{-\frac{\alpha}{6}} + \lambda_{q+1}^{\frac{35\alpha - 12m}{6}}] + M_{0}(t) [\lambda_{q+1}^{\frac{12m - 15 - 148\alpha}{24}} + \lambda_{q+1}^{\frac{ - 122\alpha + 12 - 24m}{12}}] \label{estimate 46}
\end{align} 
where we used the fact that $\frac{\alpha}{4} > \frac{4}{b}$. Therefore, taking $a > 0$ sufficiently large and $\beta > 0$ sufficiently small such that $16 \beta b < \alpha$ which implies $\lambda_{q+1}^{-\frac{\alpha}{6}} \lesssim \delta_{q+2} \lambda_{q+1}^{-\frac{\alpha}{6} + \frac{\alpha}{8}} \ll 1$ leads to 
\begin{equation}\label{[Equation (4.48h), HZZ19]}
\lVert R_{\text{linear}} \rVert_{C_{t}L_{x}^{p^{\ast}}} \leq (2\pi)^{-3 (\frac{p^{\ast}-1}{p^{\ast}})} \frac{M_{0}(t) c_{R} \delta_{q+2}}{5}.
\end{equation} 
Next, for all $t \in [0, T_{L} ]$, applying H$\ddot{\mathrm{o}}$lder's inequality on \eqref{estimate 23} gives us  
\begin{align}
 \lVert R_{\text{corrector}} \rVert_{C_{t}L^{p^{\ast}}} \lesssim& (\lVert w_{q+1}^{(c)} \rVert_{C_{t}L_{x}^{2p^{\ast}}} + \lVert w_{q+1}^{(t)} \rVert_{C_{t}L_{x}^{2p^{\ast}}}) (\lVert w_{q+1}^{(c)} \rVert_{C_{t}L_{x}^{2p^{\ast}}} + \lVert w_{q+1}^{(t)} \rVert_{C_{t}L_{x}^{2p^{\ast}}} + \lVert w_{q+1}^{(p)} \rVert_{C_{t}L_{x}^{2p^{\ast}}}) \nonumber\\
\overset{\eqref{estimate 33} \eqref{[Equation (4.41), HZZ19]} \eqref{[Equation (4.16), HZZ19]}}{\lesssim}& (M_{0}(t)^{\frac{1}{2}} \delta_{q+1}^{\frac{1}{2}} \lambda_{q+1}^{\frac{8m-10 + 203\alpha}{12}} + M_{0}(t) \delta_{q+1}\lambda_{q+1}^{\frac{11 \alpha + 4m-5}{12}}) M_{0}(t)^{\frac{1}{2}} \delta_{q+1}^{\frac{1}{2}} \lambda_{q+1}^{\frac{ - 74 \alpha - 4m +5}{24}}.
\end{align} 
Now using the fact that $\delta_{q+2}^{-1} \ll  \lambda_{q+1}^{\frac{\alpha}{8}}$ due to $16\beta b < \alpha$, that  $\alpha  <  \frac{15 - 12m}{335}$ by our choice from \eqref{alpha}, we can take $a > 0$ sufficiently large and $\beta > 0$ sufficiently small to obtain 
\begin{align}
&\lVert R_{\text{corrector}} \rVert_{C_{t}L_{x}^{p^{\ast}}} \label{[Equation (4.48i), HZZ19]}\\
\leq& M_{0}(t) \delta_{q+2} ( C \lambda_{q+1}^{\frac{\alpha}{8}} \lambda_{q+1}^{\frac{12 m -15 + 332 \alpha}{24}} + C M_{0}(t)^{\frac{1}{2}} \lambda_{q+1}^{\frac{\alpha}{8}} \lambda_{q+1}^{\frac{ 4m - 5 - 52 \alpha}{24}}) \overset{\eqref{[Equation (4.18), HZZ19]}}{\leq} (2\pi)^{-3 (\frac{p^{\ast}-1}{p^{\ast}})} \frac{M_{0}(t) \delta_{q+2} c_{R}}{5}. \nonumber
\end{align} 
Next, because $W_{(\xi)}$ is $(\mathbb{T}/r_{\bot}\lambda_{q+1})^{3}$-periodic, the minimal active frequency in $\mathbb{P}_{\neq 0} (W_{(\xi)} \otimes W_{(\xi)})$ within $R_{\text{oscillation}}^{(x)}$ of \eqref{estimate 26} is given by $r_{\bot}\lambda_{q+1}$, which implies that 
\begin{equation}\label{[Equation (4.48j), HZZ19]}
\mathbb{P}_{\neq 0} (W_{(\xi)} \otimes W_{(\xi)} ) = \mathbb{P}_{\geq \frac{r_{\bot}\lambda_{q+1}}{2}}(W_{(\xi)} \otimes W_{(\xi)}). 
\end{equation} 
Because it follows immediately from \eqref{[Equation (4.29), HZZ19]} that $\lVert D^{j} \nabla a_{(\xi)}^{2} \rVert_{C_{x}} \leq (l^{-9} M_{0}(t))l^{-5j}$, we can apply Lemma \ref{[Lemma 7.5, BV19b]} with ``$a$'' = $a_{(\xi)}^{2}$, ``$C_{a}$'' = $l^{-9} M_{0}(t)$, ``$\zeta$'' = $l^{-5}$ and ``$\kappa$'' = $r_{\bot}\lambda_{q+1}$. Here, we note that one of the hypothesis of Lemma \ref{[Lemma 7.5, BV19b]} requires that $l^{-5N} \leq (r_{\bot}\lambda_{q+1})^{N-2} = (\lambda_{q+1}^{\frac{25 - 20m}{24}})^{N-2}$, which may be verified for $N \geq 3$ large using \eqref{[Equation (4.16), HZZ19]} because our choice of $\alpha$ from \eqref{alpha} satisfies $\alpha < \frac{5-4m}{48}$. Thus, we may now deduce directly from \eqref{estimate 26} 
\begin{align}
\lVert R_{\text{oscillation}}^{(x)} \rVert_{C_{t}L_{x}^{p^{\ast}}} \overset{\eqref{[Equation (4.48j), HZZ19]}}{\lesssim}& M_{0}(t) l^{-9} r_{\bot}^{-1} \lambda_{q+1}^{-1} \lVert W_{(\xi)} \otimes W_{(\xi)} \rVert_{C_{t}L_{x}^{p^{\ast}}} \overset{ \eqref{[Equation (4.23), HZZ19]} \eqref{p}\eqref{[Equation (B.7), HZZ19]}  }{\lesssim} M_{0}(t) \lambda_{q+1}^{\frac{92\alpha + 12m - 15}{24}}.
\end{align} 
Because $\delta_{q+2}^{-1} \lambda_{q+1}^{\frac{92 \alpha + 12m -15}{24}} < \lambda_{q+1}^{\frac{95 \alpha + 12m -15}{24}}$ due to the assumption that $16\beta b < \alpha$ and that our choice of $\alpha$ from \eqref{alpha} satisfies $\alpha < \frac{15-12m}{95}$, we deduce 
\begin{align}\label{[Equation (4.48m), HZZ19]}
\lVert R_{\text{oscillation}}^{(x)} \rVert_{C_{t}L_{x}^{p^{\ast}}} \lesssim \delta_{q+2}M_{0}(t) \lambda_{q+1}^{\frac{95 \alpha + 12 m -15}{24}} \leq (2\pi)^{-3 (\frac{p^{\ast}-1}{p^{\ast}})} \frac{c_{R} \delta_{q+2}M_{0}(t)}{10}.
\end{align} 
Next we directly estimate from \eqref{estimate 26} using $\lVert \mathcal{R} \rVert_{L^{p} \mapsto L^{p}} \lesssim 1$ for $p \in (1,\infty)$ from Lemma \ref{divergence inverse operator}  
\begin{align}
\lVert R_{\text{oscillation}}^{(t)} \rVert_{C_{t}L_{x}^{p^{\ast}}}\lesssim&  \lambda_{q+1}^{ \frac{ - 1 - 28m}{24}} \sum_{\xi \in \Lambda} \lVert a_{(\xi)} \rVert_{C_{t}C_{x}} \lVert \partial_{t} a_{(\xi)} \rVert_{C_{t}C_{x}} \lVert \phi_{(\xi)} \rVert_{C_{t}L_{x}^{2p^{\ast}}}^{2} \lVert \psi_{(\xi)} \rVert_{C_{t}L_{x}^{2p^{\ast}}}^{2} \nonumber\\
\overset{\eqref{[Equation (4.16), HZZ19]} \eqref{[Equation (4.29), HZZ19]} \eqref{[Equation (B.7), HZZ19]}}{\lesssim}& \delta_{q+1} M_{0}(t) \lambda_{q+1}^{\frac{ 9 - 36 m - 4\alpha}{24}} \leq (2\pi)^{-3 (\frac{p^{\ast}-1}{p^{\ast}})} \frac{M_{0}(t) c_{R} \delta_{q+2}}{10} \label{[Equation (4.48n), HZZ19]}
\end{align} 
where we used that $\phi_{(\xi)}$ and $\psi_{(\xi)}$ are functions of variables in orthogonal directions in the first inequality. Next we estimate from \eqref{[Equation (4.19a), HZZ19]} for all $t \in [0, T_{L}]$
\begin{align}
 \lVert R_{\text{commutator 1}} \rVert_{C_{t}L_{x}^{1}} \lesssim& l (\lVert v_{q} \rVert_{C_{t,x}^{1}} + \lVert z \rVert_{C_{t}C_{x}^{1}}) ( \lVert v_{q} \rVert_{C_{t}L_{x}^{2}} + \lVert z \rVert_{C_{t,x}}) \nonumber\\
&+ l^{\frac{2}{5} - 2\delta} \lVert z \rVert_{C_{t}^{\frac{2}{5} - 2 \delta}C_{x}}  (\lVert v_{q} \rVert_{C_{t}L_{x}^{2}} + \lVert z \rVert_{C_{t,x}}) \lesssim l^{\frac{2}{5} - 2 \delta} M_{0}(t) \lambda_{q}^{4} \nonumber\\
\overset{ \eqref{[Equation (4.3), HZZ19]}  \eqref{[Equation (4.4a), HZZ19]}\eqref{[Equation (4.4b), HZZ19]} }{\lesssim}& \delta_{q+2} M_{0}(t) a^{b^{q} [ - \frac{\alpha b}{2} + \frac{10}{3} + b^{2} 2 \beta]} \lesssim  \delta_{q+2} M_{0}(t) a^{b^{q} (-\frac{8}{3})} \leq \frac{M_{0} (t) c_{R} \delta_{q+2}}{5} \label{[Equation (4.48o), HZZ19]}
\end{align} 
where we used that $-\frac{\alpha b}{2} + \frac{10}{3} + b^{2} 2 \beta < -\frac{8}{3}$ due to our assumptions of $\alpha > 16 \beta b$ and $\alpha b> 16$. Finally, we estimate for all $t \in [0, T_{L}]$ using \eqref{[Equation (4.4a), HZZ19]} at level $q + 1$ which we already verified 
 \begin{align}\label{[Equation (4.48p), HZZ19]}
\lVert R_{\text{commutator 2}} \rVert_{C_{t}L_{x}^{1}} \overset{\eqref{estimate 25} \eqref{[Equation (4.3), HZZ19]}}{\lesssim}&  [ \lVert v_{q+1} \rVert_{C_{t}L_{x}^{1}} + \lVert z \rVert_{C_{t}L_{x}^{1}}] l^{\frac{2}{5} - 2 \delta} ( \lVert z \rVert_{C_{t}^{\frac{2}{5} - 2\delta} L_{x}^{\infty}} + \lVert z \rVert_{C_{t}C_{x}^{\frac{2}{5} - 2 \delta}}) \nonumber\\
\lesssim& M_{0}(t) l^{\frac{2}{5} - 2\delta} \leq \frac{M_{0}(t) c_{R} \delta_{q+2}}{5} 
\end{align} 
where the last inequality used the fact that $l^{\frac{2}{5} - 2\delta} \lambda_{q}^{4} \ll \frac{c_{R} \delta_{q+2}}{5}$ as we showed in \eqref{[Equation (4.48o), HZZ19]}. At last, we may now readily deduce \eqref{[Equation (4.4c), HZZ19]} at level $q + 1$ using H$\ddot{\mathrm{o}}$lder's inequality as 
\begin{align}
&\lVert\mathring{R}_{q+1} \rVert_{C_{t}L_{x}^{1}} \overset{\eqref{[Equation (4.48e), HZZ19]} \eqref{[Equation (4.48o), HZZ19]} \eqref{[Equation (4.48p), HZZ19]}}{\leq} (2\pi)^{3 (\frac{p^{\ast}-1}{p^{\ast}})} [ \lVert R_{\text{linear}} \rVert_{C_{t}L_{x}^{p^{\ast}}} + \lVert R_{\text{corrector}} \rVert_{C_{t}L_{x}^{p^{\ast}}} + \lVert R_{\text{oscillation}} \rVert_{C_{t}L_{x}^{p^{\ast}}} ] \nonumber\\
& \hspace{30mm} + \frac{2M_{0}(t) c_{R} \delta_{q+2}}{5}  \overset{\eqref{[Equation (4.48h), HZZ19]} \eqref{[Equation (4.48i), HZZ19]} \eqref{[Equation (4.48m), HZZ19]} \eqref{[Equation (4.48n), HZZ19]}}{\leq} M_{0}(t) c_{R} \delta_{q+2}.  \label{[Equation (4.48q), HZZ19]}
\end{align} 

We conclude by commenting on how $(v_{q+1}, \mathring{R}_{q+1})$ is $(\mathcal{F}_{t})_{t\geq 0}$-adapted and that $v_{q+1}(0,x)$, $\mathring{R}_{q+1}(0,x)$ are both deterministic if $v_{q}(0,x), \mathring{R}_{q}(0,x)$ are deterministic. First, recall that $z$ in \eqref{[Equation (3.10), HZZ19]} is $(\mathcal{F}_{t})_{t\geq 0}$-adapted; due to the compact support of $\varphi_{l}$ in $\mathbb{R}_{+}$, it follows that $z_{l}$ in \eqref{[Equation (4.18c), HZZ19]} is $(\mathcal{F}_{t})_{t \geq 0}$-adapted. Similarly, because $(v_{q}, \mathring{R}_{q})$ are both $(\mathcal{F}_{t})_{t\geq 0}$-adapted by hypothesis, it follows that $(v_{l}, \mathring{R}_{l})$ in \eqref{[Equation (4.18c), HZZ19]} are both $(\mathcal{F}_{t})_{t\geq 0}$-adapted. Because $M_{0}(t)$ is deterministic, it follows that $\rho$ in \eqref{[Equation (4.23b), HZZ19]} is also $(\mathcal{F}_{t})_{t \geq 0}$-adapted. Due to $\rho$ and $\mathring{R}_{l}$ being $(\mathcal{F}_{t})_{t \geq 0}$-adapted, $a_{(\xi)}$ in \eqref{[Equation (4.26), HZZ19]} is also $(\mathcal{F}_{t})_{t \geq 0}$-adapted. Because $W_{(\xi)}$ from \eqref{[Equation (B.3), HZZ19]} is deterministic, this implies that $w_{q+1}^{(p)}$ in \eqref{[Equation (4.30), HZZ19]} is also $(\mathcal{F}_{t})_{t\geq 0}$-adapted. Because $V_{(\xi)}$ and $W_{(\xi)}^{(c)}$ in \eqref{[Equation (B.6), HZZ19]} are both deterministic, we see that $w_{q+1}^{(c)}$ from \eqref{[Equation (4.32), HZZ19]} is also $(\mathcal{F}_{t})_{t\geq 0}$-adapted. Moreover, because $\phi_{(\xi)}$ and $\psi_{(\xi)}$ from \eqref{[Equation (B.2c), HZZ19]} are both deterministic, it follows that $w_{q+1}^{(t)}$ is also $(\mathcal{F}_{t})_{t\geq 0}$-adapted. Consequently, $w_{q+1}$ defined in \eqref{[Equation (4.35), HZZ19]} is $(\mathcal{F}_{t})_{t \geq 0}$-adapted, which in turn implies that $v_{q+1}$ in \eqref{[Equation (4.35), HZZ19]} is also $(\mathcal{F}_{t})_{t\geq 0}$-adapted. 

It is also clear from the compact support of $\varphi_{l}$ in $\mathbb{R}_{+}$ that if $v_{q}(0,x)$ and $\mathring{R}_{q}(0,x)$ are deterministic, then so are $v_{l}(0,x), \mathring{R}_{l}(0,x)$ and $\partial_{t} \mathring{R}_{l}(0,x)$. Because $z(0,x) = 0$ from \eqref{[Equation (3.10), HZZ19]}, $R_{\text{commutator 1}}(0,x)$ in \eqref{[Equation (4.19a), HZZ19]} is also deterministic. In turn, this implies that $\rho(0,x)$ and $\partial_{t} \rho(0,x)$ from \eqref{[Equation (4.23b), HZZ19]} are also deterministic as $M_{0}(t)$ is deterministic, which leads to $a_{(\xi)}(0,x)$ and $\partial_{t} a_{(\xi)}(0,x)$ in \eqref{[Equation (4.26), HZZ19]} also being deterministic. Because $W_{(\xi)}$ from \eqref{[Equation (B.3), HZZ19]}, $\phi_{(\xi)}, \psi_{(\xi)}$ from \eqref{[Equation (B.2c), HZZ19]} are all deterministic, $R_{\text{oscillation}}(0,x)$ is also deterministic.  This implies that $w_{q+1}^{(p)} (0,x)$ and $\partial_{t} w_{q+1}^{(p)}(0,x)$ from \eqref{[Equation (4.30), HZZ19]} are also deterministic because $w_{q+1}^{(c)}(0,x)$ and $\partial_{t} w_{q+1}^{(c)}(0,x)$ from \eqref{[Equation (4.32), HZZ19]} are also deterministic given that $V_{(\xi)}$ and $W_{(\xi)}^{(c)}$ from \eqref{[Equation (B.6), HZZ19]} are both deterministic. Similarly, $a_{(\xi)}(0,x)$ being deterministic leads to $w_{q+1}^{(t)}(0,x)$ from \eqref{[Equation (4.33), HZZ19]} also being deterministic. Consequently, $w_{q+1}(0,x)$ in \eqref{[Equation (4.35), HZZ19]} is also deterministic, which in turn implies that $v_{q+1}(0,x)$ in \eqref{[Equation (4.35), HZZ19]} would also be deterministic. It follows from \eqref{estimate 20} that $R_{\text{linear}}(0,x), R_{\text{corrector}}(0,x)$ and $R_{\text{commutator 2}}(0,x)$ are all deterministic, leading us to conclude that $\mathring{R}_{q+1}(0,x)$ from \eqref{[Equation (4.48e), HZZ19]} is deterministic. This completes the proof of Proposition \ref{[Proposition 4.2, HZZ19]}.  

\section{Proof in the case of linear multiplicative noise}

\subsection{Proof of Theorem \ref{[Theorem 1.4, HZZ19]} assuming Theorem \ref{[Theorem 1.3, HZZ19]} }

From \eqref{[Equation (3.12), HZZ19]}, \eqref{[Equation (3.13), HZZ19]}, \eqref{[Equation (3.16), HZZ19]} and \eqref{[Equation (4.2), HZZ19]}, it can be understood that the stopping time in \eqref{[Equation (4.2), HZZ19]} is a function of $u$. In the case of the multiplicative noise, one will not be able to split \eqref{[Equation (1.7), HZZ19]} as we did in \eqref{[Equation (3.10), HZZ19]}-\eqref{[Equation (3.11), HZZ19]}, and the stopping time will become a function of the noise $B$. Thus, in this case  we work with a different definition of a solution, called the probabilistically weak solution. We recall $U_{1}$, $\bar{\Omega}$ and $\bar{\mathcal{B}}_{t}$ defined in section \ref{Preliminaries}, and fix any $\gamma \in (0,1)$ for the purpose of the following definition. 
\begin{define}\label{[Definition 5.1, HZZ19]}
Let $s \geq 0$ and $\xi^{\text{in}} \in L_{\sigma}^{2}$, $\theta^{\text{in}} \in U_{1}$. A probability measure $P \in \mathcal{P} (\bar{\Omega})$ is a probabilistically weak solution to \eqref{[Equation (1.1), HZZ19]} with initial condition $(\xi^{\text{in}}, \theta^{\text{in}})$ at initial time $s$ if  
\begin{itemize}
\item [] (M1) $P(\{ \xi(t) = \xi^{\text{in}}, \theta(t) = \theta^{\text{in}} \hspace{1mm} \forall \hspace{1mm} t \in [0,s]\}) = 1$ and for all $n \in \mathbb{N}$ 
\begin{equation}\label{[Equation (4.48r), HZZ19]}
P ( \{ (\xi, \theta) \in \bar{\Omega}: \int_{0}^{n} \lVert G(\xi(r)) \rVert_{L_{2} (U, L_{\sigma}^{2})}^{2} dr < \infty \} ) = 1, 
\end{equation} 
\item [] (M2) under $P$, $\theta$ is a cylindrical $(\bar{\mathcal{B}}_{t})_{t\geq s}$-Wiener process on $U$ starting from initial condition $\theta^{\text{in}}$ at initial time $s$ and for every $e_{i} \in C^{\infty} (\mathbb{T}^{3}) \cap L_{\sigma}^{2}$ and $t\geq s$, 
\begin{equation}\label{[Equation (4.48s), HZZ19]}
\langle \xi(t) - \xi(s), e_{i} \rangle + \int_{s}^{t} \langle \text{div} (\xi(r) \otimes \xi(r)) + (-\Delta)^{m} \xi(r), e_{i} \rangle dr = \int_{s}^{t} \langle e_{i}, G(\xi(r)) d\theta(r) \rangle, 
\end{equation} 
\item [] (M3) for any $q \in \mathbb{N}$, there exists a function $t \mapsto C_{t,q} \in \mathbb{R}_{+}$ for all $t \geq s$ such that 
\begin{equation}\label{[Equation (4.48t), HZZ19]}
\mathbb{E}^{P} [ \sup_{r \in [0,t]} \lVert \xi(r) \rVert_{L_{x}^{2}}^{2q} + \int_{s}^{t} \lVert \xi(r) \rVert_{H_{x}^{\gamma}}^{2} dr] \leq C_{t,q} (1+ \lVert \xi^{\text{in}} \rVert_{L_{x}^{2}}^{2q}). 
\end{equation}
\end{itemize}
The set of all such probabilistically weak solutions with the same constant $C_{t,q}$ in \eqref{[Equation (4.48t), HZZ19]} for every $q \in \mathbb{N}$ and $t \geq s$ is denoted by $\mathcal{E}(s, \xi^{\text{in}}, \theta^{\text{in}}, \{C_{t,q}\}_{q\in\mathbb{N}, t \geq s})$. 
\end{define} 

For a stopping time $\tau$, we set 
\begin{equation}
\bar{\Omega}_{\tau} \triangleq \{ \omega ( \cdot \wedge \tau(\omega)): \hspace{0.5mm} \omega \in \bar{\Omega} \} 
\end{equation} 
and denote by ($\bar{\mathcal{B}}_{\tau})$ the $\sigma$-field associated to $\tau$. 

\begin{define}\label{[Definition 5.2, HZZ19]}
Let $s \geq 0$ and $\xi^{\text{in}} \in L_{\sigma}^{2}, \theta^{\text{in}} \in U_{1}$. Let $\tau \geq s$ be a stopping time of $(\bar{\mathcal{B}}_{t})_{t \geq s}$. A probability measure $P \in \mathcal{P} (\bar{\Omega}_{\tau})$ is a probabilistically weak solution to \eqref{[Equation (1.1), HZZ19]} on $[s, \tau]$ with initial condition $(\xi^{\text{in}}, \theta^{\text{in}})$ at initial time $s$ if 
\begin{itemize}
\item [] (M1) $P(\{ \xi(t) = \xi^{\text{in}}, \theta(t) = \theta^{\text{in}} \hspace{1mm} \forall \hspace{1mm} t \in [0,s] \}) = 1$ and for all $n \in \mathbb{N}$
\begin{equation}\label{[Equation (4.48v), HZZ19]}
P ( \{ ( \xi, \theta) \in \bar{\Omega}: \int_{0}^{n \wedge \tau} \lVert G(\xi(r)) \rVert_{L_{2}(U; L_{\sigma}^{2})}^{2} dr < \infty \}) = 1, 
\end{equation}
\item  [] (M2) under $P$, $\langle \theta(\cdot \wedge \tau), l_{i} \rangle_{U}$, where $\{l_{i} \}_{i\in \mathbb{N}}$ is an orthonormal basis of $U$,  is a continuous, square-integrable martingale w.r.t. $(\bar{\mathcal{B}}_{t})_{t \geq s}$ with initial condition $\langle \theta^{\text{in}}, l_{i} \rangle_{U}$ at initial time $s$ with a quadratic variation process given by $(t \wedge \tau - s) \lVert l_{i} \rVert_{U}^{2}$ and for every $e_{i} \in C^{\infty} (\mathbb{T}^{3} ) \cap L_{\sigma}^{2}$ and $t \geq s$, 
\begin{equation}\label{[Equation (4.48w), HZZ19]}
\langle \xi(t\wedge \tau) - \xi(s), e_{i} \rangle + \int_{s}^{t \wedge \tau} \langle \text{div} (\xi(r) \otimes \xi(r)) + (-\Delta)^{m} \xi(r), e_{i} \rangle dr = \int_{s}^{t \wedge \tau} \langle e_{i}, G(\xi(r)) d\theta(r) \rangle, 
\end{equation}
\item [] (M3) for any $q \in \mathbb{N}$, there exists a function $t\mapsto C_{t,q} \in \mathbb{R}_{+}$ for all $t \geq s$ such that 
\begin{equation}\label{[Equation (4.48x), HZZ19]}
\mathbb{E}^{P} [ \sup_{r \in [0, t \wedge \tau]} \lVert \xi(r) \rVert_{L_{x}^{2}}^{2q} + \int_{s}^{t \wedge \tau} \lVert \xi(r) \rVert_{H_{x}^{\gamma}}^{2} dr] \leq C_{t,q} (1+ \lVert \xi^{\text{in}} \rVert_{L_{x}^{2}}^{2q}). 
\end{equation} 
\end{itemize} 
\end{define} 
As we already discussed referring to \cite[Corollary 4.9 in Chapter 5]{KS91}, the joint uniqueness in law for \eqref{[Equation (1.1), HZZ19]} is equivalent to the uniqueness of probabilistically weak solution in Definition \ref{[Definition 5.1, HZZ19]}, which holds precisely if probabilistically weak solutions starting from the same initial distributions are unique. Analogously to Proposition \ref{[Theorem 3.1, HZZ19]}, we have the following result concerning existence and stability of a probabilistically weak solution to \eqref{[Equation (1.1), HZZ19]}. 
\begin{proposition}\label{[Theorem 5.1, HZZ19]}
For every $(s, \xi^{\text{in}}, \theta^{\text{in}}) \in [0,\infty) \times L_{\sigma}^{2} \times U_{1}$, there exists a probabilistically weak solution $P \in \mathcal{P} (\bar{\Omega})$ to \eqref{[Equation (1.1), HZZ19]} with  initial condition $(\xi^{\text{in}}, \theta^{\text{in}})$ at initial time $s$ that satisfies Definition \ref{[Definition 5.1, HZZ19]}. Moreover, if there exists a family $(s_{n}, \xi_{n}, \theta_{n}) \subset [0,\infty) \times L_{\sigma}^{2} \times U_{1}$ such that $\lim_{n\to\infty} \lVert (s_{n}, \xi_{n}, \theta_{n}) - (s, \xi^{\text{in}}, \theta^{\text{in}}) \rVert_{\mathbb{R} \times L_{x}^{2} \times U_{1}} = 0$ and $P_{n} \in \mathcal{E} (s_{n}, \xi_{n}, \theta_{n}, \{C_{t,q}\}_{q\in\mathbb{N}, t \geq s_{n}})$ is the probabilistically weak solution corresponding to $(s_{n}, \xi_{n}, \theta_{n})$, then there exists a subsequence $\{P_{n_{k}} \}_{k\in\mathbb{N}}$ that converges weakly to some $P \in \mathcal{E}(s,\xi^{\text{in}}, \theta^{\text{in}}, \{C_{t,q}\}_{q\in\mathbb{N}, t \geq s})$. 
\end{proposition} 
\begin{proof}[Proof of Proposition \ref{[Theorem 5.1, HZZ19]}]
The proof of \cite[Theorem 5.1]{HZZ19} can be extended to our case with the fractional Laplacian diffusive term via a straight-forward modification as we did in the proof of Proposition \ref{[Theorem 3.1, HZZ19]}. 
\end{proof}

The proofs of the following results from \cite{HZZ19} do not rely on the specific form of the diffusive term and thus apply essentially directly to our case. 
\begin{lemma}\label{[Proposition 5.2, HZZ19]}
\rm{(\cite[Proposition 5.2]{HZZ19})} Let $\tau$ be a bounded stopping time of $(\bar{\mathcal{B}}_{t})_{t \geq 0}$. Then for every $\omega \in \bar{\Omega}$, there exists $Q_{\omega} \in \mathcal{P}(\bar{\Omega})$ such that 
\begin{subequations}
\begin{align}
& Q_{\omega} ( \{ \omega' \in \bar{\Omega}: \hspace{0.5mm}  ( \xi, \theta) (t, \omega') = (\xi, \theta) (t,\omega) \hspace{1mm} \forall \hspace{1mm} t \in [0, \tau(\omega)] \}) = 1, \\
& Q_{\omega} (A) = R_{\tau(\omega), \xi(\tau(\omega), \omega), \theta(\tau(\omega), \omega)} (A) \hspace{1mm} \forall \hspace{1mm} A \in \bar{\mathcal{B}}^{\tau(\omega)},  
\end{align} 
\end{subequations}
where $R_{\tau(\omega), \xi(\tau(\omega), \omega), \theta(\tau(\omega), \omega)} \in \mathcal{P} (\bar{\Omega})$ is a probabilistically weak solution to \eqref{[Equation (1.1), HZZ19]} with initial condition $(\xi(\tau(\omega), \omega), \theta(\tau(\omega), \omega))$ at initial time $\tau(\omega)$. Moreover, for every $A \in \bar{\mathcal{B}}$, the mapping $\omega \mapsto Q_{\omega}(A)$ is $\bar{\mathcal{B}}_{\tau}$-measurable. 
\end{lemma}

\begin{lemma}\label{[Proposition 5.3, HZZ19]}
\rm{(\cite[Proposition 5.3]{HZZ19})} Let $\xi^{\text{in}} \in L_{\sigma}^{2}$ and $P$ be a probabilistically weak solution to \eqref{[Equation (1.1), HZZ19]} on $[0,\tau]$ with initial condition $(\xi^{\text{in}}, 0)$ at initial time $0$ that satisfies Definition \ref{[Definition 5.2, HZZ19]}. In addition to the hypothesis of Lemma \ref{[Proposition 5.2, HZZ19]}, suppose that there exists a Borel set $\mathcal{D} \subset \bar{\Omega}_{\tau}$ such that $P(\mathcal{D}) = 0$ and for every $\omega \in \bar{\Omega}_{\tau} \setminus \mathcal{D}$, it satisfies 
\begin{equation}\label{[Equation (5.3), HZZ19]}
Q_{\omega} (\{ \omega' \in \bar{\Omega}: \hspace{0.5mm}  \tau(\omega') = \tau(\omega) \}) = 1. 
\end{equation} 
Then the probability measure $P\otimes_{\tau}R \in \mathcal{P} (\bar{\Omega})$ defined by 
\begin{equation}\label{[Equation (5.3a), HZZ19]}
P \otimes_{\tau} R (\cdot) \triangleq \int_{\bar{\Omega}} Q_{\omega} (\cdot) P(d \omega)
\end{equation} 
satisfies $P\otimes_{\tau}R \rvert_{\bar{\Omega}_{\tau}} = P \rvert_{\bar{\Omega}_{\tau}}$ and it is a probabilistically weak solution to \eqref{[Equation (1.1), HZZ19]} on $[0,\infty)$ with initial condition $(\xi^{\text{in}}, 0)$ at initial time $0$. 
\end{lemma}

Now we fix an $\mathbb{R}$-valued Wiener process $B$ on $(\Omega, \mathcal{F}, \mathbb{P})$. For $n \in \mathbb{N}, L > 1$ and 
\begin{equation}\label{new 16}
\delta \in (0, \frac{1}{12}) 
\end{equation}
we define 
\begin{subequations}
\begin{align}
&\tau_{L}^{n} (\omega) \triangleq \inf \{t \geq 0: \hspace{0.5mm}  \lvert \theta(t,\omega) \rvert > (L - \frac{1}{n})^{\frac{1}{4}} \} \wedge \inf\{t > 0: \hspace{0.5mm}  \lVert \theta(\omega) \rVert_{C_{t}^{\frac{1}{2} - 2\delta}} > (L - \frac{1}{n})^{\frac{1}{2}} \} \wedge L, \label{[Equation (5.3i), HZZ19]}\\
&\tau_{L} \triangleq \lim_{n\to\infty} \tau_{L}^{n}.  \label{[Equation (5.4), HZZ19]}
\end{align}
\end{subequations} 
It follows that $\tau_{L}^{n}$ is a stopping time of $(\bar{\mathcal{B}}_{t})_{t\geq 0}$ and thus $\tau_{L}$ is also a stopping time of $(\bar{\mathcal{B}}_{t})_{t\geq 0}$. For the fixed $(\Omega, \mathcal{F}, \textbf{P})$ we assume Theorem \ref{[Theorem 1.3, HZZ19]} and denote by $u$ the solution constructed from Theorem \ref{[Theorem 1.3, HZZ19]} on $[0, \mathfrak{t}]$ where $\mathfrak{t} = T_{L}$ for $L$ sufficiently large, and 
\begin{equation}\label{[Equation (6.3), HZZ19]}
T_{L} \triangleq \inf\{t > 0: \hspace{0.5mm}  \lvert B(t) \rvert \geq L^{\frac{1}{4}} \} \wedge \inf\{t > 0: \hspace{0.5mm} \lVert B \rVert_{C_{t}^{\frac{1}{2} - 2\delta}} \geq L^{\frac{1}{2}} \} \wedge L. 
\end{equation} 
We furthermore denote by $P$ the law of $(u,B)$ and deduce the following result similarly to Proposition \ref{[Proposition 3.7, HZZ19]}.
\begin{proposition}\label{[Proposition 5.4, HZZ19]}
Let $\tau_{L}$ be defined by \eqref{[Equation (5.4), HZZ19]}. Then $P$, the law of $(u,B)$, is a probabilistically weak solution to \eqref{[Equation (1.7), HZZ19]} on $[0, \tau_{L}]$ that satisfies Definition \ref{[Definition 5.2, HZZ19]}. 
\end{proposition}
\begin{proof}[Proof of Proposition \ref{[Proposition 5.4, HZZ19]}]
The proof is similar to that of Proposition \ref{[Proposition 3.7, HZZ19]} making use of the fact that 
\begin{equation}\label{[Equation (5.4a), HZZ19]}
\theta(t, (u, B)) = B(t) \hspace{1mm} \forall \hspace{1mm} t \in [0, T_{L}] \hspace{1mm} \textbf{P}\text{-almost surely.}
\end{equation}
\end{proof}
Next we extend $P$ on $[0,\tau_{L}]$ to $[0,\infty)$ similarly to Proposition \ref{[Proposition 3.8, HZZ19]}. 
\begin{proposition}\label{[Proposition 5.5, HZZ19]}
The probability measure $P \otimes_{\tau_{L}} R$ in \eqref{[Equation (5.3a), HZZ19]} with $\tau_{L}$ defined in \eqref{[Equation (5.4), HZZ19]} is a probabilistically weak solution to \eqref{[Equation (1.7), HZZ19]} on $[0,\infty)$ that satisfies Definition \ref{[Definition 5.1, HZZ19]}.
\end{proposition}
\begin{proof}[Proof of Proposition \ref{[Proposition 5.5, HZZ19]}]
Due to \eqref{[Equation (5.3i), HZZ19]}, $\tau_{L}$ is a bounded stopping time of $(\bar{\mathcal{B}}_{t})_{t\geq 0}$ and thus, the hypothesis of Lemma \ref{[Proposition 5.2, HZZ19]} is verified. By Proposition \ref{[Proposition 5.4, HZZ19]}, $P$ is a probabilistically weak solution to \eqref{[Equation (1.7), HZZ19]} on $[0,\tau_{L}]$ and therefore Lemma \ref{[Proposition 5.3, HZZ19]} gives us the desired result once we verify the existence of a Borel set $\mathcal{D} \subset \bar{\Omega}_{\tau}$ such that $P(\mathcal{D}) = 0$ and \eqref{[Equation (5.3), HZZ19]} holds for every $\omega \in  \bar{\Omega}_{\tau} \setminus \mathcal{D}$, and that can be achieved similarly to the proof of Proposition \ref{[Proposition 3.8, HZZ19]}. 
\end{proof}

\begin{proof}[Proof of Theorem \ref{[Theorem 1.4, HZZ19]} assuming Theorem \ref{[Theorem 1.3, HZZ19]}]
The proof is similar to that of Theorem \ref{[Theorem 1.2, HZZ19]}; we sketch it for completeness. We fix $T> 0$ arbitrarily, any $\iota \in (0,1)$ and $K > 1$ such that $\iota K^{2} \geq 1$. Then the probability measure $P \otimes_{\tau_{L}} R$ from Proposition \ref{[Proposition 5.5, HZZ19]} satisfies $P \otimes_{\tau_{L}} R (\{ \tau_{L} \geq T \})  > \iota$ due to \eqref{[Equation (5.3a), HZZ19]} which implies $\mathbb{E}^{P \otimes_{\tau_{L}} R} [ \lVert \xi(T) \rVert_{L_{x}^{2}}^{2}] > \iota K^{2} e^{T} \lVert \xi^{\text{in}} \rVert_{L_{x}^{2}}^{2}$, where $\xi^{\text{in}}$ is the deterministic initial condition constructed through Theorem \ref{[Theorem 1.3, HZZ19]}. On the other hand, via a standard Galerkin approximation scheme (e.g. \cite{FR08}), one can construct a probabilistically weak solution $\Theta$ to \eqref{[Equation (1.7), HZZ19]} starting from $\xi^{\text{in}}$ such that $\mathbb{E}^{\Theta} [ \lVert \xi(T) \rVert_{L_{x}^{2}}^{2}] \leq e^{T} \lVert \xi^{\text{in}} \rVert_{L_{x}^{2}}^{2}$. This implies the lack of joint uniqueness in law for \eqref{[Equation (1.7), HZZ19]} and consequently the non-uniqueness in law for \eqref{[Equation (1.7), HZZ19]} by \cite[Theorem C.1]{HZZ19}, which is an infinite-dimensional version of \cite[Theorem 3.1]{C03}. 
\end{proof}

\subsection{Proof of Theorem \ref{[Theorem 1.3, HZZ19]} assuming Proposition \ref{[Proposition 6.2, HZZ19]}}

We define $\Upsilon(t)\triangleq e^{B(t)}, v  \triangleq\Upsilon^{-1} u$ for $t \geq 0$. It follows from Ito's product formula (e.g., \cite[Theorem 4.4.13]{A09}) that 
\begin{equation}\label{[Equation (6.1), HZZ19]}
\partial_{t} v + \frac{1}{2} v + (-\Delta)^{m} v + \Upsilon \text{div} (v\otimes v) + \Upsilon^{-1} \nabla \pi = 0, \hspace{3mm} \nabla\cdot v =0, \hspace{3mm} t> 0.
\end{equation} 
For every $q \in \mathbb{N}_{0}$, we construct $(v_{q}, \mathring{R}_{q})$ that solves 
\begin{equation}\label{[Equation (6.2), HZZ19]}
\partial_{t} v_{q} + \frac{1}{2} v_{q} + (-\Delta)^{m} v_{q} + \Upsilon \text{div} (v_{q}\otimes v_{q}) + \nabla p_{q} = \text{div} \mathring{R}_{q}, \hspace{3mm} \nabla\cdot v_{q} =0, \hspace{3mm} t>0. 
\end{equation} 
Similarly to \eqref{[Equation (4.1a), HZZ19]} in the additive case, we continue to define $\lambda_{q} \triangleq a^{b^{q}}$ and $\delta_{q} \triangleq \lambda_{q}^{-2\beta}$ for $a > 0$ and $b \in \mathbb{N}$. Due to \eqref{[Equation (6.3), HZZ19]}, we obtain for all $L > 1, \delta \in (0, \frac{1}{12})$ and $t \in [0, T_{L}]$, 
\begin{equation}\label{[Equation (6.4), HZZ19]}
\lvert B(t) \rvert \leq L^{\frac{1}{4}} \text{ and } \lVert B \rVert_{C_{t}^{\frac{1}{2} - 2 \delta}} \leq L^{\frac{1}{2}}
\end{equation}
which immediately implies 
\begin{equation}\label{[Equation (6.5), HZZ19]}
\lVert \Upsilon \rVert_{C_{t}^{\frac{1}{2} - 2\delta}} + \lvert \Upsilon(t) \rvert + \lvert \Upsilon^{-1}(t) \rvert \leq e^{L^{\frac{1}{4}}} L^{\frac{1}{2}} + 2e^{L^{\frac{1}{4}}} \leq m_{L}^{2} \text{ where } m_{L} \triangleq \sqrt{3} L^{\frac{1}{4}} e^{\frac{1}{2} L^{\frac{1}{4}}}.
\end{equation} 
Differently from \eqref{estimate 4} we define 
\begin{equation}\label{[Equation (6.6), HZZ19]}
M_{0}(t) \triangleq e^{4Lt + 2L}.
\end{equation} 
For induction we assume that $(v_{q}, \mathring{R}_{q})$ satisfy the following bounds on $[0, T_{L}]$: 
\begin{subequations}
\begin{align}
& \lVert v_{q} \rVert_{C_{t}L_{x}^{2}} \leq m_{L} M_{0}(t)^{\frac{1}{2}} (1+ \sum_{1 \leq  r \leq q} \delta_{r}^{\frac{1}{2}}) \leq 2m_{L} M_{0}(t)^{\frac{1}{2}}, \label{[Equation (6.7a), HZZ19]}\\
& \lVert v_{q} \rVert_{C_{t,x}^{1}} \leq m_{L} M_{0}(t)^{\frac{1}{2}} \lambda_{q}^{4}, \label{[Equation (6.7b), HZZ19]}\\
& \lVert \mathring{R}_{q} \rVert_{C_{t}L_{x}^{1}} \leq M_{0}(t) c_{R} \delta_{q+1}, \label{[Equation (6.7c), HZZ19]}
\end{align}
\end{subequations}  
where $c_{R} > 0$ is again a universal constant to be determined subsequently (e.g., \eqref{[Equation (6.22), HZZ19]}, \eqref{[Equation (6.23), HZZ19]}, \eqref{estimate 50})  and we assumed again $a^{\beta b} > 3$, as formally stated in \eqref{[Equation (6.9), HZZ19]} of Proposition \ref{[Lemma 6.1, HZZ19]}, in order to deduce $\sum_{1\leq r} \delta_{r}^{\frac{1}{2}} < \frac{1}{2}$. We note that \eqref{[Equation (6.7a), HZZ19]} and \eqref{[Equation (6.7b), HZZ19]} are same respectively to \eqref{[Equation (4.4a), HZZ19]} and \eqref{[Equation (4.4b), HZZ19]} multiplied by $m_{L}$ from \eqref{[Equation (6.5), HZZ19]} while \eqref{[Equation (6.7c), HZZ19]} is identical to \eqref{[Equation (4.4c), HZZ19]}. Again, our choice of $\mathring{R}_{0}$ in the following result differs from that of \cite[Lemma 6.1]{HZZ19}.

\begin{proposition}\label{[Lemma 6.1, HZZ19]}
Let $L > 1$ and define 
\begin{equation}\label{[Equation (6.8a), HZZ19]}
v_{0} (t,x) \triangleq \frac{ m_{L} e^{2Lt+ L}}{(2\pi)^{\frac{3}{2}}} 
\begin{pmatrix}
\sin(x_{3}) & 0 & 0 
\end{pmatrix}^{T}. 
\end{equation} 
Then together with 
\begin{equation}\label{[Equation (6.8b), HZZ19]}
\mathring{R}_{0} (t,x) \triangleq \frac{m_{L} (2L + \frac{1}{2}) e^{2L t + L}}{(2\pi)^{\frac{3}{2}}} 
\begin{pmatrix}
0 & 0 & -\cos(x_{3}) \\
0 & 0 & 0\\
-\cos(x_{3}) & 0 & 0 
\end{pmatrix} 
+ \mathcal{R} (-\Delta)^{m} v_{0},
\end{equation}
it satisfies \eqref{[Equation (6.2), HZZ19]} at level $q = 0$. Moreover, \eqref{[Equation (6.7a), HZZ19]}-\eqref{[Equation (6.7c), HZZ19]} are satisfied at level $q = 0$ provided 
\begin{equation}\label{[Equation (6.9), HZZ19]}
18 (2\pi)^{\frac{3}{2}} \sqrt{3} < 2(2\pi)^{\frac{3}{2}} \sqrt{3} a^{2\beta b} \leq \frac{c_{R} e^{L}}{L^{\frac{1}{4}} (2L + 13) e^{\frac{1}{2} L^{\frac{1}{4}}}} \hspace{1mm} \text{ and } \hspace{1mm} L \leq \frac{ (2\pi)^{\frac{3}{2}} a^{4} -2}{2},
\end{equation} 
where the inequality $9 < a^{2\beta b}$ in \eqref{[Equation (6.9), HZZ19]} is assumed for the sake of second inequality in \eqref{[Equation (6.7a), HZZ19]}. Furthermore, $v_{0}(0,x)$ and $\mathring{R}_{0} (0,x)$ are both deterministic. 
\end{proposition}

\begin{proof}[Proof of Proposition \ref{[Lemma 6.1, HZZ19]}]
The proof is very similar to that of Proposition \ref{[Lemma 4.1, HZZ19]}. We provide a sketch for completeness. It may be immediately verified that $v_{0}$ and $\mathring{R}_{0}$ satisfy \eqref{[Equation (6.2), HZZ19]} with $p_{0} = 0$. Moreover, we can readily compute for all $t \in [0, T_{L}]$ 
\begin{equation}\label{[Equation (6.10a), HZZ19]}
\lVert v_{0}(t) \rVert_{L_{x}^{2}} = \frac{ m_{L} M_{0}(t)^{\frac{1}{2}}}{\sqrt{2}} \leq m_{L} M_{0}(t)^{\frac{1}{2}}, \hspace{1mm}  \lVert v_{0} \rVert_{C_{t,x}^{1}} \leq \frac{2(L+1)m_{L}  e^{2L t + L}}{(2\pi)^{\frac{3}{2}}} \overset{\eqref{[Equation (6.9), HZZ19]}} \leq \lambda_{0}^{4} m_{L} M_{0}(t)^{\frac{1}{2}}
\end{equation} 
and thus \eqref{[Equation (6.7a), HZZ19]} and \eqref{[Equation (6.7b), HZZ19]} both hold at level $q = 0$. We can compute similarly to the proof of Proposition \ref{[Lemma 4.1, HZZ19]} using the facts that $\Delta v_{0} = -v_{0}$ and $\lVert \mathcal{R} (-\Delta)^{m} f \rVert_{L_{x}^{2}} \leq 24\lVert \Delta f \rVert_{L_{x}^{2}}$, 
\begin{equation}
\lVert \mathring{R}_{0}(t) \rVert_{L_{x}^{1}} \leq m_{L} (2L + \frac{1}{2}) e^{2L t + L} 8 (2\pi)^{\frac{1}{2}} + (2\pi)^{\frac{3}{2}} 24 \lVert v_{0} \rVert_{L_{x}^{2}}  \overset{\eqref{[Equation (6.5), HZZ19]} \eqref{[Equation (6.9), HZZ19]}}{\leq} M_{0}(t) c_{R} \delta_{1}.
\end{equation} 
At last, it is clear that $v_{0}(0,x)$ is deterministic and consequently so is $\mathring{R}_{0}(0,x)$.
\end{proof}
Let us point out that 
\begin{equation}\label{estimate 38}
18 (2\pi)^{\frac{3}{2}} \sqrt{3} < \frac{c_{R} e^{L}}{ L^{\frac{1}{4}} (2L + 13)e^{\frac{1}{2} L^{\frac{1}{4}}}} 
\end{equation}
is not sufficient but certainly necessary to satisfy \eqref{[Equation (6.9), HZZ19]}. 
\begin{proposition}\label{[Proposition 6.2, HZZ19]}
Let $L > 1$ satisfy \eqref{estimate 38} and suppose that $(v_{q}, \mathring{R}_{q})$ is an $(\mathcal{F}_{t})_{t\geq 0}$-adapted solution to \eqref{[Equation (6.2), HZZ19]} that satisfies \eqref{[Equation (6.7a), HZZ19]}-\eqref{[Equation (6.7c), HZZ19]}. Then there exists a choice of parameters $a, b, \beta$ such that \eqref{[Equation (6.9), HZZ19]} is fulfilled and there exists an $(\mathcal{F}_{t})_{t\geq 0}$-adapted process $(v_{q+1}, \mathring{R}_{q+1})$ that solves \eqref{[Equation (6.2), HZZ19]}, satisfies \eqref{[Equation (6.7a), HZZ19]}-\eqref{[Equation (6.7c), HZZ19]} at level $q + 1$ and 
\begin{equation}\label{[Equation (6.13), HZZ19]}
\lVert v_{q+1}(t) - v_{q}(t) \rVert_{L_{x}^{2}} \leq m_{L} M_{0}(t)^{\frac{1}{2}} \delta_{q+1}^{\frac{1}{2}} \hspace{3mm} \forall \hspace{1mm} t \in [0,T_{L}].
\end{equation} 
Furthermore, if $v_{q}(0, x)$ and $\mathring{R}_{q}(0,x)$ are deterministic, then so are $v_{q+1}(0,x)$ and $\mathring{R}_{q+1}(0,x)$. 
\end{proposition}

We assume Proposition \ref{[Proposition 6.2, HZZ19]} and prove Theorem \ref{[Theorem 1.3, HZZ19]} now. 
\begin{proof}[Proof of Theorem  \ref{[Theorem 1.3, HZZ19]} assuming Proposition \ref{[Proposition 6.2, HZZ19]}]
The proof is similar to that of Theorem \ref{[Theorem 1.1, HZZ19]}; we sketch it for completeness. Let us fix any $T > 0, K > 1$ and $\iota \in (0,1)$, and then take $L$ that satisfies \eqref{estimate 38} and enlarge it if necessary to satisfy 
\begin{equation}\label{[Equation (6.15), HZZ19]}
(\frac{1}{\sqrt{2}} - \frac{1}{2}) e^{2L T} > (\frac{3}{2} ) e^{2L^{\frac{1}{2}}} \text{ and } L > [ \ln (K e^{\frac{T}{2}})]^{2}.
\end{equation}
We can start from $(v_{0}, \mathring{R}_{0})$ in Proposition \ref{[Lemma 6.1, HZZ19]}, and via Proposition \ref{[Proposition 6.2, HZZ19]} inductively obtain a sequence $(v_{q}, \mathring{R}_{q})$ that satisfies \eqref{[Equation (6.2), HZZ19]}, \eqref{[Equation (6.7a), HZZ19]}-\eqref{[Equation (6.7c), HZZ19]} and \eqref{[Equation (6.13), HZZ19]}. Identically to \eqref{[Equation (4.12a), HZZ19]} we can show that for any $\gamma \in (0, \frac{\beta}{4+\beta})$ and any $t \in [0, T_{L}]$, $\sum_{q\geq 0} \lVert v_{q+1}(t) - v_{q}(t) \rVert_{H_{x}^{\gamma}} \lesssim m_{L} M_{0}(t)^{\frac{1}{2}}$ by \eqref{[Equation (6.13), HZZ19]} and \eqref{[Equation (6.7b), HZZ19]}, which allows us to deduce the limiting solution $\lim_{q\to\infty} v_{q} \triangleq v \in C([0, T_{L}]; H^{\gamma}(\mathbb{T}^{3}))$ that is $(\mathcal{F}_{t})_{t\geq 0}$-adapted because each $v_{q}$ is $(\mathcal{F}_{t})_{t\geq 0}$-adapted. Because $u = \Upsilon v = e^{B_{t}} v$ where $\lvert e^{B_{t}} \rvert \leq e^{L^{\frac{1}{4}}}$ for all $t \in [0, T_{L}]$ due to \eqref{[Equation (6.4), HZZ19]}, we are able to deduce \eqref{[Equation (1.7a), HZZ19]} by choosing $\mathfrak{t} = T_{L}$ for $L$ sufficiently large. Next, because $\lim_{q\to\infty} \lVert \mathring{R}_{q} \rVert_{C_{T_{L}}L^{1}} = 0$ due to \eqref{[Equation (6.7c), HZZ19]}, $v$ is a weak solution to \eqref{[Equation (6.1), HZZ19]} on $[0, T_{L}]$.  
Moreover, we can assume $b \geq 2$ and show identically to \eqref{[Equation (4.14a), HZZ19]} that for all $t \in [0, T_{L}]$, $\lVert v(t) - v_{0}(t) \rVert_{L_{x}^{2}} \leq \frac{m_{L}}{2} M_{0}(t)^{\frac{1}{2}}$ by \eqref{[Equation (6.13), HZZ19]} and \eqref{[Equation (4.1a), HZZ19]} which in turn implies 
\begin{equation}\label{[Equation (6.14a), HZZ19]}  
e^{2 L^{\frac{1}{2}}} \lVert v(0) \rVert_{L_{x}^{2}}  \leq  e^{2L^{\frac{1}{2}}} ( \lVert v(0) - v_{0}(0) \rVert_{L_{x}^{2}} + \lVert v_{0} (0) \rVert_{L_{x}^{2}})\overset{\eqref{[Equation (6.10a), HZZ19]}}{\leq} e^{2L^{\frac{1}{2}}} (\frac{3}{2} )m_{L} M_{0}(0)^{\frac{1}{2}}.   
\end{equation} 
These lead us to, on a set $\{T_{L} \geq T \}$, 
\begin{equation}\label{[Equation (6.14), HZZ19]}
\lVert v(T) \rVert_{L_{x}^{2}} \overset{\eqref{[Equation (6.10a), HZZ19]}}{\geq} \frac{m_{L} M_{0}(T)^{\frac{1}{2}}}{\sqrt{2}} - \lVert v(T) - v_{0}(T) \rVert_{L_{x}^{2}}  \overset{ \eqref{[Equation (6.15), HZZ19]}\eqref{[Equation (6.14a), HZZ19]}  }{\geq} e^{2L^{\frac{1}{2}}} \lVert v(0) \rVert_{L_{x}^{2}}^{2}. 
\end{equation} 
Moreover, for the fixed $T > 0$, $\iota \in (0,1)$, one can take $L$ even larger to deduce $\textbf{P} (\{T_{L} \geq T \}) > \iota$. We also see that $u^{\text{in}}(x) = \Upsilon(0)v(0,x) = v(0,x)$ which is deterministic because $v_{q}(0,x)$ is deterministic for all $q \in \mathbb{N}_{0}$ by Propositions \ref{[Lemma 6.1, HZZ19]} and \ref{[Proposition 6.2, HZZ19]}. Clearly from \eqref{[Equation (6.1), HZZ19]}, $u = \Upsilon v$ is a $(\mathcal{F}_{t})_{t\geq 0}$-adapted solution to \eqref{[Equation (1.7), HZZ19]}. Finally, we can deduce on the set $\{\mathfrak{t} \geq T \}$ 
\begin{equation*}
\lVert u(T) \rVert_{L_{x}^{2}} \overset{\eqref{[Equation (6.4), HZZ19]}  \eqref{[Equation (6.14), HZZ19]}}{\geq} e^{L^{\frac{1}{2}}} \lVert u^{\text{in}} \rVert_{L_{x}^{2}} \overset{\eqref{[Equation (6.15), HZZ19]} }{>} K e^{\frac{T}{2}} \lVert u^{\text{in}} \rVert_{L_{x}^{2}}
\end{equation*} 
which verifies \eqref{[Equation (1.7a'), HZZ19]} and completes the proof of Theorem \ref{[Theorem 1.3, HZZ19]}.
\end{proof}

\subsection{Proof of Proposition \ref{[Proposition 6.2, HZZ19]}}
We now proceed with the proof of Proposition \ref{[Proposition 6.2, HZZ19]}. 
\subsubsection{Choice of parameters}
We fix $L$ sufficiently large so that it satisfies \eqref{estimate 38}. We can take the same $\alpha$ from \eqref{alpha} and $l$ from \eqref{[Equation (4.17), HZZ19]}, and that allows us to fix the special $p^{\ast}$ from \eqref{p}. We will need $b \in \mathbb{N}$ to satisfy $\alpha b > 16$ and $\alpha > 16 \beta b$, the latter of which is possible by taking $\beta > 0$ sufficiently small after fixing $b \in \mathbb{N}$ first. Now for the $\alpha> 0$ fixed, we fix $b > \frac{16}{\alpha}$ so that $\alpha b > 16$. Concerning the conditions of \eqref{[Equation (6.9), HZZ19]}, we see that the last inequality of $L \leq \frac{ (2\pi)^{\frac{3}{2}} a^{4} -2}{2}$ is certainly satisfied by taking $a$ sufficiently large while $18 (2\pi)^{\frac{3}{2}} \sqrt{3} < 2(2\pi)^{\frac{3}{2}} \sqrt{3} a^{2\beta b} \leq \frac{c_{R} e^{L}}{L^{\frac{1}{4}} (2L + 13)e^{\frac{1}{2} L^{\frac{1}{4}}}}$ is satisfied by taking $\beta> 0$ sufficiently small. Because we chose $L$ to satisfy \eqref{estimate 38}, this is possible. Thus, we now consider $L, \alpha, l$ and $b$ all fixed while reserve the freedom to take $a > 0$ as large as we wish, requiring that $a^{\frac{25m-20}{24}} \in \mathbb{N}$, and $\beta > 0$ as small as we wish. 

\subsubsection{Mollification}

Similarly to \eqref{[Equation (4.18c), HZZ19]} we mollify $v_{q}, \mathring{R}_{q}$ and $\Upsilon(t) = e^{B_{t}}$ so that 
\begin{equation}\label{[Equation (6.16b), HZZ19]}
v_{l} \triangleq (v_{q} \ast_{x} \phi_{l}) \ast_{t} \varphi_{l}, \hspace{3mm} \mathring{R}_{l} \triangleq (\mathring{R}_{q} \ast_{x} \phi_{l}) \ast_{t} \varphi_{l}, \hspace{3mm} \Upsilon_{l} \triangleq \Upsilon \ast_{t} \varphi_{l},
\end{equation} 
where we require again that $\varphi_{l}$ is compactly supported in $\mathbb{R}_{+}$. Now because $(v_{q}, \mathring{R}_{q})$ solves \eqref{[Equation (6.2), HZZ19]}, we see that 
\begin{equation}\label{[Equation (6.16c), HZZ19]}
 \partial_{t} v_{l} + \frac{1}{2} v_{l} + (-\Delta)^{m} v_{l} + \Upsilon_{l}\text{div} (v_{l} \otimes v_{l}) + \nabla p_{l}  =  \text{div}( \mathring{R}_{l}  + R_{\text{commutator 1}})
\end{equation}
where 
\begin{subequations}
\begin{align}
& p_{l} \triangleq (p_{q} \ast_{x} \phi_{l}) \ast_{t} \varphi_{l} - \frac{1}{3} (\Upsilon_{l} \lvert v_{l} \rvert^{2} - ((\Upsilon \lvert v_{q} \rvert^{2} )\ast_{x} \phi_{l} ) \ast_{t} \varphi_{l}), \label{[Equation (6.16e), HZZ19]}\\
& R_{\text{commutator 1}} \triangleq -((\Upsilon (v_{q} \mathring{\otimes} v_{q})) \ast_{x} \phi_{l} )\ast_{t} \varphi_{l} + \Upsilon_{l} (v_{l} \mathring{\otimes} v_{l}).  \label{[Equation (6.16d), HZZ19]}
\end{align}
\end{subequations} 
Similarly to \eqref{[Equation (4.20), HZZ19]}, \eqref{[Equation (4.21), HZZ19]}, \eqref{[Equation (4.22), HZZ19]} we can use the facts that $16 < \alpha b$ and $16 \beta b < \alpha$ so that $\beta < \alpha$ to estimate for all $t \in [0, T_{L}]$ and $N \geq 1$,  
\begin{subequations}
\begin{align}
& \lVert v_{q} - v_{l} \rVert_{C_{t}L_{x}^{2}} \overset{\eqref{[Equation (4.16), HZZ19]}\eqref{[Equation (6.7b), HZZ19]} }{\lesssim} \lambda_{q+1}^{-\alpha} m_{L} M_{0}(t)^{\frac{1}{2}} \leq \frac{m_{L}}{4} M_{0}(t)^{\frac{1}{2}} \delta_{q+1}^{\frac{1}{2}},  \label{[Equation (6.17), HZZ19]}\\
& \lVert v_{l} \rVert_{C_{t}L_{x}^{2}}  \overset{\eqref{[Equation (6.7a), HZZ19]}}{\leq} m_{L} M_{0}(t)^{\frac{1}{2}} (1+ \sum_{1 \leq r \leq q} \delta_{r}^{\frac{1}{2}}) \overset{\eqref{[Equation (6.9), HZZ19]}}{\leq} 2m_{L} M_{0}(t)^{\frac{1}{2}}, \label{[Equation (6.18), HZZ19]}\\
& \lVert v_{l} \rVert_{C_{t,x}^{N}}\overset{\eqref{[Equation (6.7b), HZZ19]}}{\lesssim} l^{-N +1} m_{L} M_{0}(t)^{\frac{1}{2}} \lambda_{q}^{4} 
\overset{\eqref{[Equation (4.17), HZZ19]}}{\leq} l^{-N} m_{L} M_{0}(t)^{\frac{1}{2}} \lambda_{q+1}^{-\alpha}. \label{[Equation (6.19), HZZ19]}
\end{align}
\end{subequations}
Finally, we will continue to use the same choices of $r_{\lVert}, r_{\bot}$ and $\mu$ in \eqref{[Equation (B.2a), HZZ19]}. 

\subsubsection{Perturbation}
We proceed with same definitions of $\chi$ in \eqref{[Equation (4.23a), HZZ19]} and $\rho$ in \eqref{[Equation (4.23b), HZZ19]} identically except that $M_{0}(t)$ is now defined by \eqref{[Equation (6.6), HZZ19]} instead of \eqref{estimate 4}. We point out that although our definition of $\mathring{R}_{0}$ in \eqref{[Equation (6.8b), HZZ19]} was different from $\mathring{R}_{0}$ in \eqref{[Equation (4.6), HZZ19]}, \eqref{[Equation (4.23c), HZZ19]} remains valid as its proof relied only on the definition of $\chi$. Moreover, \eqref{[Equation (4.24), HZZ19]} also remains applicable in our case as its proof only depended on the definitions of $\rho$ and $\chi$, not $M_{0}(t)$ or $\mathring{R}_{l}$. We define the following modified amplitude function
\begin{equation}\label{[Equation (6.20), HZZ19]}
\bar{a}_{(\xi)} (\omega, t, x) \triangleq \bar{a}_{\xi, q+1} (\omega, t, x) \triangleq \Upsilon_{l}^{-\frac{1}{2}} a_{(\xi)} (\omega, t, x) 
\end{equation} 
where $a_{(\xi)}(\omega, t, x)$ is defined in \eqref{[Equation (4.26), HZZ19]}. Similarly to \eqref{[Equation (4.27), HZZ19]}-\eqref{[Equation (4.28), HZZ19]} we can deduce the following identity and an estimate for all $t \in [0, T_{L}]$:
\begin{subequations}
\begin{align} 
&(2\pi)^{\frac{3}{2}} \sum_{\xi \in \Lambda} \bar{a}_{(\xi)}^{2} \fint_{\mathbb{T}^{3}} W_{(\xi)} \otimes W_{(\xi)} dx
\overset{\eqref{[Equation (6.20), HZZ19]} \eqref{[Equation (B.5), HZZ19]}}{=} \Upsilon_{l}^{-1} (\rho \text{Id} - \mathring{R}_{l}), \label{[Equation (6.21), HZZ19]}\\
&  \lVert \bar{a}_{(\xi)} \rVert_{C_{t}L_{x}^{2}}  \overset{ \eqref{[Equation (4.23c), HZZ19]} \eqref{[Equation (4.24), HZZ19]}\eqref{[Equation (6.7c), HZZ19]}  \eqref{[Equation (6.20), HZZ19]} \eqref{[Equation (A.9e), HZZ19]} } {\lesssim}
 \frac{m_{L} \sqrt{12} c_{R}^{\frac{1}{2}} \delta_{q+1}^{\frac{1}{2}} M_{0}(t)^{\frac{1}{2}} M}{8\lvert \Lambda \rvert} \leq \frac{c_{R}^{\frac{1}{4}} m_{L} M_{0}(t)^{\frac{1}{2}} \delta_{q+1}^{\frac{1}{2}}}{2 \lvert \Lambda \rvert}, \label{[Equation (6.22), HZZ19]}
\end{align} 
\end{subequations}
where we used a simple estimate of $\lVert \Upsilon_{l}^{-\frac{1}{2}} \rVert_{C_{t}} \leq m_{L}$ due to \eqref{[Equation (6.4), HZZ19]} and \eqref{[Equation (6.5), HZZ19]} and we also took $c_{R} \ll M^{-4}$. Because \eqref{[Equation (4.4c), HZZ19]} and \eqref{[Equation (6.7c), HZZ19]} are identical only except the precise definition of $M_{0}(t)$, we can use \eqref{[Equation (4.24a), HZZ19]} and consequently \eqref{[Equation (4.29), HZZ19]}. Moreover, we can take $a > 0$ sufficiently large so that $a^{26} \geq \sqrt{3} L^{\frac{1}{4}} e^{\frac{1}{2} L^{\frac{1}{4}}}$; considering $\alpha b > 16$, this gives us $m_{L} \leq l^{-1}$. Hence, we may directly apply $\partial_{t}$ and $\nabla$, use \eqref{[Equation (4.29), HZZ19]}, as well as $\lVert \Upsilon_{l}^{-\frac{1}{2}} \rVert_{C_{t}} \leq m_{L}$ and take $c_{R}^{\frac{1}{8}} < C^{-1}$ for any implicit constant $C$ to deduce the following estimates: for any $N \geq 0$, 
\begin{subequations}\label{[Equation (6.23), HZZ19]}
\begin{align}
& \lVert \bar{a}_{(\xi)} \rVert_{C_{t} C_{x}^{N}} \leq m_{L} c_{R}^{\frac{1}{4}} \delta_{q+1}^{\frac{1}{2}} M_{0}(t)^{\frac{1}{2}} l^{-2-5N},  \hspace{3mm} \lVert \bar{a}_{(\xi)} \rVert_{C_{t}^{1}C_{x}} \leq m_{L} c_{R}^{\frac{1}{8}} \delta_{q+1}^{\frac{1}{2}} M_{0}(t)^{\frac{1}{2}} l^{-7}, \\
& \lVert \bar{a}_{(\xi)} \rVert_{C_{t}^{1}C_{x}^{1}} \leq m_{L} c_{R}^{\frac{1}{8}} \delta_{q+1}^{\frac{1}{2}} M_{0}(t)^{\frac{1}{2}} l^{-12}, \hspace{7mm} \lVert \bar{a}_{(\xi)} \rVert_{C_{t}^{1}C_{x}^{2}} \leq m_{L} c_{R}^{\frac{1}{8}} \delta_{q+1}^{\frac{1}{2}} M_{0}(t)^{\frac{1}{2}} l^{-17}. 
\end{align} 
\end{subequations}
Now we define $w_{q+1}^{(p)}$ identically as we did in \eqref{[Equation (4.30), HZZ19]} with $a_{(\xi)}$ replaced by $\bar{a}_{(\xi)}$ and $M_{0}(t)$ from \eqref{[Equation (6.6), HZZ19]} within the definition of $\rho(\omega, t, x)$ so that we obtain an identity of 
\begin{equation}\label{[Equation (6.23a), HZZ19]}
\Upsilon_{l} w_{q+1}^{(p)} \otimes w_{q+1}^{(p)} + \mathring{R}_{l} 
\overset{\eqref{[Equation (6.21), HZZ19]} \eqref{[Equation (B.4), HZZ19]}}{=} \Upsilon_{l} \sum_{\xi \in \Lambda} \bar{a}_{(\xi)}^{2} \mathbb{P}_{\neq 0} (W_{(\xi)} \otimes W_{(\xi)}) + \rho \text{Id}. 
\end{equation} 
Moreover, we define $w_{q+1}^{(c)}$ exactly as in \eqref{[Equation (4.32), HZZ19]} with $a_{(\xi)}$ replaced by $\bar{a}_{(\xi)}$ while $w_{q+1}^{(t)}$ with $a_{(\xi)}$ from \eqref{[Equation (4.26), HZZ19]}, both with $M_{0}(t)$ from \eqref{[Equation (6.6), HZZ19]} within the definition of $\rho(\omega, t, x)$. We emphasize here that the definition of $w_{q+1}^{(t)}$ remained with $a_{(\xi)}$ rather than $\bar{a}_{(\xi)}$,  and it will be crucial in deriving the identity \eqref{[Equation (6.33c), HZZ19]}. At last, we define 
\begin{equation}\label{[Equation (6.23b), HZZ19]}
w_{q+1} \triangleq w_{q+1}^{(p)} + w_{q+1}^{(c)} + w_{q+1}^{(t)} \text{ and } v_{q+1} \triangleq v_{l} + w_{q+1}, 
\end{equation} 
which are both divergence-free by the same reasoning as in the proof of Theorem \ref{[Theorem 1.1, HZZ19]}, and additionally $w_{q+1}$ has its mean zero. Next, similarly to \eqref{[Equation (4.37), HZZ19]} using  Lemma \ref{[Lemma 7.4, BV19b]}, \eqref{[Equation (4.38), HZZ19]} and \eqref{[Equation (4.39), HZZ19]}, we can deduce the following estimates for all $t \in [0,T_{L}]$: 
\begin{subequations}\label{estimate 50}
\begin{align}
& \lVert w_{q+1}^{(p)} \rVert_{C_{t}L_{x}^{2}} \lesssim \sup_{s \in [0,t]} \sum_{\xi \in \Lambda} \frac{ m_{L} c_{R}^{\frac{1}{4}} \delta_{q+1}^{\frac{1}{2}} M_{0}(s)^{\frac{1}{2}}}{2 \lvert \Lambda \rvert} \lVert W_{(\xi)}(s) \rVert_{L_{x}^{2}}\leq \frac{1}{2} m_{L} M_{0} (t)^{\frac{1}{2}} \delta_{q+1}^{\frac{1}{2}}, \label{[Equation (6.24), HZZ19]} \\
& \lVert w_{q+1}^{(p)} \rVert_{C_{t}L_{x}^{p}} \leq \sup_{s \in [0,t]} \sum_{\xi \in \Lambda} \lVert \bar{a}_{(\xi)}(s) \rVert_{L_{x}^{\infty}} \lVert W_{(\xi)}(s) \rVert_{L_{x}^{p}} \overset{\eqref{[Equation (6.23), HZZ19]} \eqref{[Equation (B.7c), HZZ19]}}{\lesssim} m_{L} M_{0}(t)^{\frac{1}{2}} \delta_{q+1}^{\frac{1}{2}} l^{-2} r_{\bot}^{\frac{2}{p} - 1} r_{\lVert}^{\frac{1}{p} - \frac{1}{2}}, \label{[Equation (6.25), HZZ19]}\\
& \lVert w_{q+1}^{(c)} \rVert_{C_{t}L_{x}^{p}}  \overset{\eqref{[Equation (6.23), HZZ19]} \eqref{[Equation (B.7c), HZZ19]}}{\lesssim} m_{L} \delta_{q+1}^{\frac{1}{2}} M_{0}(t)^{\frac{1}{2}} l^{-12} r_{\bot}^{\frac{2}{p}} r_{\lVert}^{\frac{1}{p} - \frac{3}{2}},  \label{[Equation (6.26), HZZ19]}
\end{align} 
\end{subequations}
where $c_{R}$ was taken sufficiently small in \eqref{[Equation (6.24), HZZ19]} to offset any implicit constant. Finally, let us simply mention that \eqref{[Equation (4.40), HZZ19]} still applies to current $w_{q+1}^{(t)}$ with $M_{0}(t)$ defined by \eqref{[Equation (6.6), HZZ19]} because the definition of $w_{q+1}^{(t)}$ in the current case is same as that in the proof of Theorem \ref{[Theorem 1.1, HZZ19]} except the definition of $M_{0}(t)$ which played no role in that estimate. Since $\alpha$ from \eqref{alpha} satisfies both $24 \alpha + \frac{20m-25}{24} < 0, 8 \alpha + \frac{4m-5}{8} < 0$, we may take $a > 0$ sufficiently large to deduce from \eqref{[Equation (6.23b), HZZ19]} 
\begin{align}
 \lVert w_{q+1} \rVert_{C_{t}L_{x}^{2}}  \overset{\eqref{[Equation (4.16), HZZ19]} \eqref{[Equation (4.40), HZZ19]}  \eqref{[Equation (6.24), HZZ19]} \eqref{[Equation (6.26), HZZ19]}  }{\leq}&
 m_{L} M_{0}(t)^{\frac{1}{2}} \delta_{q+1}^{\frac{1}{2}} [\frac{1}{2} + C (\lambda_{q+1}^{24 \alpha + \frac{20m-25}{24}} +  M_{0}(t)^{\frac{1}{2}} \delta_{q+1}^{\frac{1}{2}} \lambda_{q+1}^{8\alpha + \frac{4m-5}{8}})] \nonumber \\
 &\leq \frac{3m_{L} M_{0}(t)^{\frac{1}{2}} \delta_{q+1}^{\frac{1}{2}}}{4}.  \label{[Equation (6.28), HZZ19]}
\end{align} 
It follows from \eqref{[Equation (6.23b), HZZ19]}, \eqref{[Equation (6.28), HZZ19]} and \eqref{[Equation (6.18), HZZ19]} that \eqref{[Equation (6.7a), HZZ19]} at level $q+1$ is verified. Moreover, \eqref{[Equation (6.23b), HZZ19]}, \eqref{[Equation (6.28), HZZ19]} and \eqref{[Equation (6.17), HZZ19]} also imply \eqref{[Equation (6.13), HZZ19]}. 

Next, similarly to \eqref{[Equation (4.43), HZZ19]}, \eqref{[Equation (4.44), HZZ19]}, for $t \in [0, T_{L}]$, due to \eqref{[Equation (B.7c), HZZ19]} we can estimate 
\begin{subequations}
\begin{align}
& \lVert w_{q+1}^{(p)} \rVert_{C_{t,x}^{1}} 
\leq \sum_{\xi \in \Lambda} \lVert \bar{a}_{(\xi)} \rVert_{C_{t,x}^{1}} \lVert W_{(\xi)} \rVert_{C_{t,x}^{1}} \overset{\eqref{[Equation (6.23), HZZ19]}}{\lesssim} m_{L} M_{0}(t)^{\frac{1}{2}} l^{-7} r_{\bot}^{-1} r_{\lVert}^{-\frac{1}{2}} \lambda_{q+1}^{2m},  \label{[Equation (6.29), HZZ19]}\\
&  \lVert w_{q+1}^{(c)} \rVert_{C_{t,x}^{1}}  \overset{\eqref{[Equation (6.23), HZZ19]}}{\lesssim} m_{L} M_{0}(t)^{\frac{1}{2}} l^{-17} r_{\lVert}^{-\frac{3}{2}} \lambda_{q+1} [ \frac{r_{\bot} \mu}{r_{\lVert}} + 1] 
\lesssim m_{L} M_{0}(t)^{\frac{1}{2}} l^{-17} r_{\lVert}^{-\frac{3}{2}} \lambda_{q+1}^{2m}. \label{[Equation (6.30), HZZ19]}
\end{align}
\end{subequations} 
Moreover, we simply mention that the estimate \eqref{[Equation (4.45), HZZ19]} with $M_{0}(t)$ defined by \eqref{[Equation (6.6), HZZ19]} remains valid in the current case because $w_{q+1}^{(t)}$ is defined identically as in the proof of Theorem \ref{[Theorem 1.1, HZZ19]}, only with $M_{0}(t)$ defined by \eqref{[Equation (6.6), HZZ19]}.  Using \eqref{[Equation (6.23b), HZZ19]}, \eqref{[Equation (6.29), HZZ19]}, \eqref{[Equation (6.30), HZZ19]} and \eqref{[Equation (4.45), HZZ19]} with $M_{0}(t)$ defined by \eqref{[Equation (6.6), HZZ19]}, as well as \eqref{estimate 16} - \eqref{estimate 19}, imply that \eqref{[Equation (6.7b), HZZ19]} at level $q+1$ is satisfied. 

Next it can be verified that the following analog of \eqref{[Equation (4.32a), HZZ19]} continues to hold
\begin{equation}\label{estimate 39}
\sum_{\xi \in \Lambda} \text{curl curl} (\bar{a}_{(\xi)} V_{(\xi)} ) \overset{\eqref{[Equation (B.6), HZZ19]}}{=} w_{q+1}^{(c)} + w_{q+1}^{(p)}.  
\end{equation} 
Using this identity, we can readily estimate $w_{q+1}^{(p)} + w_{q+1}^{(c)}$ in the norms of $C_{t}L^{p}, C_{t}W^{1,p}$,  and $C_{t}W^{2,p}$ and interpolate again similarly to \eqref{estimate 29} to deduce for all $t \in [0, T_{L}]$ and $p \in (1,\infty)$, 
\begin{equation}\label{estimate 40}  
\lVert w_{q+1}^{(p)} + w_{q+1}^{(c)} \rVert_{C_{t}W_{x}^{2m-1, p}} \lesssim m_{L} M_{0}(t)^{\frac{1}{2}} r_{\bot}^{\frac{2}{p} - 1} r_{\lVert}^{\frac{1}{p} - \frac{1}{2}} l^{-2} \lambda_{q+1}^{2m-1}
\end{equation} 
while $\lVert w_{q+1}^{(t)} \rVert_{C_{t}W_{x}^{2m-1, p}} \lesssim \mu^{-1} M_{0}(t) r_{\bot}^{\frac{2}{p} - 2} r_{\lVert}^{\frac{1}{p} -1} l^{-4} \lambda_{q+1}^{2m-1}$ from \eqref{estimate 30} remains valid with $M_{0}(t)$ defined by \eqref{[Equation (6.6), HZZ19]} because the definition of $w_{q+1}^{(t)}$ remained same as that from the proof of Theorem \ref{[Theorem 1.1, HZZ19]} with the only exception being the definition of $M_{0}(t)$. 

\subsubsection{Reynolds stress} 
Due to \eqref{[Equation (6.2), HZZ19]}, \eqref{[Equation (6.16c), HZZ19]} and \eqref{[Equation (6.23b), HZZ19]} we see that 
\begin{align}
 \text{div} \mathring{R}_{q+1} - \nabla p_{q+1} =& \underbrace{\frac{1}{2} w_{q+1} + (-\Delta)^{m} w_{q+1} + \partial_{t} (w_{q+1}^{(p)} + w_{q+1}^{(c)}) + \Upsilon_{l} \text{div} (v_{l} \otimes w_{q+1} + w_{q+1} \otimes v_{l})}_{\text{div} (R_{\text{linear}} ) + \nabla p_{\text{linear}}} \nonumber  \\
& + \underbrace{\Upsilon_{l}\text{div} ((w_{q+1}^{(c)} + w_{q+1}^{(t)}) \otimes w_{q+1} + w_{q+1}^{(p)} \otimes (w_{q+1}^{(c)} + w_{q+1}^{(t)}))}_{\text{div} (R_{\text{corrector}}) + \nabla p_{\text{corrector}}} \nonumber \\
&+ \underbrace{\text{div} (\Upsilon_{l} w_{q+1}^{(p)} \otimes w_{q+1}^{(p)} + \mathring{R}_{l}) + \partial_{t} w_{q+1}^{(t)}}_{\text{div} (R_{\text{oscillation}}) + \nabla p_{\text{oscillator}}}  \nonumber \\
&+ \underbrace{( \Upsilon - \Upsilon_{l}) \text{div} (v_{q+1} \otimes v_{q+1})}_{\text{div} (R_{\text{commutator 2}}) + \nabla p_{\text{commutator 2}}} + \text{div} (R_{\text{commutator 1}}) - \nabla p_{l} \label{estimate 53}
\end{align} 
with 
\begin{subequations}\label{[Equation (6.33b), HZZ19]}
\begin{align}
& R_{\text{linear}} \triangleq \mathcal{R} ( \frac{1}{2} w_{q+1} + (-\Delta)^{m} w_{q+1} + \partial_{t} (w_{q+1}^{(p)} + w_{q+1}^{(c)})) + \Upsilon_{l} (v_{l} \mathring{\otimes} w_{q+1} + w_{q+1} \mathring{\otimes}v_{l}), \\
& p_{\text{linear}} \triangleq \Upsilon_{l} (\frac{2}{3}) (v_{l} \cdot w_{q+1}), \\
& R_{\text{corrector}} \triangleq \Upsilon_{l}((w_{q+1}^{(c)} + w_{q+1}^{(t)}) \mathring{\otimes} w_{q+1} + w_{q+1}^{(p)} \mathring{\otimes} (w_{q+1}^{(c)} + w_{q+1}^{(t)})), \\
& p_{\text{corrector}} \triangleq \frac{\Upsilon_{l}}{3} ((w_{q+1}^{(c)} + w_{q+1}^{(t)}) \cdot w_{q+1} + w_{q+1}^{(p)} \cdot (w_{q+1}^{(c)} + w_{q+1}^{(t)})) , \\
&R_{\text{oscillation}} \triangleq \sum_{\xi \in \Lambda} \mathcal{R} (\nabla a_{(\xi)}^{2} \mathbb{P}_{\neq 0} (W_{(\xi)} \otimes W_{(\xi)})) -  \mu^{-1} \sum_{\xi \in \Lambda} \mathcal{R} (\partial_{t} a_{(\xi)}^{2} \phi_{(\xi)}^{2} \psi_{(\xi)}^{2} \xi), \label{[Equation (6.33c), HZZ19]}\\
&p_{\text{oscillation}} \triangleq \rho + \Delta^{-1} \text{div}[ \mu^{-1} \sum_{\xi \in \Lambda } \mathbb{P}_{\neq 0} \partial_{t} (a_{(\xi)}^{2} \phi_{(\xi)}^{2} \psi_{(\xi)}^{2} \xi)], \\
& R_{\text{commutator 2}} \triangleq (\Upsilon - \Upsilon_{l}) (v_{q+1} \mathring{\otimes} v_{q+1}), \\
& p_{\text{commutator 2}}  \triangleq \frac{ \Upsilon - \Upsilon_{l}}{3} \lvert v_{q+1} \rvert^{2}, 
\end{align}
\end{subequations} 
where we made use of the following identity concerning $R_{\text{oscillation}}$: 
\begin{align*}
& \text{div} (\Upsilon_{l} w_{q+1}^{(p)} \otimes w_{q+1}^{(p)} + \mathring{R}_{l}) + \partial_{t} w_{q+1}^{(t)} 
\overset{  \eqref{[Equation (4.33), HZZ19]} \eqref{[Equation (6.23a), HZZ19]} \eqref{estimate 7}  }{=}  \sum_{\xi \in \Lambda} \mathbb{P}_{\neq 0} (\nabla a_{(\xi)}^{2} \mathbb{P}_{\neq 0} (W_{(\xi ) } \otimes W_{(\xi)})) + \nabla \rho \\
& \hspace{30mm} + \nabla \Delta^{-1} \text{div} \mu^{-1} \sum_{\xi \in \Lambda } \mathbb{P}_{\neq 0} \partial_{t} (a_{(\xi)}^{2} \phi_{(\xi)}^{2} \psi_{(\xi)}^{2} \xi) - \mu^{-1} \sum_{\xi \in \Lambda} \mathbb{P}_{\neq 0} (\partial_{t} a_{(\xi)}^{2} \phi_{(\xi)}^{2} \psi_{(\xi)}^{2} \xi). 
\end{align*} 
Quite naturally from \eqref{estimate 53}, we define again $\mathring{R}_{q+1} \triangleq R_{\text{linear}} + R_{\text{oscillation}} + R_{\text{corrector}} + R_{\text{commutator 1}} + R_{\text{commutator 2}}$. We observe that because there is only $a_{(\xi)}$, and not $\bar{a}_{(\xi)}$, in $R_{\text{oscillation}}$, our estimates \eqref{[Equation (4.48m), HZZ19]}-\eqref{[Equation (4.48n), HZZ19]} directly apply to our case as our choices of $\alpha$ and $l$ are same as those in the proof of Theorem \ref{[Theorem 1.1, HZZ19]}. Now with our same choice of $p^{\ast}$ from \eqref{p}, we compute from \eqref{[Equation (6.33b), HZZ19]} using \eqref{[Equation (6.7b), HZZ19]} and \eqref{estimate 39}
\begin{equation}\label{estimate 45}
 \lVert R_{\text{linear}} \rVert_{C_{t}L_{x}^{p^{\ast}}} \lesssim \lVert w_{q+1} \rVert_{C_{t}W_{x}^{2m-1, p^{\ast}}} + \sum_{\xi \in \Lambda} \lVert \partial_{t} \text{curl} (\bar{a}_{(\xi)} V_{(\xi)}) \rVert_{C_{t}L_{x}^{p^{\ast}}} + m_{L}^{3} M_{0}(t)^{\frac{1}{2}} \lVert w_{q+1} \rVert_{C_{t}L_{x}^{p^{\ast}}}. 
\end{equation}
We can estimate similarly to \eqref{estimate 41} and \eqref{estimate 31} (see also \eqref{estimate 46}) 
\begin{subequations}
\begin{align}
& \lVert w_{q+1} \rVert_{C_{t}W_{x}^{2m-1, p^{\ast}}}  \overset{\eqref{[Equation (4.16), HZZ19]} \eqref{[Equation (6.23b), HZZ19]} \eqref{estimate 40} }{\lesssim}  m_{L} M_{0}(t)^{\frac{1}{2}} \lambda_{q+1}^{- \frac{61 \alpha}{6}} + M_{0}(t) \lambda_{q+1}^{\frac{12 m -15 - 148\alpha}{24}}, \label{estimate 42}\\
& \sum_{\xi \in \Lambda} \lVert \partial_{t} \text{curl} (\bar{a}_{(\xi)} V_{(\xi)} ) \rVert_{C_{t}L_{x}^{p^{\ast}}} \overset{\eqref{[Equation (4.16), HZZ19]} \eqref{[Equation (6.23), HZZ19]} \eqref{[Equation (B.7c), HZZ19]}  }{\lesssim}  m_{L} M_{0}(t)^{\frac{1}{2}} \lambda_{q+1}^{-\frac{\alpha}{6}} + m_{L} M_{0}(t)^{\frac{1}{2}} \lambda_{q+1}^{\frac{59 \alpha - 12m}{6}}, \label{estimate 43}\\
& m_{L}^{3} M_{0}(t)^{\frac{1}{2}} \lVert w_{q+1} \rVert_{C_{t}L_{x}^{p^{\ast}}} \overset{\eqref{[Equation (6.23b), HZZ19]}  \eqref{estimate 50}  }{\lesssim} m_{L}^{4} M_{0}(t) l^{-2} r_{\bot}^{\frac{2}{p^{\ast}} -1} r_{\lVert}^{\frac{1}{p^{\ast}} - \frac{1}{2}}  \lesssim m_{L}^{4} M_{0}(t) \lambda_{q+1}^{\frac{ -122 \alpha - 24m + 12}{12}}  \label{estimate 44}
\end{align}
\end{subequations}
where the last inequality made use of the assumption that $b > \frac{16}{\alpha}$. We apply \eqref{estimate 42}, \eqref{estimate 43}, \eqref{estimate 44} to \eqref{estimate 45}, observe that $1 < \delta_{q+2} \lambda_{q+1}^{\frac{\alpha}{8}}$ due to $\alpha > 16 \beta b$ and $\alpha$ from \eqref{alpha} satisfies
\begin{align*}
&(- \frac{61\alpha}{6} + \frac{\alpha}{8}) \vee (\frac{12 m - 15 - 145 \alpha}{24}) \vee (- \frac{\alpha}{6} + \frac{\alpha}{8}) \\
& \hspace{5mm} \vee (\frac{59 \alpha - 12m}{6} + \frac{\alpha}{8}) \vee (\frac{-122 \alpha - 24 m + 12}{12} + \frac{\alpha}{8}) < 0; 
\end{align*} 
thus, taking $a > 0$ sufficiently large give us 
\begin{align}\label{[Equation (6.33f), HZZ19]}
\lVert R_{\text{linear}} \rVert_{C_{t}L_{x}^{p^{\ast}}} \leq (2\pi)^{-3 (\frac{p^{\ast}-1}{p^{\ast}})} \frac{M_{0}(t)c_{R} \delta_{q+2}}{5}. 
\end{align} 
Next, concerning $R_{\text{corrector}}$ from \eqref{[Equation (6.33b), HZZ19]} we first observe that because $\alpha$ from \eqref{alpha} satisfies $\alpha < \frac{25 - 20m}{480} \wedge \frac{5-4m}{32} $, by taking $a > 0$ sufficiently large we can estimate  
\begin{align}
&\lVert w_{q+1} \rVert_{C_{t}L_{x}^{2p^{\ast}}} + \lVert w_{q+1}^{(p)} \rVert_{C_{t}L_{x}^{2p^{\ast}}}  \overset{\eqref{[Equation (4.16), HZZ19]} \eqref{[Equation (6.25), HZZ19]}  \eqref{[Equation (6.26), HZZ19]} }{\lesssim} m_{L} M_{0}(t)^{\frac{1}{2}} \delta_{q+1}^{\frac{1}{2}} \lambda_{q+1}^{\frac{5- 4m - 74 \alpha}{24}}\label{estimate 51} \\
& \hspace{32mm} \times [ 1+ \lambda_{q+1}^{\frac{480 \alpha + 20m - 25}{24}} + M_{0}(t)^{\frac{1}{2}} \delta_{q+1}^{\frac{1}{2}} \lambda_{q+1}^{\frac{32 \alpha + 4 m- 5}{8}}] \lesssim m_{L} M_{0}(t)^{\frac{1}{2}} \delta_{q+1}^{\frac{1}{2}} \lambda_{q+1}^{\frac{5-4m-74\alpha}{24}}.  \nonumber 
 \end{align}  
Using that $\alpha < \frac{15 - 12m}{335}$, this allows us to estimate from \eqref{[Equation (6.33b), HZZ19]} 
\begin{align}
&\lVert R_{\text{corrector}} \rVert_{C_{t}L_{x}^{p^{\ast}}}  \overset{\eqref{[Equation (6.5), HZZ19]}}{\lesssim} m_{L}^{2} ( \lVert w_{q+1}^{(c)} \rVert_{C_{t}L_{x}^{2p^{\ast}}} + \lVert w_{q+1}^{(t)} \rVert_{C_{t}L_{x}^{2p^{\ast}}}) (\lVert w_{q+1} \rVert_{C_{t}L_{x}^{2p^{\ast}}} + \lVert w_{q+1}^{(p)} \rVert_{C_{t}L_{x}^{2p^{\ast}}})  \label{[Equation (6.33g), HZZ19]}\\
&\hspace{15mm} \overset{\eqref{[Equation (4.16), HZZ19]}  \eqref{[Equation (6.26), HZZ19]} }{\lesssim}  m_{L}^{4} M_{0}(t) \delta_{q+1}^{\frac{1}{2}} \lambda_{q+1}^{\frac{ 5- 4m - 74\alpha}{24}} ( \delta_{q+1}^{\frac{1}{2}} \lambda_{q+1}^{\frac{203 \alpha + 8m - 10}{12}} + \delta_{q+1}M_{0}(t)^{\frac{1}{2}} \lambda_{q+1}^{\frac{11 \alpha + 4m -5}{12}})  \nonumber\\
& \hspace{15mm} \leq M_{0}(t)\delta_{q+2} [m_{L}^{4} \lambda_{q+1}^{\frac{335 \alpha + 12m - 15}{24}} + m_{L}^{4} M_{0}(t)^{\frac{1}{2}} \lambda_{q+1}^{ \frac{ - 49 \alpha + 4m-5}{24}}]  \leq (2\pi)^{-3 (\frac{p^{\ast}-1}{p^{\ast}})} \frac{M_{0}(t) c_{R} \delta_{q+2}}{5}. \nonumber
\end{align} 
Next, for any $t \in [0, T_{L}]$, we can use the fact that $l^{\frac{1}{2} - 2 \delta} \lambda_{q}^{4} \ll \delta_{q+2} \lambda_{q}^{-\frac{8}{3}}$ due to the assumptions of $\delta \in (0, \frac{1}{12}), \alpha > 16 \beta b$ and $\alpha b > 16$, and estimate from \eqref{[Equation (6.16d), HZZ19]}  
\begin{equation}\label{[Equation (6.33h), HZZ19]}
 \lVert R_{\text{commutator 1}} \rVert_{C_{t}L_{x}^{1}} \overset{\eqref{[Equation (6.5), HZZ19]}\eqref{[Equation (6.7a), HZZ19]} \eqref{[Equation (6.7b), HZZ19]} }{\lesssim}  m_{L}^{4} l^{\frac{1}{2} - 2 \delta} M_{0}(t) \lambda_{q}^{4} \leq \frac{c_{R} M_{0}(t)}{5} \delta_{q+2}. 
\end{equation} 
Moreover, we know $\lvert \Upsilon_{l}(t) - \Upsilon (t) \rvert  \leq l^{\frac{1}{2} - 2 \delta} m_{L}^{2}$ due to \eqref{[Equation (6.5), HZZ19]}. Therefore, it follows from \eqref{[Equation (6.33b), HZZ19]} that using again $l^{\frac{1}{2} - 2\delta} \ll \delta_{q+2} \lambda_{q}^{-\frac{20}{3}}$, 
\begin{equation}\label{[Equation (6.34a), HZZ19]}
 \lVert R_{\text{commutator 2}} \rVert_{C_{t}L_{x}^{1}} \lesssim l^{\frac{1}{2} - 2\delta} m_{L}^{2} \lVert v_{q+1} \rVert_{C_{t}L_{x}^{2}}^{2} \overset{\eqref{[Equation (6.18), HZZ19]} \eqref{[Equation (6.23b), HZZ19]}  \eqref{[Equation (6.28), HZZ19]}}{\lesssim} l^{\frac{1}{2} - 2 \delta} m_{L}^{4} M_{0}(t) \leq \frac{M_{0}(t) c_{R} \delta_{q+2}}{5}. 
\end{equation} 
Collecting the estimates \eqref{[Equation (6.33f), HZZ19]}, \eqref{[Equation (6.33g), HZZ19]}-\eqref{[Equation (6.34a), HZZ19]}, \eqref{[Equation (4.48m), HZZ19]}-\eqref{[Equation (4.48n), HZZ19]} implies that \eqref{[Equation (6.7c), HZZ19]} is satisfied at level $q+1$ identically as we showed in \eqref{[Equation (4.48q), HZZ19]}. At last, essentially identical arguments at the end of the proof of Proposition \ref{[Proposition 4.2, HZZ19]} show that  $(v_{q}, \mathring{R}_{q})$ being $(\mathcal{F}_{t})_{t\geq 0}$-adapted leads to $(v_{q+1}, \mathring{R}_{q+1})$ being $(\mathcal{F}_{t})_{t\geq 0}$-adapted, and that $v_{q}(0,x), \mathring{R}_{q}(0,x)$ being deterministic implies $v_{q+1}(0,x), \mathring{R}_{q+1}(0,x)$ being deterministic. 

\section{Appendix}
In this Appendix, for convenience of readers we collect results which have been used throughout this manuscript, as well as proofs which we chose to postpone. 

\subsection{Past results}
The following important definition, along with proofs of its properties, was introduced initially in \cite[Definition 4.2, Lemma 4.3]{DS13} and has proven to be useful (e.g., \cite[Definition 1.4, Lemma 1.5]{BDIS15}, \cite[Lemma 5]{LT20}, \cite[Equation (2.26)]{BCV18}): 
\begin{lemma}\label{divergence inverse operator}
\rm{(\cite[Equation (5.34)]{BV19b})} For any $v \in C^{\infty}(\mathbb{T}^{3})$ that has mean zero, define 
\begin{equation}\label{estimate 5}
(\mathcal{R}v)_{kl} \triangleq ( \partial_{k}\Delta^{-1} v^{l} + \partial_{l} \Delta^{-1} v^{k}) - \frac{1}{2} (\delta_{kl} + \partial_{k} \partial_{l} \Delta^{-1}) \text{div} \Delta^{-1} v
\end{equation} 
for $k, l \in \{1,2,3\}$. Then $\mathcal{R} v(x)$ is a symmetric trace-free matrix for each $x \in \mathbb{T}^{3}$, that  satisfies $\text{div} (\mathcal{R} v) = v$. When $v$ is not mean zero, we overload the notation and denote by $\mathcal{R}v \triangleq \mathcal{R} (v - \int_{\mathbb{T}^{3}}  vdx )$. Moreover, $\mathcal{R}$ satisfies the classical Calder$\acute{\mathrm{o}}$n-Zygmund and Schauder estimates: $\lVert (-\Delta)^{\frac{1}{2}} \mathcal{R} \rVert_{L_{x}^{p} \mapsto L_{x}^{p}} + \lVert \mathcal{R} \rVert_{L_{x}^{p} \mapsto L_{x}^{p}}  + \lVert \mathcal{R} \rVert_{C_{x} \mapsto C_{x}} \lesssim 1$ for all $p \in (1, \infty)$. 
\end{lemma} 

\begin{lemma}\label{[Lemma 7.4, BV19b]}
\rm{(\cite[Lemma 7.4]{BV19b})} Fix integers $N, \kappa \geq 1$ and let $\zeta > 1$ satisfy 
\begin{equation}\label{[Equation (7.41), BV19b]}
\frac{ 2 \pi \sqrt{3} \zeta}{\kappa} \leq \frac{1}{3} \text{ and } \zeta^{4} \frac{ (2 \pi \sqrt{3} \zeta)^{N}}{\kappa^{N}} \leq 1. 
\end{equation} 
Let $p \in \{1, 2\}$ and $f$ be a $\mathbb{T}^{3}$-periodic function such that there exists a constant $C_{f} > 0$ so that $\lVert D^{j} f \rVert_{L_{x}^{p}} \leq C_{f} \zeta^{j}$ for all $0 \leq j \leq N + 4$. In addition, let $g$ be a $(\mathbb{T}/\kappa)^{3}$-periodic function. Then $\lVert f g \rVert_{L_{x}^{p}} \lesssim C_{f} \lVert g \rVert_{L_{x}^{p}}$  
where the implicit constant is universal. 
\end{lemma}

\begin{lemma}\label{[Lemma 7.5, BV19b]}
\rm{(\cite[Lemma 7.5]{BV19b})} Fix parameters $1 \leq \zeta < \kappa$, $p \in (1, 2]$ and assume that there exists an $N \in \mathbb{N}$ such that $\zeta^{N}\leq \kappa^{N-2}$. Let $a \in C^{N} (\mathbb{T}^{3})$ be such that there exists $C_{a} > 0$ that  satisfies $\lVert D^{j} a \rVert_{C_{x}} \leq C_{a} \zeta^{j}$ for all $0 \leq j \leq N$. Assume furthermore that $f \in L^{p}(\mathbb{T}^{3})$ satisfies $\int_{\mathbb{T}^{3}} a(x) \mathbb{P}_{\geq \kappa} f(x) dx= 0$. Then we have 
\begin{equation*}
\lVert (-\Delta)^{\frac{1}{2}}  (a \mathbb{P}_{\geq \kappa} f) \rVert_{L_{x}^{p}} \lesssim_{p, N} C_{a} \frac{\lVert f\rVert_{L_{x}^{p}}}{\kappa}.
\end{equation*} 
\end{lemma}

The following are some notations and results concerning intermittent jets from \cite[Appendix B]{HZZ19}, originally from \cite[Section 7.4]{BV19b} (also \cite[Section 4]{BCV18}). Lemma \ref{[Lemma B.1, HZZ19]} is a geometric lemma in variation of \cite[Lemma 3.2]{DS13}, which proved to be useful on many other occasions (e.g., \cite[Lemma 1.3]{BDIS15}, \cite[Proposition 3.2]{BV19a}).
\begin{lemma}\label{[Lemma B.1, HZZ19]}
\rm{(\cite[Lemma 6.6]{BV19b})}
Let $\overline{B_{\frac{1}{2}}(\text{Id})}$ denote the closed ball of radius $\frac{1}{2}$ around an identity matrix in the space of $3\times 3$ symmetric matrices. Then there exists $\Lambda \subset \mathbb{S}^{2} \cap \mathbb{Q}^{3}$ such that for each $\xi \in \Lambda$, there exist $C^{\infty}$ smooth functions $\gamma_{\xi}: \hspace{0.5mm}  B_{\frac{1}{2}} (\text{Id}) \mapsto \mathbb{R}$ which obey 
\begin{equation}\label{[Equation (A.9d), HZZ19]}
R = \sum_{\xi \in \Lambda} \gamma_{\xi}^{2} (R) (\xi \otimes \xi) 
\end{equation}
for every symmetric matrix $R$ that satisfies $\lvert R - \text{Id} \rvert \leq \frac{1}{2}$. 
\end{lemma} 
We note that the precise statement from \cite[Lemma 6.6]{BV19b} consists of mutually disjoint sets $\{\Lambda_{i}\}_{i=0,1} \subset \mathbb{S}^{2} \cap \mathbb{Q}^{3}$; however, as pointed out in \cite[Section 7.4]{BV19b}, one can just choose one of $\Lambda_{0}$ or $\Lambda_{1}$ and relabel it as $\Lambda$ in case of the NS equations. Define a constant 
\begin{equation}\label{[Equation (A.9e), HZZ19]}
M \triangleq C_{\Lambda} \sup_{\xi \in \Lambda}(\lVert \gamma_{\xi} \rVert_{C^{0}} + \lVert \nabla \gamma_{\xi} \rVert_{C^{0}}) \text{ where } C_{\Lambda} \triangleq 8 \lvert \Lambda \rvert (1+ 8\pi^{3})^{\frac{1}{2}} 
\end{equation}
and $\lvert \Lambda \rvert$ denotes the cardinality of the set $\Lambda$. For every $\xi \in \Lambda$, let $A_{\xi} \in \mathbb{S}^{2} \cap \mathbb{Q}^{3}$ be an orthogonal vector to $\xi$. It follows that for each $\xi \in \Lambda$, $\{\xi, A_{\xi}, \xi \times A_{\xi} \} \subset \mathbb{S}^{2} \cap \mathbb{Q}^{3}$ forms an orthonormal basis for $\mathbb{R}^{3}$. Furthermore, we denote by $n_{\ast}$ the smallest natural number such that $\{n_{\ast}, \xi, n_{\ast}A_{\xi}, n_{\ast} \xi \times A_{\xi} \} \subset \mathbb{Z}^{3}$ for every $\xi \in \Lambda$. 

Now let $\Phi: \hspace{0.5mm}  \mathbb{R}^{2} \mapsto \mathbb{R}^{2}$ be a smooth function with support contained in a ball of radius one. We normalize $\Phi$ so that $\phi \triangleq - \Delta \Phi$ obeys 
\begin{equation}\label{[Equation (B.1), HZZ19]}
 \int_{\mathbb{R}^{2}} \phi^{2} ( x_{1}, x_{2}) dx_{1} dx_{2} = 4\pi^{2}. 
\end{equation} 
It follows that $\phi$ has mean zero. We define $\psi: \hspace{0.5mm}  \mathbb{R} \mapsto \mathbb{R}$ to be a smooth, mean zero function with support in the ball of radius one such that 
\begin{equation}\label{[Equation (B.2), HZZ19]}
 \int_{\mathbb{R}}\psi^{2} (x_{3}) dx_{3} = 2\pi. 
\end{equation} 
Let $\phi_{r_{\bot}}, \Phi_{r_{\bot}}$ and $\psi_{r_{\lVert}}$ be the rescaled cutoff functions 
\begin{equation}\label{[Equation (B.2b), HZZ19]}
\phi_{r_{\bot}} (x_{1}, x_{2}) \triangleq \frac{ \phi( \frac{x_{1}}{r_{\bot}}, \frac{x_{2}}{r_{\bot}})}{r_{\bot}}, \Phi_{r_{\bot}} (x_{1}, x_{2}) \triangleq \frac{ \Phi( \frac{x_{1}}{r_{\bot}}, \frac{x_{2}}{r_{\bot}})}{r_{\bot}} \text{ and } \psi_{r_{\lVert}} (x_{3}) \triangleq \frac{ \psi( \frac{x_{3}}{r_{\lVert }})}{r_{\lVert }^{\frac{1}{2}}}
\end{equation} 
so that $\phi_{r_{\bot}} = - r_{\bot}^{2} \Delta \Phi_{r_{\bot}}$ in which we will assume $r_{\bot}, r_{\lVert} > 0$ to satisfy $r_{\bot} \ll r_{\lVert} \ll 1$ and $r_{\bot}^{-1} \ll \lambda$ from \eqref{[Equation (B.2a), HZZ19]} where the last assumption $r_{\bot}^{-1} \ll \lambda$ will be used for derivative estimates of \eqref{[Equation (B.7a), HZZ19]}-\eqref{[Equation (B.7c), HZZ19]}. By an abuse of notation, we periodize $\phi_{r_{\bot}}, \Phi_{r_{\bot}}$ and $\psi_{r_{\lVert}}$ so that they are treated as functions defined on $\mathbb{T}^{2}, \mathbb{T}^{2}$ and $\mathbb{T}$, respectively. For a large real number $\lambda$ such that $\lambda r_{\bot} \in \mathbb{N}$, and a large time oscillation parameter $\mu > 0$, for every $\xi \in \Lambda$ we introduce 
\begin{subequations}\label{[Equation (B.2c), HZZ19]}
\begin{align}
& \psi_{(\xi)} (t,x) \triangleq \psi_{\xi, r_{\bot}, r_{\lVert}, \lambda, \mu} (t,x) \triangleq \psi_{r_{\lVert}} (n_{\ast} r_{\bot} \lambda(x \cdot \xi + \mu t)), \\
& \Phi_{(\xi)} (x) \triangleq \Phi_{\xi, r_{\bot}, \lambda} (x) \triangleq \Phi_{r_{\bot}} (n_{\ast} r_{\bot} \lambda (x - a_{\xi}) \cdot A_{\xi}, n_{\ast} r_{\bot} \lambda (x- a_{\xi}) \cdot (\xi \times A_{\xi})), \\
& \phi_{(\xi)} (x) \triangleq \phi_{\xi, r_{\bot}, \lambda} (x) \triangleq \phi_{r_{\bot}} (n_{\ast} r_{\bot} \lambda (x- a_{\xi}) \cdot A_{\xi}, n_{\ast} r_{\bot} \lambda (x- a_{\xi}) \cdot (\xi \times A_{\xi} )), 
\end{align} 
\end{subequations}
where $a_{\xi} \in \mathbb{R}^{3}$ are shifts which ensure that the functions $\{ \Phi_{(\xi)}\}_{\xi \in \Lambda}$ have mutually disjoint support. In order to ensure that such shifts exist, it suffices to require $r_{\bot}$ to be sufficiently small depending on $\Lambda$. We can now define intermittent jets $W_{(\xi)}: \hspace{0.5mm}  \mathbb{T}^{3} \times \mathbb{R} \mapsto \mathbb{R}^{3}$ by 
\begin{equation}\label{[Equation (B.3), HZZ19]}
W_{(\xi)} (t,x) \triangleq W_{\xi, r_{\bot}, r_{\lVert}, \lambda, \mu} (t,x) \triangleq \xi \psi_{(\xi)} (t,x) \phi_{(\xi)} (x). 
\end{equation} 
It follows that $W_{(\xi)}$ has mean zero, it is $(\mathbb{T}/r_{\bot}\lambda)^{3}$-periodic, and 
\begin{equation}\label{[Equation (B.4), HZZ19]}
W_{(\xi)} \otimes W_{(\xi')} = 0 \hspace{3mm} \forall \hspace{1mm} \xi, \xi' \in \Lambda \text{ such that } \xi \neq \xi'.
\end{equation} 
Due to \eqref{[Equation (B.1), HZZ19]}-\eqref{[Equation (B.2b), HZZ19]} we also have 
\begin{equation}\label{[Equation (B.4a), HZZ19]}
\fint_{\mathbb{T}^{3}} W_{(\xi)} (t,x) \otimes W_{(\xi)} (t,x) dx = \xi \otimes \xi. 
\end{equation} 
Lemma \ref{[Lemma B.1, HZZ19]} and \eqref{[Equation (B.4), HZZ19]} imply 
\begin{equation}\label{[Equation (B.5), HZZ19]}
\sum_{\xi \in \Lambda} \gamma_{\xi}^{2}(R) \fint_{\mathbb{T}^{3}} W_{(\xi)} (t,x) \otimes W_{(\xi)} (t,x) dx = R
\end{equation} 
for every symmetric matrix $R$ such that  $\lvert R - \text{Id} \rvert \leq \frac{1}{2}$. The following identity also holds: 
\begin{equation}\label{estimate 7}
\text{div} (W_{(\xi)} \otimes W_{(\xi)}) =  \mu^{-1} \partial_{t} (\phi_{(\xi)}^{2} \psi_{(\xi)}^{2} \xi). 
\end{equation} 
Although $W_{(\xi)}$ is not divergence-free, assuming $r_{\bot} \ll r_{\lVert}$ as we did in \eqref{[Equation (B.2a), HZZ19]}, we can define 
\begin{equation}\label{[Equation (B.6), HZZ19]}
W_{(\xi)}^{(c)} \triangleq \frac{ \nabla \psi_{(\xi)} }{n_{\ast}^{2} \lambda^{2}} \times \text{curl} (\Phi_{(\xi)} \xi) = \text{curl curl} V_{(\xi)} - W_{(\xi)} \text{ with } V_{(\xi)} (t,x) \triangleq \frac{ \xi \psi_{(\xi)} (t,x)}{n_{\ast}^{2} \lambda^{2}} \Phi_{(\xi)} (x), 
\end{equation} 
from which it follows that $\text{div} (W_{(\xi)} + W_{(\xi)}^{(c)}) = 0$. Finally, it was shown in \cite[Section 7.4]{BV19b} that for $N, M \geq 0$ and $p \in [1, \infty]$, 
\begin{subequations}\label{[Equation (B.7), HZZ19]}
\begin{align}
 \lVert \nabla^{N} \partial_{t}^{M} \psi_{(\xi)} \rVert_{L^{p}} \lesssim& r_{\lVert}^{\frac{1}{p} - \frac{1}{2}} \left( \frac{r_{\bot} \lambda}{r_{\lVert}} \right)^{N} \left( \frac{ r_{\bot} \lambda \mu}{r_{\lVert}}\right)^{M}, \label{[Equation (B.7a), HZZ19]}\\
 \lVert \nabla^{N} \phi_{(\xi)} \rVert_{L^{p}} + \lVert \nabla^{N} \Phi_{(\xi)} \rVert_{L^{p}} \lesssim& r_{\bot}^{\frac{2}{p} - 1} \lambda^{N}, \label{[Equation (B.7b), HZZ19]}\\
 \lVert \nabla^{N} \partial_{t}^{M} W_{(\xi)} \rVert_{L^{p}} + \frac{r_{\lVert}}{r_{\bot}} \lVert \nabla^{N} \partial_{t}^{M} W_{(\xi)}^{(c)} \rVert_{L^{p}} + \lambda^{2} \lVert \nabla^{N} \partial_{t}^{M} V_{(\xi)} \rVert_{L^{p}} \lesssim& r_{\bot}^{\frac{2}{p} - 1} r_{\lVert}^{\frac{1}{p} - \frac{1}{2}} \lambda^{N} \left( \frac{r_{\bot} \lambda \mu}{r_{\lVert} } \right)^{M}, \label{[Equation (B.7c), HZZ19]}
\end{align} 
\end{subequations} 
where it should be emphasized that the implicit constants are independent of $\lambda, r_{\bot}, r_{\lVert}$ and $\mu$. For the derivative estimates, the assumption of $r_{\bot}^{-1} \ll \lambda$ of \eqref{[Equation (B.2a), HZZ19]} was used.

\subsection{Proof of Proposition \ref{[Theorem 3.1, HZZ19]}}

We refer to the proof of \cite[Theorem 4.1]{FR08} (and \cite[Theorem 3.1]{Y19}) for the proof of the existence of a martingale solution. The proof of the following result from \cite{HZZ19} does not depend on the explicit form of the diffusive term and thus applies directly to our case. 
\begin{lemma}\label{[Lemma A.1, HZZ19]} 
\rm{ (\cite[Lemma A.1]{HZZ19})} Let $\{(s_{n}, \xi_{n}) \}_{n\in\mathbb{N}} \subset [0,\infty) \times L_{\sigma}^{2}$  satisfy $\lim_{n\to\infty} \lVert (s_{n}, \xi_{n}) - (s, \xi^{\text{in}}) \rVert_{\mathbb{R} \times L_{x}^{2}} = 0$ and $\{P_{n}\}_{n\in\mathbb{N}}$ be a family of probability measures on $\Omega_{0}$ satisfying for all $n \in \mathbb{N}$, $P_{n}(\{\xi(t) = \xi_{n} \hspace{1mm} \forall \hspace{1mm} t \in [0, s_{n}]\}) = 1$ and for some $\gamma, \kappa > 0$ and any $T > 0$,  
\begin{equation}\label{[Equation (A.2), HZZ19]} 
\sup_{n\in\mathbb{N}} \mathbb{E}^{P_{n}} [ \lVert \xi \rVert_{C([0,T]; L_{x}^{2})} + \sup_{r, t \in [0,T]: \hspace{0.5mm}  r \neq t} \frac{ \lVert \xi(t) - \xi(r) \rVert_{H_{x}^{-3}}}{\lvert t- r \rvert^{\kappa}} + \lVert \xi \rVert_{L^{2}([s_{n}, T]; H_{x}^{\gamma})}^{2}] < \infty.  
\end{equation} 
Then $\{P_{n}\}_{n\in\mathbb{N}}$ is tight in $\mathbb{S} \triangleq C_{\text{loc}} ([0,\infty); H^{-3}(\mathbb{T}^{3})) \cap L^{2}_{\text{loc}}(0,\infty; L_{\sigma}^{2})$. 
\end{lemma}
We fix $\{P_{n}\} \subset \mathcal{C} (s_{n}, \xi_{n}, \{C_{t,q}\}_{q\in\mathbb{N}, t \geq s_{n}})$ and define $F(\xi) \triangleq - \mathbb{P} \text{div} (\xi \otimes \xi) - (-\Delta)^{m} \xi$. By (M2) of Definition \ref{[Definition 3.1, HZZ19]} we know that for all $n \in \mathbb{N}$ and $t \in [s_{n}, \infty)$, 
\begin{equation}\label{[Equation (A.4b), HZZ19]}
\xi(t) = \xi_{n} + \int_{s_{n}}^{t} F(\xi(r)) dr + M_{t, s_{n}}^{\xi} \hspace{3mm} P_{n}\text{-a.s.},
\end{equation}
where a mapping $t\mapsto M_{t,s_{n}}^{\xi, i} \triangleq \langle M_{t, s_{n}}^{\xi}, e_{i} \rangle$, $\xi \in \Omega_{0}$, is a continuous, square-integrable martingale w.r.t. $(\mathcal{B}_{t})_{t\geq s_{n}}$ with $\langle \langle M_{t, s_{n}}^{\xi, i} \rangle \rangle = \int_{s_{n}}^{t} \lVert G(\xi(r))^{\ast} e_{i} \rVert_{U}^{2} dr$. We compute for any $p \in (1, \infty)$, 
\begin{align} 
\mathbb{E}^{p_{n}} [ \sup_{r, t \in [s_{n},T]: r \neq t} \frac{ \lVert \int_{r}^{t} F(\xi(l)) dl \rVert_{H_{x}^{-3}}^{p}}{\lvert t- r \rvert^{p-1}}]  \overset{\text{(M3)}}{\lesssim_{T,p}} (1+ \lVert \xi_{n} \rVert_{L_{x}^{2}}^{2p})  \label{[Equation (A.4c), HZZ19]}
\end{align} 
where we used an estimate of $\lVert \xi \rVert_{H_{x}^{2m-3}} \lesssim 1+ \lVert \xi \rVert_{L_{x}^{2}}^{2}$ for the diffusive term, and therefore the implicit constant is independent of $n$. On the other hand, making use of \eqref{[Equation (1.7c), HZZ19]}, (M2), (M3) and Kolmogorov's test (e.g., \cite[Theorem 3.3]{DZ14}) gives for any $\alpha \in (0, \frac{p-1}{2p})$, 
\begin{equation}\label{[Equation (A.4e), HZZ19]}
\mathbb{E}^{P_{n}} [ \sup_{r, t \in [0, T]: r \neq t} \frac{ \lVert M_{t,s_{n}}^{\xi} - M_{r, s_{n}}^{\xi} \rVert_{L_{x}^{2}}}{\lvert t-r \rvert^{\alpha}} ] \lesssim_{p} C_{t,p}(1+ \lVert \xi_{n} \rVert_{L_{x}^{2}}^{2p}). 
\end{equation} 
Making use of \eqref{[Equation (A.4b), HZZ19]}-\eqref{[Equation (A.4e), HZZ19]} leads to for all $\kappa \in (0, \frac{1}{2})$, 
\begin{equation}\label{[Equation (A.5), HZZ19]}
\sup_{n \in \mathbb{N}} \mathbb{E}^{P_{n}} [ \sup_{r, t \in [0,T]: r \neq t} \frac{ \lVert \xi(t) - \xi(r) \rVert_{H_{x}^{-3}}}{\lvert t-r \rvert^{\kappa}} ] < \infty. 
\end{equation} 
Due to (M1), \eqref{[Equation (1.7h), HZZ19]} and \eqref{[Equation (A.5), HZZ19]}, we are able to apply Lemma \ref{[Lemma A.1, HZZ19]} and deduce that $\{P_{n}\}$ is tight in $\mathbb{S}$. By Prokhorov's and Skorokhod's representation theorems we deduce that $P_{n}$ converges weakly to some $P \in \mathcal{P}(\Omega_{0})$ and that there exists a probability space $(\tilde{\Omega}, \tilde{\mathcal{F}}, \tilde{P})$ and $\mathbb{S}$-valued random variables $\{ \tilde{\xi}_{n}\}_{n\in\mathbb{N}}$ and $\tilde{\xi}$ such that 
\begin{equation}\label{[(i) and (ii), HZZ19]}
\tilde{\xi}_{n} \text{ has the law } P_{n} \hspace{1mm} \forall \hspace{1mm} n \in \mathbb{N}, \tilde{\xi}_{n} \to \tilde{\xi} \text{ in } \mathbb{S} \hspace{1mm}  \tilde{P}\text{-a.s. and }\tilde{\xi} \text{ has the law } P. 
\end{equation} 
Making use of \eqref{[(i) and (ii), HZZ19]} and (M1) for $P_{n}$ immediately leads to 
\begin{equation}\label{[Equation (A.5a), HZZ19]}
P( \{ \xi(t) = \xi^{\text{in}} \hspace{1mm} \forall \hspace{1mm} t \in [0,s] \}) = 1,
\end{equation} 
which implies (M1) for $P$. Next, it follows that for every $e_{i} \in C^{\infty}(\mathbb{T}^{3})$, $\tilde{P}$-a.s. 
\begin{equation}\label{[Equation (A.5b), HZZ19]}
\langle \tilde{\xi}_{n}(t), e_{i} \rangle \to \langle \tilde{\xi}(t), e_{i} \rangle, \hspace{3mM} \int_{s_{n}}^{t} \langle F(\tilde{\xi}_{n}(r)), e_{i} \rangle dr \to \int_{s}^{t} \langle F(\tilde{\xi} (r)), e_{i} \rangle dr.
\end{equation} 
The first convergence follows from (2) as in \cite{HZZ19}. Concerning the second convergence, we can rely on the convergence of  
\begin{equation*}
 \mathbb{E}^{\tilde{P}} [ \int_{s_{n}}^{s} \langle - (-\Delta)^{m} \tilde{\xi}_{n}, e_{i} \rangle dr]  \to 0 \text{ and }  \mathbb{E}^{\tilde{P}} [ \int_{s}^{t} \langle (-\Delta)^{m}( \tilde{\xi}_{n} - \tilde{\xi}), e_{i} \rangle dr]  \to 0
\end{equation*} 
as $n\to\infty$ to handle the diffusive term within $F$. Hereafter, one can trace the proof of \cite[Theorem 3.1]{HZZ19} identically to verify that $P$ satisfies (M2) and (M3) and conclude that $P \in \mathcal{C} ( s, \xi_{0}, \{C_{t,q} \}_{q \in \mathbb{N}, t \geq s })$ as the specific form of the diffusive term is never used; thus, we omit them and refer to \cite{HZZ19}. 

\subsection{Continuation of the proof of Proposition \ref{[Proposition 3.7, HZZ19]}}

By Theorem \ref{[Theorem 1.1, HZZ19]} we know that there exists a process $u$ that is a weak solution to \eqref{[Equation (1.4), HZZ19]} on $[0, T_{L}]$ such that $u(\cdot \wedge T_{L}) \in \Omega_{0}$ $\textbf{P}$-almost surely. By \eqref{[Equation (3.13a), HZZ19]}, \eqref{[Equation (3.12), HZZ19]} and \eqref{[Equation (1.4), HZZ19]} we can deduce 
\begin{equation}\label{[Equation (3.16), HZZ19]}
Z^{u}(t) = z(t) \hspace{1mm} \forall \hspace{1mm} t \in [0, T_{L}] \hspace{1mm} \textbf{P}\text{-almost surely.}
\end{equation} 
By Proposition \ref{[Proposition 3.6, HZZ19]} we know that the trajectories $t \mapsto \lVert z(t) \rVert_{H_{x}^{\frac{5+\sigma}{2}}}$ and $t \mapsto \lVert z \rVert_{C_{t}^{\frac{2}{5} -  \delta} H_{x}^{\frac{3+ \sigma}{2}}}$ for any $\delta \in (0, \frac{2}{5})$ are both $\textbf{P}$-a.s. continuous. It follows immediately from \eqref{[Equation (3.16), HZZ19]} and \eqref{[Equation (3.13b), HZZ19]} that $\tau_{L}(u) \leq T_{L}$ $\textbf{P}$-a.s. while from \eqref{[Equation (4.2), HZZ19]} that $T_{L} \leq \tau_{L}(u)$ $\textbf{P}$-a.s. and hence 
\begin{equation}\label{[Equation (3.15), HZZ19]}
\tau_{L}(u) = T_{L} \hspace{1mm} \textbf{P}\text{-almost surely}. 
\end{equation}
Next we verify that $P$ is a martingale solution to \eqref{[Equation (1.4), HZZ19]} on $[0, \tau_{L}]$. By Theorem \ref{[Theorem 1.1, HZZ19]} $u^{\text{in}}$ is deterministic; along with \eqref{[Equation (1.7c), HZZ19]} and \eqref{[Equation (1.5), HZZ19]}, this shows that (M1) is satisfied. (M3) also follows from \eqref{[Equation (1.5), HZZ19]} by modifying the constant $C_{L}$ in \eqref{[Equation (4.13), HZZ19]} to satisfy the upper bound of $C_{\tau_{L},q} (1+ \lVert u^{\text{in}} \rVert_{L_{x}^{2}}^{2q})$. Finally, in order to verify (M2), we let $s \leq t$ and $g$ be any bounded, $\mathbb{R}$-valued, $\mathcal{B}_{s}$-measurable and continuous function on $\Omega_{0}$. By Theorem \ref{[Theorem 1.1, HZZ19]}, we know that $u(\cdot \wedge T_{L})$ is $(\mathcal{F}_{t})_{t \geq 0}$-adapted. As $g$ is $\mathcal{B}_{s}$-measurable, this implies that $g(u(\cdot \wedge \tau_{L}(u)))$ is $\mathcal{F}_{s}$-measurable by \eqref{[Equation (3.15), HZZ19]}.  From the fact that $M_{t \wedge \tau_{L}(u), 0}^{u,i} \overset{\eqref{[Equation (3.12), HZZ19]}}{=} \langle B_{t \wedge \tau_{L} (u)}, e_{i} \rangle$ is an martingale w.r.t. $(\mathcal{F}_{t})_{t\geq 0}$ such that $\langle \langle M_{t \wedge \tau_{L} (u), 0}^{u,i} \rangle \rangle = \lVert G e_{i} \rVert_{L_{x}^{2}}^{2} (t \wedge \tau_{L} (u))$, it follows that $\mathbb{E}^{P} [ M_{t \wedge \tau_{L}, 0}^{i} g ] = \mathbb{E}^{P} [ M_{s \wedge \tau_{L}, 0}^{i} g]$, which implies that $M_{t \wedge \tau_{L}, 0}^{i}$ is a martingale w.r.t. $(\mathcal{B}_{t})_{t\geq 0}$ under $P$. Similarly, using the fact that $\langle \langle M_{t \wedge \tau_{L} (u), 0}^{u,i} \rangle \rangle = \lVert Ge_{i} \rVert_{L_{x}^{2}}^{2} (t \wedge \tau_{L} (u))$ implies that $(M_{t \wedge \tau_{L} (u), 0}^{u,i})^{2} - (t \wedge \tau_{L} (u)) \lVert G e_{i} \rVert_{L_{x}^{2}}^{2}$ is a martingale w.r.t. $(\mathcal{F}_{t})_{t \geq 0}$ under $\textbf{P}$, we can show that 
\begin{equation*}
 \mathbb{E}^{P} [ (( M_{t \wedge \tau_{L}, 0}^{i})^{2} - (t \wedge \tau_{L}) \lVert G e_{i} \rVert_{L_{x}^{2}}^{2}) g ] = \mathbb{E}^{P} [ ((M_{s \wedge \tau_{L}, 0}^{i})^{2} - (s \wedge \tau_{L}) \lVert G e_{i} \rVert_{L_{x}^{2}}^{2}) g] 
\end{equation*} 
which implies $\langle \langle M_{t \wedge \tau_{L}, 0}^{i} \rangle \rangle = (t \wedge \tau_{L}) \lVert G e_{i} \rVert_{L_{x}^{2}}^{2}$ and consequently $M_{t \wedge \tau_{L}, 0}^{i}$ is square-integrable. 

\subsection{Continuation of the poof of Proposition \ref{[Proposition 3.8, HZZ19]}}
First, it follows from \eqref{[Equation (3.16), HZZ19]} and \eqref{[Equation (3.15), HZZ19]} that there exists a $P$-measurable set $\mathcal{D} \subset \Omega_{0}$ such that $P(\mathcal{D}) = 0$ and for any $\omega \in \Omega_{0} \setminus \mathcal{D}$ and any $\delta \in (0, \frac{2}{5})$, 
\begin{equation}\label{[Equation (3.17), HZZ19]}
Z^{\omega} (\cdot \wedge \tau_{L} (\omega)) \in CH_{x}^{\frac{5+\sigma}{2}} \cap C_{\text{loc}}^{\frac{2}{5} - \delta}H_{x}^{\frac{3+ \sigma}{2}}.  
\end{equation}
On the other hand, we recall $M_{t,0}^{\omega}$ from \eqref{[Equation (3.12), HZZ19]} and define for every $\omega' \in \Omega_{0}, \omega \in \Omega_{0} \setminus \mathcal{D}$, 
\begin{align}
\mathbb{Z}_{\tau_{L} (\omega)}^{\omega'} (t) \triangleq& M_{t,0}^{\omega'} - e^{- (t- t \wedge \tau_{L}(\omega)) (-\Delta)^{m}} M_{t \wedge \tau_{L} (\omega), 0}^{\omega'}  - \int_{t \wedge \tau_{L} (\omega)}^{t} \mathbb{P} (-\Delta)^{m} e^{- (t-s) (-\Delta)^{m}} M_{s,0}^{\omega'} ds \nonumber \\
=& M_{t,0}^{\omega'} - M_{t \wedge \tau_{L} (\omega), 0}^{\omega'}  - \int_{t \wedge \tau_{L} (\omega)}^{t} \mathbb{P} (-\Delta)^{m} e^{-(t-s)(-\Delta)^{m}} [M_{s,0}^{\omega'} - M_{s\wedge \tau_{L} (\omega), 0}^{\omega'} ] ds \label{[Equation (3.17b), HZZ19]}
\end{align} 
and see that 
\begin{align}
Z^{\omega'} (t) - Z^{\omega'}(t \wedge \tau_{L}(\omega)) 
\overset{\eqref{[Equation (3.13), HZZ19]}}{=}&  \mathbb{Z}_{\tau_{L}(\omega)}^{\omega'} (t) + (e^{-(t - t \wedge \tau_{L} (\omega))(-\Delta)^{m}} - \text{Id}) Z^{\omega'} (t \wedge \tau_{L}(\omega)). \label{[Equation (3.17a), HZZ19]}
\end{align} 
From the proof of Proposition \ref{[Proposition 3.7, HZZ19]}, we know that $M_{t \wedge \tau_{L}, 0}^{i}$ is a martingale w.r.t. $(\mathcal{B}_{t})_{t\geq 0}$ under $P$ and hence by \eqref{[Equation (3.17b), HZZ19]} we deduce that $\mathbb{Z}_{\tau_{L} (\omega)}^{\omega'}$ is $\mathcal{B}^{\tau_{L} (\omega)}$-measurable. Next,  \eqref{[Equation (3.17a), HZZ19]} and \eqref{[Equation (3.4), HZZ19]} imply that $Q_{\omega}$ from Lemma \ref{[Proposition 3.2, HZZ19]} satisfies for any $\delta \in (0, \frac{2}{5})$ 
\begin{align}
&Q_{\omega} ( \{ \omega' \in \Omega_{0}: \hspace{0.5mm}  Z_{\cdot}^{\omega'} \in CH_{x}^{ \frac{5+ \sigma}{2}} \cap C_{\text{loc}}^{\frac{2}{5} - \delta} H_{x}^{\frac{3+ \sigma}{2}} \}) = \delta_{\omega} ( \{ \omega' \in \Omega_{0}: \hspace{0.5mm}  Z_{\cdot \wedge \tau_{L} (\omega)}^{\omega'} \in CH_{x}^{\frac{5+ \sigma}{2}} \cap C_{\text{loc}}^{\frac{2}{5} - \delta} H_{x}^{\frac{3+ \sigma}{2}} \}) \nonumber \\
& \hspace{34mm} \times R_{\tau_{L} (\omega), \xi(\tau_{L} (\omega), \omega)} ( \{ \omega' \in \Omega_{0}: \hspace{0.5mm}  \mathbb{Z}_{\tau_{L} (\omega)}^{\omega'} \in CH_{x}^{\frac{5+ \sigma}{2}} \cap C_{\text{loc}}^{\frac{2}{5} - \delta} H_{x}^{ \frac{3+ \sigma}{2}} \}). \label{[Equation (3.17c), HZZ19]}
\end{align} 
By \eqref{[Equation (3.17), HZZ19]}, the first quantity is one for all $\omega \in \Omega_{0} \setminus \mathcal{D}$.  On the other hand, we can write 
\begin{align}
 \int_{0}^{t} \mathbb{P} e^{-(t-s) (-\Delta)^{m}} d (M_{s,0}^{\omega'} - M_{s \wedge \tau_{L} (\omega), 0}^{\omega'})  
\overset{ \eqref{[Equation (3.12), HZZ19]}  \eqref{[Equation (3.17b), HZZ19]}}{=} \mathbb{Z}_{\tau_{L}(\omega)}^{\omega'} (t). \label{[Equation (3.17d), HZZ19]}
\end{align} 
As we deduced \eqref{[Equation (3.11a), HZZ19]} from \eqref{[Equation (3.10), HZZ19]}, the fact that the process $\omega' \mapsto M_{\cdot, 0}^{\omega'} - M_{\cdot \wedge \tau_{L} (\omega), 0}^{\omega'}$ is a $GG^{\ast}$-Wiener process starting from zero at  $\tau_{L} (\omega)$ w.r.t. $(\mathcal{B}_{t})_{t \geq 0}$ under $R_{\tau_{L} (\omega), \xi(\tau_{L} (\omega), \omega)}$ due to Lemma \ref{[Proposition 3.2, HZZ19]} and the hypothesis of $\text{Tr} ((-\Delta)^{ \frac{5}{2} - m + 2 \sigma} GG^{\ast}) < \infty$ allow us to deduce from \eqref{[Equation (3.17d), HZZ19]} that $
\mathbb{Z}_{\tau_{L} (\omega)}^{\omega'} \in CH_{x}^{ \frac{5+ \sigma}{2}} \cap C_{\text{loc}}^{\frac{2}{5} - \delta} H_{x}^{\frac{3+ \sigma}{2}}$ for $R_{\tau_{L} (\omega), \xi(\tau_{L} (\omega), \omega)}$-a.e. $\omega'$ and any $\delta \in (0, \frac{2}{5})$. Applying this to \eqref{[Equation (3.17c), HZZ19]} shows that for all $\omega \in \Omega_{0} \setminus \mathcal{D}$ there exists a measurable set $N_{\omega}$ such that $Q_{\omega} (N_{\omega}) = 0$ while the trajectory $t \mapsto Z^{\omega'}(t)$ lies in $CH_{x}^{\frac{5+ \sigma}{2}} \cap C_{\text{loc}}^{\frac{2}{5} - \delta} H_{x}^{\frac{3+ \sigma}{2}}$ for all $\omega' \in \Omega_{0} \setminus N_{\omega}$ for any $\delta \in (0, \frac{2}{5})$. We can also deduce 
\begin{equation}\label{estimate 3}
\tau_{L} (\omega') \overset{\eqref{[Equation (3.13b), HZZ19]}\eqref{[Equation (3.14), HZZ19]} }{=} \bar{\tau}_{L} (\omega') \hspace{1mm} \forall \hspace{1mm} \omega' \in \Omega_{0} \setminus N_{\omega}
\end{equation}
where 
\begin{align}
\bar{\tau}_{L} (\omega') \triangleq \inf \{ t \geq 0: \hspace{0.5mm}  \lVert Z^{\omega'} (t) \rVert_{H_{x}^{\frac{5+ \sigma}{2}}} \geq \frac{L^{\frac{1}{4}}}{C_{S}}  \wedge \inf\{t \geq 0: \hspace{0.5mm}  \lVert Z^{\omega'} \rVert_{C_{t}^{\frac{2}{5} - 2 \delta} H_{x}^{\frac{3+ \sigma}{2}}} \geq \frac{L^{\frac{1}{2}}}{C_{S}} \} \wedge L \label{[Equation (3.17g), HZZ19]}
\end{align} 
for $\delta \in (0, \frac{1}{30})$, which in turn implies that for $t < L$, 
\begin{align}
& \{ \omega' \in \Omega_{0}\setminus N_{\omega}: \hspace{0.5mm}  \tau_{L} (\omega') \leq t \} = \{ \omega' \in \Omega_{0}\setminus N_{\omega}: \hspace{0.5mm}  \sup_{s \in \mathbb{Q}: \hspace{0.5mm}  s \leq t} \lVert Z^{\omega'} (s) \rVert_{H_{x}^{\frac{5+ \sigma}{2}}} \geq \frac{L^{\frac{1}{4}}}{C_{S}} \} \nonumber\\
& \hspace{23mm} \cup \{ \omega' \in \Omega_{0}\setminus N_{\omega}: \hspace{0.5mm}  \sup_{s_{1}, s_{2} \in \mathbb{Q} \cap [0,t]: \hspace{0.5mm}  s_{1} \neq s_{2}} \frac{ \lVert Z^{\omega'} (s_{1}) - Z^{\omega'} (s_{2}) \rVert_{H_{x}^{\frac{3+ \sigma}{2}}}}{\lvert s_{1} - s_{2} \rvert^{\frac{2}{5} - 2 \delta}} \geq \frac{L^{\frac{1}{2}}}{C_{S}} \}. \label{[Equation (3.18), HZZ19]}
\end{align} 
It follows from \eqref{[Equation (3.18), HZZ19]} that $ \{ \omega' \in \Omega_{0}\setminus N_{\omega}: \hspace{0.5mm}  \tau_{L} (\omega') = \tau_{L}(\omega) \} \subset (\Omega_{0}\setminus N_{\omega}) \cap \mathcal{B}_{\tau_{L} (\omega)}^{0}$. This leads to for all $\omega \in \Omega_{0} \setminus \mathcal{D}$, 
\begin{align}\label{[Equation (3.19), HZZ19]}
Q_{\omega} ( \{\omega' \in \Omega_{0}: \hspace{0.5mm}  \tau_{L} (\omega') = \tau_{L}(\omega)\} ) = 1 
\end{align} 
and thus verifies \eqref{[Equation (3.5), HZZ19]} and completes the proof.

\section*{Acknowledgements} 

The author is grateful to Prof. Akif Ibraguimov for some comments on the work by Shnirelman \cite{S97}. The author also expresses deep gratitude to the Editor and the Referees for their valuable comments that have improved this manuscript significantly.

\end{document}